\documentclass[reqno,a4paper,11pt]{amsart}

\usepackage[utf8]{inputenc}
\usepackage[english]{babel}

\calclayout

\usepackage{verbatim}
\usepackage{float}
\usepackage{enumerate}
\usepackage{graphicx}
\usepackage{subcaption}

\usepackage{amsmath}
\usepackage{amsfonts}
\usepackage{amssymb}
\usepackage{amsthm}
\usepackage{mathtools}
\usepackage{bm}
\usepackage{esint}
\usepackage{upgreek}

\makeatletter
\def\thm@space@setup{%
  \thm@preskip=12pt plus 3pt minus 3pt
  \thm@postskip=\thm@preskip
}
\makeatother

\usepackage{tikz}
\usetikzlibrary{decorations.pathreplacing,decorations.markings,calc}
\usetikzlibrary{arrows}

\tikzset{
  on each segment/.style={
    decorate,
    decoration={
      show path construction,
      moveto code={},
      lineto code={
        \path [#1]
        (\tikzinputsegmentfirst) -- (\tikzinputsegmentlast);
      },
      curveto code={
        \path [#1] (\tikzinputsegmentfirst)
        .. controls
        (\tikzinputsegmentsupporta) and (\tikzinputsegmentsupportb)
        ..
        (\tikzinputsegmentlast);
      },
      closepath code={
        \path [#1]
        (\tikzinputsegmentfirst) -- (\tikzinputsegmentlast);
      },
    },
  },
  mid arrow/.style={postaction={decorate,decoration={
        markings,
        mark=at position .5 with {\arrow[#1]{stealth}}
      }}},
  lftarrow/.style={postaction={decorate,decoration={
        markings,
        mark=at position .25 with {\arrow[#1]{stealth}}
      }}},
  rgtarrow/.style={postaction={decorate,decoration={
        markings,
        mark=at position .75 with {\arrow[#1]{stealth}}
      }}},
  rmid arrow/.style={postaction={decorate,decoration={
        markings,
        mark=at position .3 with {\arrowreversed[#1]{stealth}}
      }}},
  vmid arrow/.style 2 args={postaction={decorate,decoration={
        markings,
        mark=at position #2 with {\arrow[#1]{stealth}}
      }}},
  vrmid arrow/.style 2 args={postaction={decorate,decoration={
        markings,
        mark=at position #2 with {\arrowreversed[#1]{stealth}}
      }}},
}

\makeatletter
\def\grd@save@target#1{\def\grd@target{#1}}
\def\grd@save@start#1{\def\grd@start{#1}}
\tikzset{
  grid with coordinates/.style={
    to path={%
      \pgfextra{%
        \edef\grd@@target{(\tikztotarget)}%
        \tikz@scan@one@point\grd@save@target\grd@@target\relax
        \edef\grd@@start{(\tikztostart)}%
        \tikz@scan@one@point\grd@save@start\grd@@start\relax
        \draw[minor help lines] (\tikztostart) grid (\tikztotarget);
        \draw[major help lines] (\tikztostart) grid (\tikztotarget);
        \grd@start
        \pgfmathsetmacro{\grd@xa}{\the\pgf@x/1cm}
        \pgfmathsetmacro{\grd@ya}{\the\pgf@y/1cm}
        \grd@target
        \pgfmathsetmacro{\grd@xb}{\the\pgf@x/1cm}
        \pgfmathsetmacro{\grd@yb}{\the\pgf@y/1cm}
        \pgfmathsetmacro{\grd@xc}{\grd@xa + \pgfkeysvalueof{/tikz/grid with coordinates/major step}}
        \pgfmathsetmacro{\grd@yc}{\grd@ya + \pgfkeysvalueof{/tikz/grid with coordinates/major step}}
        \foreach \x in {\grd@xa,\grd@xc,...,\grd@xb}
          \node[anchor=north] at (\x,\grd@ya) {\pgfmathprintnumber{\x}};
        \foreach \y in {\grd@ya,\grd@yc,...,\grd@yb}
          \node[anchor=east] at (\grd@xa,\y) {\pgfmathprintnumber{\y}};
      }
    }
  },
  minor help lines/.style={
    help lines,
    step=\pgfkeysvalueof{/tikz/grid with coordinates/minor step}
  },
  major help lines/.style={
    help lines,
    line width=\pgfkeysvalueof{/tikz/grid with coordinates/major line width},
    step=\pgfkeysvalueof{/tikz/grid with coordinates/major step}
  },
  grid with coordinates/.cd,
  minor step/.initial=.2,
  major step/.initial=1,
  major line width/.initial=2pt,
}
\makeatother

\usepackage{xcolor}
\definecolor{lightgreen}{RGB}{100,200,100}
\definecolor{apricot}{rgb}{0.98,0.81,0.69}

\usepackage[backgroundcolor=lightgreen]{todonotes}

\usepackage{seqsplit}
\usepackage{needspace}

\usepackage[pdftex,colorlinks]{hyperref}
\hypersetup{
  colorlinks=true,
  linkcolor=blue,
  citecolor=red
}

\usepackage[
backend=biber,
style=numeric-comp,
sorting=nyt,
giveninits=true,
maxbibnames=99
]{biblatex}
\usepackage{csquotes}
\AtEveryBibitem{%
  \clearfield{issn}%
  \clearfield{isbn}%
  \clearlist{language}%
  \iffieldundef{doi}{}{\clearfield{url}}%
}
\renewbibmacro{in:}{%
  \ifentrytype{article}{}{\printtext{\bibstring{in}\intitlepunct}}}
\newcommand*{\doihrefaux}[2]{\href{https://doi.org/#1}{\nolinkurl{#2}}}
\ExplSyntaxOn
\cs_new_protected:Npn \doihref #1
  {
    \str_set:Nn \l_tmpa_str {#1}
    \str_set:Nx \l_tmpb_str { \c_percent_str 23 }
    \exp_args:NNVV \str_replace_all:Nnn \l_tmpa_str \c_hash_str \l_tmpb_str
    \exp_args:NV \doihrefaux \l_tmpa_str {#1}
  }
\ExplSyntaxOff
\DeclareFieldFormat{doi}{%
  \mkbibacro{DOI}\addcolon\space
  \ifhyperref{\doihref{#1}}{\nolinkurl{#1}}}

\makeatletter
\def\l@subsection{\@tocline{2}{0pt}{2.5pc}{5pc}{}}
\makeatother

\DeclareMathOperator{\diag}{diag}

\DeclareMathOperator{\ee}{e}
\DeclareMathOperator{\supp}{supp}

\DeclareMathOperator{\sign}{sign}

\DeclareMathOperator{\I}{I}
\DeclareMathOperator{\CT}{C}
\DeclareMathOperator{\U}{U}

\renewcommand{\Re}{\mathop{\rm Re}}
\renewcommand{\Im}{\mathop{\rm Im}}

\newcommand{\C}{\mathbb{C}}
\newcommand{\R}{\mathbb{R}}
\newcommand{\Z}{\mathbb{Z}}

\newcommand{\boh}{\mathit{o}}
\newcommand{\Boh}{\mathcal{O}}

\newcommand{\ii}{\mathrm{i}}
\newcommand{\dd}{\mathrm{d}}

\renewcommand{\bm}{\mathbf}
\newcommand{\mcal}{\mathcal}
\newcommand{\SW}{\mathcal{W}}

\newcommand{\msf}{\mathsf}

\newcommand{\wh}{\widehat}
\newcommand{\wt}{\widetilde}

\renewcommand{\sp}{\boldsymbol\sigma}

\newcommand*{\defeq}{\mathrel{\vcenter{\baselineskip0.5ex \lineskiplimit0pt
  \hbox{\scriptsize.}\hbox{\scriptsize.}}}=}

\newcommand\restr[2]{%
  \left.\kern-\nulldelimiterspace
  #1\vphantom{\big|}
  \right|_{#2}
}

\newtheorem{theorem}{Theorem}[section]

\newtheorem{proposition}[theorem]{Proposition}
\newtheorem{lemma}[theorem]{Lemma}
\newtheorem{corollary}[theorem]{Corollary}

\theoremstyle{definition}
\newtheorem{rhp}[theorem]{RHP}
\newtheorem{definition}[theorem]{Definition}

\newtheorem{assumption}[theorem]{Assumption}

\theoremstyle{remark}
\newtheorem{remark}[theorem]{Remark}

\numberwithin{equation}{section}

\begin{document}

\title[MVOPs with non-scalar varying exponential weights]{Asymptotics of matrix orthogonal polynomials with non-scalar varying exponential weights}

\author[V.~Alves]{Victor Alves}
\address{School of Mathematical Sciences, Fudan University, 220 Handan Road, Shanghai 200433, China}
\email{victoralves@fudan.edu.cn}

\date{}



\begin{abstract}
    We study the $2\times2$ matrix orthogonal polynomials $\bm P_N$ on the real line with respect to the varying weights $\ee^{-N\bm Q(z)}\bm B(z)$, where $\bm Q$ and $\bm B$ are real symmetric matrix polynomials that commute with a Hermite-type matrix weight, $\bm Q$ is not scalar and $\bm B$ is positive definite. We are motivated by the non-Abelian Toda deformation of the Hermite-type weight, in the sense of Dea\~no, Morey and Rom\'an. The weight depends on $N$ through a matrix exponential that is not scalar, so that its two eigenvalues have different exponential rates. The matrix orthogonality reduces to two scalar orthogonality problems on the real line, with a varying weight that has an essential singularity at the origin. Assuming that the associated equilibrium measure is one-cut regular, we obtain the strong asymptotics of $\bm P_N$ and of $\det \bm P_N$ by the Deift--Zhou steepest descent method. Our main result is the limiting zero distribution of the entries of $\bm P_N$. These zeros are complex in general, and they accumulate on a curve where the contributions of the two eigenvalues of the weight have the same size. When this curve crosses the real segment that attracts the zeros of $\det\bm P_N$, the zeros of the entries display a critical behavior, and a positive fraction of them accumulates on this segment.
\end{abstract}

\vspace*{-1.6cm}

\maketitle

\section{Introduction}\label{sec:intro}

\subsection{Background} \hfill

    Orthogonal polynomials on the real line with respect to varying weights $\ee^{-NV(x)}$ are by now well understood as $N \to \infty$. The classical example is given by the Hermite polynomials, for which the asymptotics of Plancherel and Rotach \cite{PlancherelRotach29} describe the behavior in the oscillatory region, near its endpoints and outside it. For a general real analytic $V$, the strong asymptotics were obtained by Deift, Kriecherbauer, McLaughlin, Venakides and Zhou \cite{Egregium,deift_exp_weights}, by a combination of the Riemann--Hilbert characterization of orthogonal polynomials due to Fokas, Its and Kitaev \cite{FokasItsKitaev90,FokasItsKitaev91} with the Deift--Zhou steepest descent method, see also \cite{Deift1999}. The central object in this analysis is the equilibrium measure of an external field $V$ from logarithmic potential theory \cite{DKM98,SaffTotik2024}, which gives the limiting distribution of the zeros.

    Matrix orthogonal polynomials (MVOPs) go back to Krein \cite{Krein49}. Given a $d\times d$ matrix weight $\bm W$ on $\R$, positive definite and with finite moments, the monic MVOPs $\bm P_N$ are the matrix polynomials of degree $N$ with leading coefficient $\bm I$ such that
    \[
        \int_\R \bm P_N(z)\,\bm W(z)\,z^k\,\dd z = \bm 0, \qquad
        k = 0, 1, \ldots, N-1.
    \]
    Many features of the scalar theory find analogues in the matrix setting, such as the three-term recurrence relation, now with matrix coefficients, and the Christoffel--Darboux formula, see \cite{Berg,DPS,Miranian}.

    Many examples come from the matrix analogue of a classical problem of Bochner. Bochner \cite{Bochner29} showed that, up to affine changes of variable, the Hermite, Laguerre and Jacobi polynomials are the only polynomials orthogonal with respect to a positive weight on the real line that are eigenfunctions of a second order differential operator. Dur\'an \cite{Duran97} investigated the matrix version of this problem, that is, the matrix weights whose MVOPs are eigenfunctions of a second order differential operator with matrix coefficients. Such a weight is of interest only if it does not reduce to scalar weights, that is, if it is not of the form $\bm M\bm D\bm M^\ast$ with $\bm M$ constant and $\bm D$ diagonal, since the MVOPs of $\bm M\bm D\bm M^\ast$ are $\bm M\bm D_N\bm M^{-1}$, with $\bm D_N$ a diagonal matrix of scalar orthogonal polynomials.

    The first examples that do not reduce to scalar weights came from the matrix spherical functions associated to the complex projective plane \cite{GPT02} and from a general construction of Dur\'an and Gr\"unbaum \cite{DG04,DG05}, which produced, among others, the Hermite-type weights $\ee^{-z^2}\ee^{z\bm A}\ee^{z\bm A^\ast}$ with $\bm A$ nilpotent. Their MVOPs, the matrix Hermite polynomials, have Rodrigues formulas and explicit expressions \cite{DG05,IKR2}, and they are related to Burchnall-type identities, the non-Abelian Toda lattice \cite{IKR2} and non-commutative Painlev\'e equations \cite{CafassoI}. Ladder relations for the more general weights $\ee^{-v(z)}\ee^{z\bm A}\ee^{z\bm A^\ast}$, with $v$ a polynomial of even degree, were studied in \cite{DER}. Analogues of Laguerre and Gegenbauer type are also known \cite{KoelinkRlaguerre,KdlRR}, and Casper and Yakimov \cite{Casper2} solved the matrix Bochner problem under natural assumptions.

    Much less is known about the behavior of MVOPs as the degree tends to infinity. The Riemann--Hilbert characterization extends to MVOPs \cite{GdIMF11}, but the Riemann--Hilbert problem has size $2d\times2d$, and its steepest descent analysis requires additional structure. The cases worked out so far rely on a factorization of the weight. Dea\~no, Kuijlaars and Rom\'an \cite{Deano2023a} studied MVOPs on $[-1,1]$, and Dea\~no and Rom\'an \cite{DeanoRoman2025} obtained the strong asymptotics of MVOPs on the real line for varying weights of the form $\ee^{-Nv(z)}\bm M(z)$, with a scalar potential $v$ whose equilibrium measure is supported on one interval and a matrix part $\bm M$ that does not depend on $N$. Both works use a matrix Szeg\H{o} function, and Lahiry \cite{Lahiry26} extended the approach of \cite{DeanoRoman2025} to allow equilibrium measures with disconnected support. In these works the weight either does not depend on $N$, or depends on $N$ only through the scalar factor $\ee^{-Nv}$.

    For MVOPs arising from tiling models \cite{Duits2021,Groot2021a}, the weight is a power of a fixed rational matrix-valued function, with exponent proportional to the degree, and the orthogonality is non-Hermitian on a closed contour. The eigenvalues of the weight are naturally functions on a Riemann surface. Charlier \cite{Charlier2021} showed that, for rational weights that are diagonalizable, the matrix Christoffel--Darboux kernel is equivalent to a scalar reproducing kernel on such a Riemann surface. Groot and Kuijlaars \cite{Groot2021a} reduced the $2\times2$ MVOPs of a hexagon tiling model to (scalar) orthogonal polynomials. From a steepest descent analysis of these polynomials they obtained the limiting zero distribution of $\det\bm P_n$ and, under a geometric condition supported by numerical evidence, that of the $(1,2)$ entry, which involves a balayage measure.

\subsection{Overview} \hfill

    In this paper we consider $2\times2$ varying weights $\ee^{-N\bm Q}\bm B$ on the real line, where $\bm Q$ and $\bm B$ are real symmetric matrix polynomials, $\bm Q$ is not scalar, and $\bm B$ is positive definite and does not depend on $N$. Our motivation is a non-Abelian Toda deformation, in the sense of Dea\~no, Morey and Rom\'an \cite{DMR}, which deforms the Hermite-type weight $\bm H_{\msf a}$ into $\ee^{-\msf t\bm\Lambda}\bm H_{\msf a}$ and takes this form when $\msf t$ is proportional to $N$ (Section~\ref{sec:intro-toda}). To our knowledge, this is the first time that the asymptotics of MVOPs on the real line are obtained for varying matrix exponential weights whose eigenvalues have different exponential rates.

    As in the Toda example, we assume that $\bm Q$ and $\bm B$ commute with $\bm H_{\msf a}$. Up to a constant rotation and an affine change of variable, this assumption covers all the weights in which $\bm Q$ and $\bm B$ commute with a real symmetric matrix polynomial of degree one with non-constant eigenvectors. Moreover, the approach should extend to other weights whose common eigenvectors are parametrized by a curve of genus zero, possibly with additional orthogonality conditions (see Remark~\ref{rem:other-weights}). As in \cite{Groot2021a}, the genus zero parametrization allows us to reduce the matrix orthogonality to scalar orthogonality on the real line (Section~\ref{sec:intro-setup}), which we analyze by the Deift--Zhou steepest descent method.

    Since $\ee^{-N\bm Q}$ is not scalar, its two eigenvalues have different exponential rates, and each entry of $\bm P_N$ combines two contributions of different exponential sizes. Our main results are the following.
    \begin{itemize}
        \item We obtain the limiting zero distribution of each entry of $\bm P_N$ (Theorem~\ref{thm:intro-zeros-general}), assuming that the equilibrium measure associated to the external field $V$ in \eqref{eq:intro-V} is one-cut regular. In contrast with \cite[Theorem~2.5]{Groot2021a}, no extra condition is needed. The zeros are complex in general, and they accumulate on the curve where the two contributions have the same size and, possibly, on part of the real segment that attracts the zeros of $\det\bm P_N$.
        \item The zeros of the entries display a critical behavior: when this curve crosses the real segment, a positive fraction of them accumulates on the segment, and otherwise all of them accumulate on the curve. Figure~\ref{fig:zeros-entries} illustrates both situations for the Toda deformed weight (Remark~\ref{rem:critical}).
        \item Away from the zeros, each entry of $\bm P_N$ is dominated by one of the two contributions, and the leading term of $\bm P_N$ is a rank-one matrix (Theorem~\ref{thm:intro-asymptotics}), whereas $\det\bm P_N$ behaves like a scalar orthogonal polynomial of degree $2N$ with real zeros (Theorem~\ref{thm:det-asymptotics}).
    \end{itemize}
    We also show that the Toda deformed weight satisfies our assumption for large values of the deformation parameter (Theorem~\ref{thm:intro-eq}).

\subsection{Decomposition on the spectral curve}\label{sec:intro-setup} \hfill

    For a fixed $\msf a > 0$, the $2\times2$ Hermite-type weight of Dur\'an and Gr\"unbaum \cite{DG04,DG05}, see also \cite{IKR2}, is
    \begin{equation}\label{eq:intro-A-Lambda}
        \bm H_{\msf a}(z)
            \defeq \ee^{-z^2} \begin{pmatrix}
            1   &   \msf a z    \\
            \msf a z    &   1 + \msf a^2 z^2
        \end{pmatrix}
            = \ee^{-z^2} \left[(1 + \msf a^2 z^2) \bm I  +\msf a z \bm \Lambda(z) \right], \qquad \bm \Lambda(z) = \begin{pmatrix}
                - \msf a z  &   1   \\
                1   &   0
            \end{pmatrix},
    \end{equation}
    which is of the form $\ee^{-z^2}\ee^{z\bm A}\ee^{z\bm A^\ast}$ with $\bm A = \left(\begin{smallmatrix} 0 & 0 \\ \msf a & 0 \end{smallmatrix}\right)$. A matrix polynomial commutes with $\bm H_{\msf a}$ if, and only if, it commutes with $\bm \Lambda$, and the real symmetric $2\times2$ matrix polynomials with this property are
    \begin{equation}\label{eq:intro-Q}
        \bm Q(z) = Q_0(z)\,\bm I + Q_1(z)\,\bm\Lambda(z), \qquad Q_0, Q_1 \in \R[z].
    \end{equation}
    Given two such matrix polynomials $\bm Q$ and $\bm B$, with $\bm B$ positive definite on $\R$, we consider the matrix weight
    \begin{equation}\label{eq:intro-weight-Q}
        \bm W_{\bm Q}(z) \defeq \ee^{-\bm Q(z)}\,\bm B(z), \qquad z \in \R,
    \end{equation}
    and its monic MVOPs $\bm P_N$.

    The common eigenvectors of $\bm Q$ and $\bm B$ are those of $\bm\Lambda$, namely $(1,x)^T$, where $x$ solves $x^2 - \msf a zx - 1 = 0$. This equation defines a spectral curve $\mcal R$ of genus zero, a double cover of $\C$ with global coordinate $x \in \C\setminus\{0\}$. Over every point $z \in \C$ lie two points $x$ and $\iota x = \iota(x)$, related by
    \begin{equation}\label{eq:intro-z-iota}
        z = \msf z(x) \defeq \frac{1}{\msf a}\left(x - \frac1x\right), \qquad \iota(x) \defeq -\frac1x.
    \end{equation}
    The involution $\iota$ exchanges the two points $x$ and $\iota x$ of $\mcal R$ lying over $z$.
    
    Writing $\bm B = B_0\,\bm I + B_1\,\bm\Lambda$, the eigenvalues of $\bm Q(\msf z(x))$ and $\bm B(\msf z(x))$ for the eigenvector $(1,x)^T$ are the real Laurent polynomials
    \begin{equation}\label{eq:intro-Qx}
        \msf Q(x) \defeq Q_0(\msf z(x)) + \frac{Q_1(\msf z(x))}{x}, \qquad \msf B(x) \defeq B_0(\msf z(x)) + \frac{B_1(\msf z(x))}{x},
    \end{equation}
    and $\bm B$ is positive definite on $\R$ if, and only if, $\msf B > 0$ on $\R\setminus\{0\}$.

\begin{figure}[ht]
        \centering
        \captionsetup{font=small,margin=0pt}
        \includegraphics{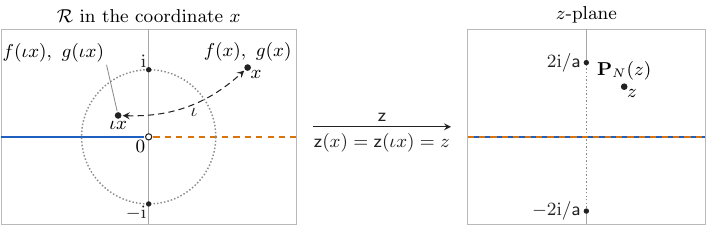}
        \caption{The spectral curve $\mcal R$ in the coordinate $x$ (left) and the $z$-plane (right). The functions $f$ and $g$ are defined on $\mcal R$, while $\bm P_N$ is defined on the $z$-plane. The projection $\msf z$ is two-to-one: the points $x$ and $\iota x = -1/x$ lie over the same point $z$, and by \eqref{eq:intro-two-sheets} each entry of $\bm P_N(z)$ combines the values of $f$, or of $g$, at $x$ and $\iota x$. The involution $\iota$ fixes $\pm\ii$, which lie over the branch points $\pm 2\ii/\msf a$, and exchanges the interior and the exterior of the unit circle (dotted), which $\msf z$ maps onto the segment $[-2\ii/\msf a, 2\ii/\msf a]$. The negative half-line (solid) and the positive half-line (dashed) are both mapped onto the real line of the $z$-plane, which is therefore covered twice. The point $x = 0$ (white dot) does not belong to $\mcal R$.}
        \label{fig:spectral-curve}
    \end{figure}

    Putting $\bm P_N = (\bm P_{ij})_{i,j=1}^2$ and pairing its rows with the eigenvector $(1,x)^T$, we consider two Laurent polynomials
    \[
        f(x) = f_N(x) \defeq \bm P_{11}(\msf z(x)) + x\,\bm P_{12}(\msf z(x)),
        \qquad
        g(x) = g_N(x) \defeq \bm P_{21}(\msf z(x)) + x\,\bm P_{22}(\msf z(x)),
    \]
    which determine $\bm P_N$ (Lemma~\ref{lem:divided-differences}),
    \begin{equation}\label{eq:intro-two-sheets}
        \bm P_N(\msf z(x)) = \frac{1}{x - \iota x}
        \begin{pmatrix}
            x f(\iota x) -  f(x) \iota x &   f(x) - f(\iota x) \\
            x g(\iota x) -  g(x) \iota x &   g(x) - g(\iota x)
        \end{pmatrix},
        \qquad x \neq 0, \pm\ii,
    \end{equation}
    so that each entry of $\bm P_N$ combines the values of $f$, or of $g$, at the two points of $\mcal R$ over $z$ (see Figure~\ref{fig:spectral-curve}). Moreover, the orthogonality of $\bm P_N$ is equivalent to the orthogonality of $f$ and $g$ on the real line of the $x$-plane (Theorem~\ref{thm:orth-fg}):
    \[
        \int_\R f(x)\,x^k\,w_{\bm Q}(x)\,\dd x = \int_\R g(x)\,x^k\,w_{\bm Q}(x)\,\dd x = 0, \quad -N + 1 \leq k \leq N, \qquad w_{\bm Q}(x) \defeq \lambda(x)\,\ee^{-\msf Q(x)},
    \]
    where $\lambda(x) \defeq x^{-2}\msf B(x)$. This decomposition is the analogue, for the weights considered here, of the correspondence between matrix orthogonality and scalar orthogonality on a Riemann surface established by Charlier \cite{Charlier2021} for rational diagonalizable weights with non-Hermitian orthogonality, and of the reduction of Groot and Kuijlaars \cite{Groot2021a} for a hexagon tiling model. It makes the entries of $\bm P_N$ accessible to the Deift--Zhou steepest descent method, since $f$ and $g$ are orthogonal Laurent polynomials.

\subsection{Main results} \hfill

    From now on we consider the varying weight $\bm W_{N\bm Q} = \ee^{-N\bm Q}\bm B$, and $\bm P_N$, $f$ and $g$ refer to this weight. The scalar weight becomes $\lambda(x)\,\ee^{-N\msf Q(x)}$, and the numerator $p(x) = x^Nf(x)$, a polynomial of degree at most $2N$, is orthogonal to $x^m$, $1 \leq m \leq 2N$, with respect to
    \begin{equation}\label{eq:intro-scalar-weight}
        x^{-2N}\lambda(x)\,\ee^{-N\msf Q(x)} = \lambda(x)\,\ee^{-2NV(x)},
    \end{equation}
    a varying weight with external field
    \begin{equation}\label{eq:intro-V}
        V(x) \defeq \frac12\,\msf Q(x) + \log|x|, \quad x \in \R\setminus\{0\}, \qquad V(0) \defeq +\infty.
    \end{equation}
    The factor $\lambda$ does not depend on $N$ and acts as a Szeg\H{o} factor. The field $V$ is semiclassical (i.e., $V'$ is rational) with only one pole, at $x = 0$, so that the weight $\ee^{-2NV}$ has an essential singularity there. Weights with an essential singularity at the origin appear in several problems of random matrix theory \cite{ChenIts10,MezzadriMo09,XuDaiZhao14,BrightmoreMezzadriMo15,AtkinClaeysMezzadri16}, where the main interest lies in the interaction of the singularity with the eigenvalues near the origin, described by Painlev\'e III equations in double scaling limits. In our setting the singular part of $V$ is of the same order as the rest of the field, and the origin lies at a positive distance from the support of the equilibrium measure. This is also the setting of McLaughlin, Vartanian and Zhou \cite{MVZ06,MVZ08,MVZ08b}, who obtained the strong asymptotics of Laurent polynomials orthogonal with respect to varying exponential weights.

    The asymptotic behavior of $f$ and $g$ is governed by the equilibrium measure $\mu$ of $V$, the unique minimizer of the energy functional
    \[
        \I^V(\sigma) \defeq \iint\log\frac{1}{|x - y|}\,\dd\sigma(x)\,\dd\sigma(y) + \int V\,\dd\sigma
    \]
    among the probability measures on $\R$, where we assume that $V$ is \textit{admissible}, that is, that the weight $\ee^{-V/2}$ is admissible in the sense of \cite[Definition~I.1.3]{SaffTotik2024}. Then $\mu$ has compact support, $0 \notin \supp \mu$, and $\mu$ is characterized by the Euler--Lagrange conditions
    \[
        \begin{split}
            2\int \log |x-y|^{-1} \, \dd \mu(y) + V(x) &= \ell, \quad x \in \supp \mu, \\
        2 \int \log |x-y|^{-1} \, \dd \mu(y) + V(x) &\geq \ell, \quad x \in \R \setminus \supp \mu,
        \end{split}
    \]
    which in general hold quasi-everywhere, but here hold everywhere, since the density of $\mu$ is bounded (Lemma~\ref{lem:quad}) and $V$ is continuous on $\R\setminus\{0\}$. For real analytic external fields the support of the equilibrium measure consists of finitely many intervals \cite{DKM98}, and the same holds for $V$ (Lemma~\ref{lem:quad}). We work under the following assumption.

    \begin{assumption}\label{ass:onecut}
        The equilibrium measure $\mu$ is \textit{one-cut regular}, with $\supp \mu = [a_1,b_1] \subset (-\infty, 0)$. That is, the Euler--Lagrange inequality is strict outside $\supp \mu$, and the density of $\mu$ is positive in $(a_1,b_1)$ and vanishes as a square root at $a_1$ and $b_1$.
    \end{assumption}

    \begin{remark}\label{rem:scalar-Q}
        The one-cut condition forces $\bm Q$ to be non-scalar. Indeed, if $Q_1 \equiv 0$ in \eqref{eq:intro-Q}, then $\msf Q = Q_0\circ\msf z = \msf Q\circ\iota$, so that $V(\iota x) + 2\log|x| = V(x)$ and $\I^V(\iota_*\sigma) = \I^V(\sigma)$ for every probability measure $\sigma$ with compact support in $\R\setminus\{0\}$, where $\iota_*\sigma$ is the push-forward of $\sigma$ by $\iota$. By uniqueness, $\mu$ is invariant under $\iota$, which exchanges the two half-lines, so that $\supp\mu$, which does not contain $0$, meets both of them and is not one-cut. This applies to the classical Toda deformation $\ee^{\msf a\msf t z}\bm H_{\msf a}(z)$, for which $\bm Q(z) = (z^2 - \msf a\msf t z)\,\bm I$. Weights with a scalar potential are treated in \cite{DeanoRoman2025,Lahiry26}.
    \end{remark}

    \begin{remark}\label{rem:symmetries}
        Under the one-cut assumption there is no loss of generality in assuming $\supp \mu = [a_1,b_1] \subset (-\infty, 0)$. Indeed, putting $\sp_3 = \diag (1,-1)$ we have $\sp_3\bm\Lambda(-z)\sp_3 = -\bm\Lambda(z)$, and thus the weight $\check{\bm W}(z) = \sp_3\bm W_{N\bm Q}(-z)\sp_3$ is also of the form \eqref{eq:intro-weight-Q}. Its associated MVOPs are $(-1)^N\sp_3\bm P_N(-z)\sp_3$ and its associated scalar external field is $V(-x)$, so that the equilibrium measure associated to $\check{\bm W}$ is the push-forward of $\mu$ by $x \mapsto -x$.
    \end{remark}

    We start with the asymptotics of $f = f_N$, $g = g_N$ and $\det \bm P_N$, which are proven in Sections~\ref{sec:asymptotics-pq} and~\ref{sec:det}.

    \begin{theorem}[Strong asymptotics for $f$ and $g$]\label{thm:intro-asymptotics-fg}
        There are functions $A_f$ and $A_g$, analytic and non-vanishing in $\C\setminus[a_1,b_1]$, with finite non-zero limits at infinity, given explicitly in \eqref{eq:amplitudes}, such that for every compact set $K$ of $\C \setminus \left( [a_1,b_1] \cup \{0\}\right)$
        \begin{align*}
            f_N(x) &= \left(-\frac{1}{\msf a}\right)^N\frac{\ee^{2N(\msf g(x) - \msf g(0))}}{x^N}\,A_f(x)\left(1 + \Boh\left(\frac1N\right)\right), \\
            g_N(x) &= \msf a^{-N}\,\frac{\ee^{2N\msf g(x)}}{x^{N-1}}\,A_g(x)\left(1 + \Boh\left(\frac1N\right)\right),
        \end{align*}
        uniformly for $x \in K$ as $N \to +\infty$, where $\msf g(x) \defeq \int \log (x-s) \, \dd \mu(s)$.
    \end{theorem}

    \begin{theorem}[Strong asymptotics for $\det \bm P_N$] \label{thm:det-asymptotics}
        Let $\msf z_*\mu$ be the push-forward of $\mu$ by $\msf z$, a probability measure on the segment $\msf z([a_1,b_1]) \subset \R$. There is a function $\wh A_{\det}$, analytic and non-vanishing in $\overline{\C}\setminus\msf z([a_1,b_1])$, positive on $\R\setminus\msf z([a_1,b_1])$, with $\wh A_{\det}(\infty) = 1$, given explicitly in \eqref{eq:A-det-hat}, such that for every neighborhood $O$ of $\msf z([a_1,b_1])$
        \begin{equation}\label{eq:det-asymptotics}
            \det\bm P_N(\zeta) = \ee^{2N\wh{\msf g}(\zeta)}\,\wh A_{\det}(\zeta)\left(1 + \Boh\left(\frac1N\right)\right), \quad \wh{\msf g}(\zeta) \defeq \int \log(\zeta - t)\,\dd(\msf z_*\mu)(t),
        \end{equation}
        uniformly for $\zeta \in \C\setminus  O$ as $N \to +\infty$. Moreover, the normalized zero counting measure of $\det \bm P_N$ converges weakly to $\msf z_*\mu$ as $N \to \infty$.
    \end{theorem}

    Theorem~\ref{thm:det-asymptotics} says that asymptotically $\det\bm P_N$ behaves like an orthogonal polynomial of degree $2N$ whose zeros are distributed according to $\msf z_*\mu$. The factor $\bm B$ enters Theorems~\ref{thm:intro-asymptotics-fg} and~\ref{thm:det-asymptotics} only through the functions $A_f$, $A_g$ and $\wh A_{\det}$.

    On each entry of $\bm P_N$, the two terms evaluated at $x$ and $\iota x$ in \eqref{eq:intro-two-sheets} are the contributions of the two points of the spectral curve over $\zeta = \msf z(x)$. Their relative size is measured by the function
    \begin{equation}\label{eq:intro-W}
        \begin{split}
            \SW(x) \defeq \int\log\left|\frac{1 + sx}{x - s}\right|\dd\mu(s).
        \end{split}
    \end{equation}
    Indeed, by Theorem~\ref{thm:intro-asymptotics-fg}, $|f(\iota x)/f(x)| = \ee^{2N\SW(x) + \Boh(1)}$ for every $x \neq 0$ such that $x$ and $\iota x$ are not in $[a_1,b_1]$, and the same holds for $g$. Hence $f(x)$ dominates $f(\iota x)$ exponentially where $\SW(x) < 0$, the opposite happens where $\SW(x) > 0$, and both contributions are comparable only near
    \[
        \Gamma_{\SW} \defeq \{ x \in \C : \SW(x) = 0 \}.
    \]
    We therefore expect the zeros of the entries of $\bm P_N(\zeta)$ outside $\msf z([a_1,b_1])$ to accumulate on
    \[
        \Gamma \defeq \msf z (\Gamma_{\SW}\setminus\{0\}).
    \]
    For $\zeta \notin \Gamma$, the point of $\mcal R$ over $\zeta$ at which $\SW < 0$ dominates in \eqref{eq:intro-two-sheets}, except possibly over $\msf z([a_1,b_1]\cap\{\SW < 0\})$, where this point lies on the support of $\mu$ and $f$ and $g$ oscillate. The next theorem makes this precise.

    \begin{theorem}[Strong asymptotics for $\bm P_N$]\label{thm:intro-asymptotics}
        For $\zeta \in \C\setminus\Gamma$ let $x = x_-(\zeta)$ be the solution of $\msf z(x) = \zeta$ with $\SW(x) < 0$. Let $K$ be a compact subset of $\C\setminus\left(\Gamma\cup\msf z([a_1,b_1]\cap\{\SW < 0\})\right)$. Then there exists $c > 0$ such that
        \[
            \bm P_N(\zeta) = \frac{1}{x - \iota x}
            \begin{pmatrix}
                -\iota x\,f(x)\left(1 + \Boh\left(\ee^{-cN}\right)\right) & f(x)\left(1 + \Boh\left(\ee^{-cN}\right)\right) \\[1ex]
                -\iota x\,g(x)\left(1 + \Boh\left(\ee^{-cN}\right)\right) & g(x)\left(1 + \Boh\left(\ee^{-cN}\right)\right)
            \end{pmatrix},\quad  x = x_-(\zeta), \quad N \to +\infty,
        \]
        uniformly for $\zeta \in K$. In particular, the leading term of $\bm P_N(\zeta)$ is the rank-one matrix
        \[
            \frac{1}{x - \iota x}\left(f(x),\,g(x)\right)^T\left(-\iota x,\,1\right).
        \]
    \end{theorem}

    Theorem~\ref{thm:intro-asymptotics} is proven in Section~\ref{sec:proof-asymptotics}. In particular, for $N$ large no entry of $\bm P_N$ vanishes on $K$, so that the zeros of the entries can only accumulate on $\Gamma\cup\msf z([a_1,b_1]\cap\{\SW < 0\})$. In contrast with the zeros of $\det\bm P_N$, they are in general complex, and the next theorem gives their limiting distribution.

    \begin{theorem}\label{thm:intro-zeros-general}
        Let $\nu$ be the positive measure on $\Gamma_{\SW}$ with density $\frac{1}{2\pi}|\nabla\SW|$ with respect to arc length. Then, for each entry of $\bm P_N$, the normalized zero counting measure converges weakly to the probability measure
        \begin{equation}\label{eq:intro-zeros-general}
            \rho \defeq \msf z_*\nu + 2\,\msf z_*\big(\restr{\mu}{\{\SW < 0\}}\big).
        \end{equation}
    \end{theorem}

    Theorem~\ref{thm:intro-zeros-general} is proven in Section~\ref{sec:proof-zeros-general}. The two terms in \eqref{eq:intro-zeros-general} correspond to two ways in which the entries of $\bm P_N$ vanish: over a point of $\Gamma$ the two points of the spectral curve give contributions of the same size, and over a point of $\msf z([a_1,b_1]\cap\{\SW < 0\})$ the dominant point lies on the support of $\mu$, where $f$ and $g$ oscillate.

    \begin{remark}[Critical behavior]\label{rem:critical}
        By Lemma~\ref{lem:W-positive-somewhere}, $\SW > 0$ somewhere on $[a_1,b_1]$, and, apart from the transition case discussed below, two situations can occur.
        \begin{itemize}
            \item If $\SW > 0$ on $[a_1,b_1]$, then $\rho = \msf z_*\nu$, so that all the zeros of the entries accumulate on $\Gamma$, which does not meet $\msf z([a_1,b_1])$. If moreover $\SW(0) \neq 0$, then $\rho$ is the balayage onto $\Gamma$ of $\msf z_*\mu$, the limit of the zeros of $\det\bm P_N$ (Remark~\ref{rem:balayage}).
            \item If $\SW$ changes sign on $[a_1,b_1]$, then $\Gamma_{\SW}$ crosses $[a_1,b_1]$, and $\Gamma$ reaches the real line inside $\msf z([a_1,b_1])$. A fraction $2\mu(\{\SW < 0\}) > 0$ of the zeros of each entry then accumulates on the segment $\msf z([a_1,b_1]\cap\{\SW < 0\})$, where the zeros of $\det\bm P_N$ also accumulate.
        \end{itemize}
        Figure~\ref{fig:zeros-entries} illustrates both situations for the Toda deformed weight of Section~\ref{sec:intro-toda}. The transition between them happens when $\Gamma_{\SW}$ reaches the support of $\mu$, and there the zeros of the entries display a critical behavior: a positive fraction of them moves from $\Gamma$ to the real line. The local behavior of the zeros near the points where $\Gamma$ meets the real line is not covered by our results, and we leave it for future work.
    \end{remark}

    \begin{remark}
        The convergence of the zeros of $\det\bm P_N$ in Theorem~\ref{thm:det-asymptotics} is the analogue of \cite[Theorem~2.4]{Groot2021a}. In \cite[Theorem~2.5]{Groot2021a} the limiting zero distribution of the $(1,2)$ entry, which involves a balayage, is obtained under an inequality on their $\msf g$-function, proven under a geometric condition supported by numerical evidence. This inequality means that the analogue of $\SW$ is positive on the right half-plane, so that the imaginary axis plays the role of $\Gamma_{\SW}$. Theorems~\ref{thm:intro-asymptotics} and~\ref{thm:intro-zeros-general} require no condition of this type, apply to all the entries of $\bm P_N$, and do not require $\Gamma_{\SW}$ to be known a priori. As in \cite[Theorem~2.5]{Groot2021a}, where the imaginary axis crosses the support of the analogue of $\mu$, the limit has a part on the support of $\mu$ when $\Gamma_{\SW}$ crosses it. A balayage also appears in our setting, see Remark~\ref{rem:balayage}.
    \end{remark}

    \begin{remark}\label{rem:other-weights}
        The commutation with $\bm\Lambda$ enters our approach through two facts: $\bm Q$ and $\bm B$ share the eigenvectors of a fixed matrix polynomial, which gives the decomposition \eqref{eq:intro-two-sheets}, and these eigenvectors are parametrized by a curve of genus zero with a global coordinate $x$, which makes $f$ and $g$ accessible to the steepest descent method. This structure is not specific to the Hermite-type weight.

        \emph{Degree one.} Let $\wt{\bm\Lambda}$ be a real symmetric $2\times2$ matrix polynomial of degree one, not of the form $\theta\,\bm I$. If its eigenvectors are constant, every weight that commutes with $\wt{\bm\Lambda}$ reduces to two scalar weights by a constant change of basis. Otherwise, there is a constant rotation $\bm O$ such that $\bm O^T\wt{\bm\Lambda}(z)\bm O = \theta_0(z)\bm I + \theta_1\bm\Lambda(\alpha z + \beta)$, for some $\msf a > 0$, a scalar polynomial $\theta_0$ and constants $\theta_1, \alpha \neq 0$ and $\beta$. Then, up to this rotation and an affine change of variable, the weights $\ee^{-N\bm Q}\bm B$ with $\bm Q$ and $\bm B$ commuting with $\wt{\bm\Lambda}$ are the weights considered here. Only the description of the zeros of the entries changes, since the entries in the rotated basis are linear combinations of the entries of $\bm P_N$.

        \emph{Other settings.} We expect the reduction of matrix orthogonality to scalar orthogonality of Section~\ref{sec:prelim} to extend to other varying weights $\ee^{-N\bm Q}\bm B$ in which $\bm Q$ and $\bm B$ share the eigenvectors of a fixed matrix polynomial, possibly with additional orthogonality conditions. When the corresponding spectral curve has positive genus, the steepest descent analysis has to be carried out on this curve, see \cite{Charlier2021,Bertola-2023,Bertola2021a,Bertola2022b,Bertola2022a,Kuijlaars25,Alves26} for related work.
    \end{remark}

    \begin{remark}[Recurrence coefficients]\label{rem:recurrence}
        It is natural to ask for the asymptotics of the recurrence coefficients. For each $N$, the monic MVOPs $\bm P_{n,N}$, $n \geq 0$, of the weight $\bm W_{N\bm Q}$ satisfy a three-term recurrence relation \cite{Berg,DPS}
        \[
            \zeta\,\bm P_{n,N}(\zeta) = \bm P_{n+1,N}(\zeta) + \boldsymbol\beta_{n,N}\,\bm P_{n,N}(\zeta) + \boldsymbol\gamma_{n,N}\,\bm P_{n-1,N}(\zeta), \qquad n \geq 1,
        \]
        with $\bm P_{N,N} = \bm P_N$, and our methods also give the asymptotics of $\boldsymbol\beta_N \defeq \boldsymbol\beta_{N,N}$ and $\boldsymbol\gamma_N \defeq \boldsymbol\gamma_{N,N}$. With $u \defeq \sqrt{a_1b_1}$, $m \defeq \frac{a_1 + b_1}{2}$ and $\bm\Delta_N \defeq \diag\big(1,(-1)^N\ee^{-2N\SW(0)}\big)$, one can show that, as $N \to \infty$,
        \begin{align*}
            \bm\Delta_N^{-1}\,\boldsymbol\beta_N\,\bm\Delta_N &= \frac{1}{\msf a}\begin{pmatrix} \dfrac{u - m}{u^2} & \dfrac{1 + u^2}{u^2}\,A_f(\infty) \\[2ex] -(1 + u^2)\,A_g(0) & -(u - m)\end{pmatrix} + \Boh\left(\frac1N\right), \displaybreak[1] \\[1ex]
            \bm\Delta_N^{-1}\,\boldsymbol\gamma_N\,\bm\Delta_N &= -\frac{u + m}{4\msf a^2u^2}\begin{pmatrix} \dfrac{2u^3 + u - m}{u^2} & -2(1 + u^2)\,A_f(\infty) \\[2ex] -2(1 + u^2)\,A_g(0) & u\left(2 + u(u - m)\right)\end{pmatrix} + \Boh\left(\frac1N\right),
        \end{align*}
        where $A_f$ and $A_g$ are as in Theorem~\ref{thm:intro-asymptotics-fg}. In contrast with the case of a scalar exponential factor \cite{DeanoRoman2025}, where the limits are multiples of the identity, if $\SW(0) \neq 0$ then in each of $\boldsymbol\beta_N$ and $\boldsymbol\gamma_N$ one off-diagonal entry grows like $\ee^{2N|\SW(0)|}$ and the other one decays at the same rate.

        The limits are obtained by evaluating the recurrence relation on the spectral curve with a formula analogous to \eqref{eq:intro-two-sheets}, reducing to quantities in terms of Laurent polynomials like $f$ and $g$. We do not include the details here, since the paper is already long and technical.
    \end{remark}

\subsection{The motivating example: a non-Abelian Toda deformation}\label{sec:intro-toda} \hfill

    This work was motivated by the following example. In the scalar case, the deformation $\ee^{-\msf t z}w(z)$ of a weight $w$ leads to the Toda lattice for the recurrence coefficients. Likewise, exponential deformations of matrix weights lead to non-Abelian versions of the Toda lattice \cite{Bruschi,GeK,Miranian,IKR2,BFGA_2020}. Dea\~no, Morey and Rom\'an \cite{DMR} introduced deformations $\ee^{-\msf t\bm\Lambda}\bm W$ of a matrix weight $\bm W$, in which the role of the variable in the classical Toda deformation is played by a symmetric operator of order zero $\bm\Lambda$, that is, a matrix polynomial such that $\bm\Lambda\bm W = \bm W\bm\Lambda^\ast$. They showed that the coefficients of the corresponding difference operators, which are expressed in terms of the recurrence coefficients, satisfy non-Abelian Toda-type equations. We consider
    \begin{equation}\label{eq:intro-weight}
        \bm W_{\msf a, \msf t}(z) \defeq \ee^{-\msf t\bm\Lambda(z)}\,\bm H_{\msf a}(z), \qquad \msf t \in \R, \quad \msf a > 0,
    \end{equation}
    with $\bm\Lambda$ given by \eqref{eq:intro-A-Lambda}. The example of a Hermite-type weight worked out in \cite[Section~5]{DMR} uses instead the Casimir operator $\bm J - z\bm A$, with $\bm J = \diag(1,2)$, and the corresponding MVOPs are explicit \cite[(5.5)--(5.6)]{DMR}.

    For $\msf t = 0$ the MVOPs of \eqref{eq:intro-weight} are explicit in terms of Hermite polynomials \cite{DG05}. For general $\msf t$ the weight can be identified with the setup of this paper in the scaling regime
    \begin{equation}\label{eq:scaling}
        \msf a \mapsto \frac{\msf a_0}{\sqrt N}, \qquad \msf t \mapsto N\msf t_0, \qquad z = \sqrt N\,\zeta,
    \end{equation}
    with $\msf a_0 > 0$ and $\msf t_0$ fixed. Under \eqref{eq:scaling}, $\bm W_{\msf a, \msf t}$ becomes the varying weight $\bm W_{N\bm Q}(\zeta) = \ee^{-N\bm Q(\zeta)}\bm B(\zeta)$ with
    \begin{equation}\label{eq:intro-Q-toda}
        \bm Q(\zeta) = \zeta^2\,\bm I + \msf t_0\,\bm\Lambda(\zeta), \qquad \msf Q(x) = \msf z(x)^2 + \frac{\msf t_0}{x} = \frac{1}{\msf a_0^2}\left(x - \frac1x\right)^2 + \frac{\msf t_0}{x}, \qquad
        \bm B(\zeta) = (1 + \msf a_0^2\zeta^2)\,\bm I + \msf a_0\zeta\,\bm\Lambda(\zeta),
    \end{equation}
    where $\bm\Lambda$ and $\msf z$ are now defined with $\msf a_0$ in place of $\msf a$. Theorems~\ref{thm:intro-asymptotics-fg}--\ref{thm:intro-zeros-general} apply to this example whenever the equilibrium measure $\mu = \mu_{\msf a_0,\msf t_0}$ of the external field
    \begin{equation}\label{eq:intro-V-toda}
        V(x) = V_{\msf a_0,\msf t_0}(x) = \frac{1}{2\msf a_0^2}\left(x - \frac1x\right)^2 + \frac{\msf t_0}{2x} + \log|x|, \qquad x \in \R\setminus\{0\},
    \end{equation}
    satisfies Assumption~\ref{ass:onecut}. This is the case for large values of $\msf t_0$.

    \begin{theorem}\label{thm:intro-eq}
        Fix $\msf a_0 > 0$. For $\msf t_0$ sufficiently large the equilibrium measure $\mu_{\msf a_0,\msf t_0}$ of \eqref{eq:intro-V-toda} is one-cut regular, with support contained in $\left(-\frac{1}{\msf a_0},0\right)$. In particular, Assumption~\ref{ass:onecut} holds.
    \end{theorem}

    Theorem~\ref{thm:intro-eq} is proved in Appendix~\ref{app:onecut}.

    Figure~\ref{fig:zeros-entries} illustrates these results for $\msf a_0 = 3$. For $\msf t_0 = 3$ (left), $\SW > 0$ on $[a_1,b_1]$. The sheet that dominates in \eqref{eq:intro-two-sheets} changes across $\Gamma_{\SW}$. This gives the rank-one leading term and makes the zeros of the entries accumulate on $\Gamma$, while the zeros of the determinant stay on the real line. For $\msf t_0 = 2$ (right), $\SW$ changes sign on $[a_1,b_1]$, and a positive fraction of the zeros of each entry accumulates on the real line. This does not happen for large $\msf t_0$: by Lemma~\ref{lem:localization}, $\supp\mu \subset (-1,0)$ for $\msf t_0$ sufficiently large, and then $\SW > 0$ on $\supp\mu$ by \eqref{eq:intro-W}.

\subsection{Outline} \hfill

    In Section~\ref{sec:prelim} we present the decomposition of matrix orthogonality into scalar orthogonality on the spectral curve for the weight $\bm W_{\bm Q}$. Section~\ref{sec:steepest} contains the steepest descent analysis, which we use in Section~\ref{sec:asympt1} to prove our main results. Appendix~\ref{app:eq-measure} presents some simple results on one-cut regular equilibrium measures used in this work, and proves Theorem~\ref{thm:intro-eq}. Appendix~\ref{app:prelim} collects the facts from potential theory that we use in Section~\ref{sec:asympt1}.

    \emph{Notation.} Matrices are written in boldface. $\bm E_{ij}$ denotes the $2\times2$ matrix with entry $1$ at the position $(i,j)$ and $0$ elsewhere, $\sp_3 = \diag(1,-1)$, and $\C_\pm \defeq \{z \in \C : \pm\Im z > 0\}$. Throughout, $\msf z$ and $\iota$ are as in \eqref{eq:intro-z-iota}.

    \begin{figure}[ht]
        \centering
        \captionsetup{font=small,margin=0pt}
        \includegraphics{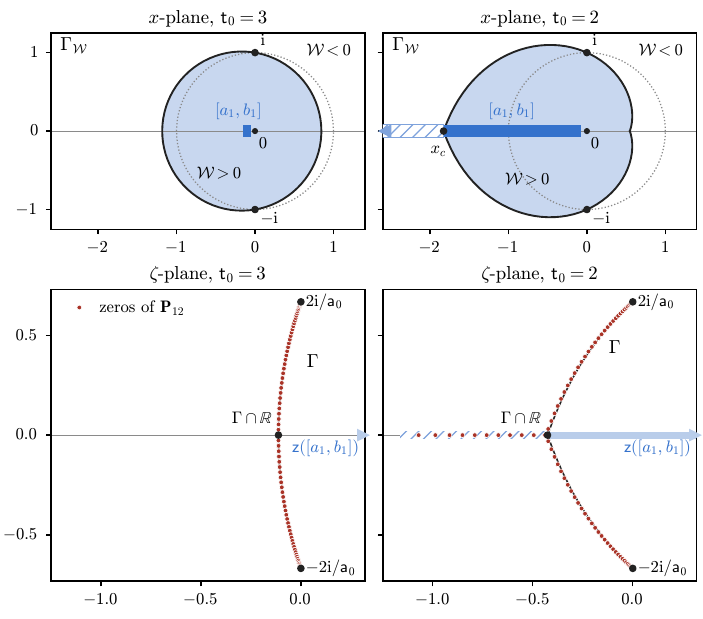}
        \caption{Zeros of the $(1,2)$ entry of $\bm P_N$, with $N = 80$, for the Toda deformed Hermite-type weight with $\msf a_0 = 3$ and $\msf t_0 = 3$ (left) or $\msf t_0 = 2$ (right). In both cases $\mu$ is one-cut regular (numerically verified). Top: the curve $\Gamma_{\SW}$ in the $x$-plane and the support $[a_1,b_1]$ of $\mu$; the dotted curve is the unit circle. For $\msf t_0 = 3$, $[a_1,b_1] \approx [-0.155,-0.053]$ and $\SW > 0$ on $[a_1,b_1]$. For $\msf t_0 = 2$, $[a_1,b_1] \approx [-3.76,-0.072]$, $\Gamma_{\SW}$ crosses it at $x_c \approx -1.83$, and $\SW < 0$ on the hatched part $[a_1,x_c)$. Bottom: the arc $\Gamma = \msf z(\Gamma_{\SW}\setminus\{0\})$, the segment $\msf z([a_1,b_1])$, which attracts the zeros of $\det\bm P_N$ (Theorem~\ref{thm:det-asymptotics}) and for $\msf t_0 = 3$ lies to the right of the window, and the zeros of $\bm P_{12}$. For $\msf t_0 = 3$ all of them accumulate on $\Gamma$. For $\msf t_0 = 2$, $\Gamma$ meets the real line at $\msf z(x_c)$, inside $\msf z([a_1,b_1])$, and a fraction $2\mu([a_1,x_c]) \approx 0.13$ of them accumulates on the hatched segment $\msf z([a_1,x_c])$ (Remark~\ref{rem:critical}). By Theorem~\ref{thm:intro-zeros-general}, the zeros of the other entries have the same limiting distribution.}
        \label{fig:zeros-entries}
    \end{figure}

\clearpage
\section{From matrix to scalar orthogonality}\label{sec:prelim}

    Throughout the paper, $\bm Q$ and $\bm B$ are real symmetric matrix polynomials of the form \eqref{eq:intro-Q}, and $\bm B(z)$ is positive definite for every $z \in \R$.

    Since $\bm Q$ and $\bm B$ commute with $\bm \Lambda$, the real symmetric matrices $\ee^{-\bm Q(z)}$ and $\bm B(z)$ commute, and the weight
    \begin{equation}\label{eq:WQ}
        \bm W_{\bm Q}(z) \defeq \ee^{-\bm Q(z)}\,\bm B(z) = \ee^{-\bm Q(z)/2}\,\bm B(z)\,\ee^{-\bm Q(z)/2}, \qquad z \in \R,
    \end{equation}
    is real, symmetric and positive definite on $\R$.

    We do not indicate the dependence of $\bm W_{\bm Q}$ on $\bm B$, which is fixed throughout. In the varying regime that we consider in Section~\ref{sec:steepest}, only $\bm Q$ is multiplied by $N$.

    We denote by $\bm P_N = (\bm P_{ij})_{i,j=1}^2$ the $N$-th monic MVOP for $\bm W_{\bm Q}$. Since $\bm W_{\bm Q}$ is positive definite, $\bm P_N$ exists and is unique provided $\bm W_{\bm Q}$ has finite moments, see for instance \cite{Berg,DPS}. A priori there is no guarantee that this is the case. Soon we verify this when $V$ given by \eqref{eq:intro-V} is admissible.

    For every $z \in \C$, the eigenvalues of $\bm\Lambda(z)$ are the numbers $\frac1x$, where $x$ is a solution of the quadratic equation $x^2 - \msf a z x - 1 = 0$, and in this case $(1,x)^T$ is an eigenvector. This leads us to consider the spectral curve
    \begin{equation}\label{eq:spectral-curve}
        \mcal R \defeq \left\{(z,x) \in \C^2 : x^2 - \msf a zx - 1 = 0\right\}.
    \end{equation}
    As $x \neq 0$ on $\mcal R$, the defining equation can be solved for $z$: $(z,x) \in \mcal R$ if, and only if, $x \neq 0$ and $z = \msf z(x)$, with $\msf z$ as in \eqref{eq:intro-z-iota}. Thus $x$ is a global coordinate on $\mcal R$, which is a sphere with the two points $x = 0$ and $x = \infty$, both lying over $z = \infty$, removed. In this coordinate the projection $(z,x) \mapsto z$ is $x \mapsto \msf z(x)$, which is two-to-one: the two solutions of $x^2 - \msf a zx - 1 = 0$ have product $-1$, so that
    \[
        \msf z(x') = \msf z(x) \iff x' \in \{x,\iota x\}.
    \]
    The involution $\iota$ exchanges the two sheets, and its fixed points $x = \pm\ii$ lie over the branch points $\msf z(\pm\ii) = \pm\frac{2\ii}{\msf a}$. We shall often use the identities
    \begin{equation}\label{eq:iota-identities}
        x + \iota x = \msf a\,\msf z(x), \qquad x\,\iota x = -1, \qquad x - \iota x = \frac{1 + x^2}{x}, \qquad \msf z'(x) = \frac{x - \iota x}{\msf a x}.
    \end{equation}
    On the real line, $\msf z' > 0$, and $\msf z$ maps each of the half-lines $(-\infty,0)$ and $(0,+\infty)$ increasingly onto $\R$, while $\iota$ maps $(0,+\infty)$ increasingly onto $(-\infty,0)$. Hence the real line of the $z$-plane is covered twice, once by $x > 0$ and once by $x < 0$.

    We start by compiling the spectral decomposition of the matrices we consider in this work.

    \begin{lemma}\label{lem:eigen-curve}
        For $x \in \C\setminus\{0\}$ and $z = \msf z(x)$,
        \begin{equation}\label{eq:eigen-curve}
            \begin{split}
            \bm\Lambda(z)\begin{pmatrix} 1 \\ x\end{pmatrix} &= \frac1x\begin{pmatrix} 1 \\ x\end{pmatrix}, \qquad
            \bm Q(z)\begin{pmatrix} 1 \\ x\end{pmatrix} = \msf Q(x)\begin{pmatrix} 1 \\ x\end{pmatrix}, \\
            \bm B(z)\begin{pmatrix} 1 \\ x\end{pmatrix} &= \msf B(x)\begin{pmatrix} 1 \\ x\end{pmatrix}, \qquad
            \bm W_{\bm Q}(z)\begin{pmatrix} 1 \\ x\end{pmatrix} = \omega_{\bm Q}(x)\begin{pmatrix} 1 \\ x\end{pmatrix},
            \end{split}
        \end{equation}
        where $\msf Q$ and $\msf B$ are the Laurent polynomials \eqref{eq:intro-Qx} and $\omega_{\bm Q}(x) \defeq \msf B(x)\,\ee^{-\msf Q(x)}$. Consequently $\bm B(z)$ is positive definite for all $z \in \R$ if, and only if, $\msf B(x)$ is positive for all $x \in \R \setminus \{0\}$. Moreover,  for $x \in \C\setminus\{0,\ii,-\ii\}$ the matrix
        \begin{equation}\label{eq:E-curve}
            \bm E(x) \defeq \begin{pmatrix} 1 & 1 \\ x & \iota x\end{pmatrix} \quad \text{is invertible}, \qquad
            \bm E(x)^{-1} = \frac{1}{x - \iota x}\begin{pmatrix} -\iota x & 1 \\ x & -1\end{pmatrix},
        \end{equation}
        and
        \begin{equation}\label{eq:spectral-decomposition}
            \bm W_{\bm Q}(\msf z(x)) = \bm E(x)\begin{pmatrix} \omega_{\bm Q}(x) & 0 \\ 0 & \omega_{\bm Q}(\iota x)\end{pmatrix}\bm E(x)^{-1}.
        \end{equation}
    \end{lemma}

    \begin{proof}
        The first identity in \eqref{eq:eigen-curve} is what motivated us to consider $\mcal R$. We repeat it here for convenience. The second and the third follow from the first one since $\bm Q = Q_0\,\bm I + Q_1\,\bm\Lambda$ and $\bm B = B_0\,\bm I + B_1\,\bm\Lambda$. By the second identity, $(1,x)^T$ is an eigenvector of $\ee^{-\bm Q(z)}$ with eigenvalue $\ee^{-\msf Q(x)}$, and the fourth identity follows from the third.

        Applying \eqref{eq:eigen-curve} at the points $x$ and $\iota x$, which lie over the same $z$, we get \[\bm W_{\bm Q}(z)\bm E(x) = \bm E(x)\diag(\omega_{\bm Q}(x),\omega_{\bm Q}(\iota x)).\] Finally, $\det\bm E(x) = \iota x - x$ vanishes only at $x = \pm\ii$, and the formula for $\bm E(x)^{-1}$ is checked directly.
    \end{proof}

    The functions $\msf Q(x)$ and $\msf B(x)$ are Laurent polynomials. We write
    \begin{equation}\label{eq:Q-Laurent}
        \begin{split}
        \msf Q(x)
            &=  \sum_{j=-d_-}^{d_+}\msf Q_j\,x^j, \qquad \msf Q_{-d_-} \neq 0, \quad \msf Q_{d_+} \neq 0,
        \end{split}
    \end{equation}
    and it is straightforward to check that $V$ defined in \eqref{eq:intro-V} is admissible if, and only if, $d_-$ and $d_+$ are even and positive with $\msf Q_{d_+} > 0$ and $\msf Q_{-d_-} > 0$. In particular, since $\msf B$ is a Laurent polynomial, $\omega_{\bm Q}(x)$ decays faster than any power of $x$ as $x \to 0$ and as $|x| \to \infty$, and $\bm W_{\bm Q}$ has finite moments by \eqref{eq:spectral-decomposition}.

    Since $\msf z(\iota x) = \msf z(x)$, we have
    \begin{equation}\label{eq:PN-laurent}
        \bm P_N(\msf z(x))\,\bm E(x) = \begin{pmatrix} f(x) & f(\iota x) \\ g(x) & g(\iota x)\end{pmatrix}, \qquad x \in \C\setminus\{0\},
    \end{equation}
    where
    \begin{equation}\label{eq:fg-def}
        f(x) \defeq \bm P_{11}(\msf z(x)) + x\,\bm P_{12}(\msf z(x)), \qquad
        g(x) \defeq \bm P_{21}(\msf z(x)) + x\,\bm P_{22}(\msf z(x)), \qquad x \in \C\setminus\{0\}.
    \end{equation}
    By \eqref{eq:intro-z-iota}, the functions $f$ and $g$ are Laurent polynomials in $x$. Equation \eqref{eq:PN-laurent} allows us to characterize the entries of $\bm P_N$ and its determinant in terms of $f$, $g$ and $\iota$.

    \begin{lemma}\label{lem:divided-differences}
        For $x \in \C\setminus\{0,\ii,-\ii\}$ and $z = \msf z(x)$
        \begin{equation}\label{eq:divided-differences}
            \bm P_{11}(z) = \frac{x\,f(\iota x) - \iota x\,f(x)}{x - \iota x}, \qquad
            \bm P_{12}(z) = \frac{f(x) - f(\iota x)}{x - \iota x},
        \end{equation}
        the same identities hold for $\bm P_{21}$ and $\bm P_{22}$ with $f$ replaced by $g$, and
        \begin{equation}\label{eq:det-fg}
            \det\bm P_N(z) = \frac{f(\iota x)\,g(x) - f(x)\,g(\iota x)}{x - \iota x}.
        \end{equation}
    \end{lemma}

    \begin{proof}
        By \eqref{eq:PN-laurent}, $\bm P_N(\msf z(x)) = \begin{pmatrix} f(x) & f(\iota x) \\ g(x) & g(\iota x)\end{pmatrix}\bm E(x)^{-1}$. The formulas \eqref{eq:divided-differences} follow from the expression of $\bm E(x)^{-1}$ in Lemma~\ref{lem:eigen-curve}, and \eqref{eq:det-fg} follows by taking determinants, since $\det\bm E(x) = \iota x - x$.
    \end{proof}

    We now investigate how the orthogonality relations satisfied by $\bm P_N$ manifest on the functions $f$ and $g$. Let
    \begin{equation}\label{eq:RL-weight}
        w_{\bm Q}(x) \defeq \lambda(x)\,\ee^{-\msf Q(x)}, \qquad \lambda(x) \defeq \frac{\msf B(x)}{x^2}, \qquad x \in \R\setminus\{0\},
    \end{equation}
    so that $\omega_{\bm Q}(x) = x^2w_{\bm Q}(x)$ in Lemma~\ref{lem:eigen-curve}. The weight $w_{\bm Q}$ is positive and decays faster than any (integer) power of $x$ as $x \to 0$ and as $|x| \to \infty$. In particular, $\int_\R h(x)\,w_{\bm Q}(x)\,\dd x$ converges absolutely for every Laurent polynomial $h$.

    \begin{theorem}\label{thm:orth-fg}
        The Laurent polynomials $f$ and $g$ satisfy the orthogonality conditions
        \[
            \int_\R f(x)\,x^k\,w_{\bm Q}(x)\,\dd x = \int_\R g(x)\,x^k\,w_{\bm Q}(x)\,\dd x = 0, \qquad -N + 1 \leq k \leq N.
        \]
    \end{theorem}

    \begin{proof}
        Fix $0 \leq j \leq N - 1$. The orthogonality of $\bm P_N$ gives $\int_\R \bm P_N(z)\,\bm W_{\bm Q}(z)\,z^j\,\dd z = \bm 0$. We parametrize the real line of the $z$-plane by $z = \msf z(x)$ with $x \in (0,+\infty)$, so that the two points over $z$ are $x$ and $\iota x \in (-\infty,0)$. By the spectral decomposition of $\bm W_{\bm Q}$ in \eqref{eq:spectral-decomposition} and  \eqref{eq:PN-laurent} we write
        \[
            \bm P_N(z)\,\bm W_{\bm Q}(z) = \begin{pmatrix} f(x) & f(\iota x) \\ g(x) & g(\iota x)\end{pmatrix}\begin{pmatrix} \omega_{\bm Q}(x) & 0 \\ 0 & \omega_{\bm Q}(\iota x)\end{pmatrix}\bm E(x)^{-1}.
        \]
        By \eqref{eq:E-curve} and \eqref{eq:iota-identities}, the first row of $\bm P_N(z)\,\bm W_{\bm Q}(z)\,\msf z'(x)$ is
        \[
            \frac{f(x)\,\omega_{\bm Q}(x)}{\msf a x}\begin{pmatrix} -\iota x & 1\end{pmatrix} + \frac{f(\iota x)\,\omega_{\bm Q}(\iota x)}{\msf a x}\begin{pmatrix} x & -1\end{pmatrix}
            = \frac{1}{\msf a}\left[f(x)\,w_{\bm Q}(x)\begin{pmatrix} 1 & x\end{pmatrix} + f(\iota x)\,w_{\bm Q}(\iota x)\begin{pmatrix} 1 & \iota x\end{pmatrix}\frac{1}{x^2}\right],
        \]
        where we used $\omega_{\bm Q} = x^2w_{\bm Q}$, $x\,\iota x =- 1$ in the first term, and $(\iota x)^2(x,-1)/x = x^{-2}(1,\iota x)$ in the second one. The change of variables $z = \msf z(x)$, followed in the second term by $y = \iota x$, gives
        \[
            \bm 0 = \frac{1}{\msf a}\int_0^{+\infty} f(x)\begin{pmatrix} 1 & x\end{pmatrix}\msf z(x)^j\,w_{\bm Q}(x)\,\dd x + \frac{1}{\msf a}\int_{-\infty}^0 f(y)\begin{pmatrix} 1 & y\end{pmatrix}\msf z(y)^j\,w_{\bm Q}(y)\,\dd y
        \]
        for the first row of $\int_\R \bm P_N(z)\,\bm W_{\bm Q}(z)\,z^j\,\dd z$. Hence
        \[
            \int_\R f(x)\,\msf z(x)^j\,w_{\bm Q}(x)\,\dd x = \int_\R f(x)\,x\,\msf z(x)^j\,w_{\bm Q}(x)\,\dd x = 0, \qquad 0 \leq j \leq N - 1.
        \]
        The $2N$ Laurent polynomials $(\msf a\msf z(x))^j = (x - x^{-1})^j$ and $x\,(\msf a\msf z(x))^j = x\,(x - x^{-1})^j$, $0 \leq j \leq N - 1$, belong to the $2N$-dimensional space spanned by $x^k$, $-N + 1 \leq k \leq N$, and they are linearly independent. Hence they span that space, which proves the statement for $f$. The same analysis for the second row gives the statement for $g$.
    \end{proof}

    The argument is reversible: for every $2\times2$ matrix polynomial $\bm P$, with $f$ and $g$ defined from $\bm P$ as in \eqref{eq:fg-def}, the conditions of Theorem~\ref{thm:orth-fg} are equivalent to $\int_\R \bm P(z)\,\bm W_{\bm Q}(z)\,z^j\,\dd z = \bm 0$ for $0 \leq j \leq N - 1$. In particular, a monic matrix polynomial $\bm P$ of degree $N$ is the MVOP $\bm P_N$ if, and only if, its functions $f$ and $g$ satisfy the conditions of Theorem~\ref{thm:orth-fg}.

    Since $f$ and $g$ are Laurent polynomials, we focus our analysis on their numerators. Let
    \begin{equation}\label{eq:pq-def}
        p(x) \defeq x^N f(x), \qquad q(x) \defeq x^{N-1}g(x).
    \end{equation}

    \begin{corollary}\label{cor:orth-pq}
        The polynomials $p$ and $q$ of \eqref{eq:pq-def} satisfy the orthogonality conditions
        \begin{equation}\label{eq:orth-conditions-pq}
            \int_\R x\,p(x)\,x^k\,x^{-2N}w_{\bm Q}(x)\,\dd x = \int_\R q(x)\,x^k\,x^{-2(N-1)}w_{\bm Q}(x)\,\dd x = 0, \qquad 0 \leq k \leq 2N - 1.
        \end{equation}
        In particular,
        \begin{itemize}
            \item $p$ is a polynomial of degree at most $2N$ which is orthogonal to all polynomials of degree at most $2N - 1$ with respect to the real-valued weight $x^{-2N+1}w_{\bm Q}(x)$;
            \item $\msf a^Nq$ is the monic orthogonal polynomial of degree $2N$ for the positive weight $x^{-2(N-1)}w_{\bm Q}(x)$.
        \end{itemize}
        Moreover,
        \begin{equation}\label{eq:pq-at-zero}
            p(0) = \frac{(-1)^N}{\msf a^N}.
        \end{equation}
    \end{corollary}

    \begin{proof}
        By Theorem~\ref{thm:orth-fg} and \eqref{eq:pq-def}, for $-N + 1 \leq j \leq N$,
        \[
            0 = \int_\R f(x)\,x^j\,w_{\bm Q}(x)\,\dd x = \int_\R p(x)\,x^{j-N}\,w_{\bm Q}(x)\,\dd x = \int_\R p(x)\,x^{j+N-1}\,x^{-2N+1}\,w_{\bm Q}(x)\,\dd x,
        \]
        and with $k = j + N - 1$ this is the orthogonality of $p$ to $x^k$, $0 \leq k \leq 2N - 1$, with respect to $x^{-2N+1}w_{\bm Q}(x)$, which is the first identity in \eqref{eq:orth-conditions-pq}. In the same way, for $-N + 1 \leq j \leq N$,
        \[
            0 = \int_\R g(x)\,x^j\,w_{\bm Q}(x)\,\dd x = \int_\R q(x)\,x^{j-N+1}\,w_{\bm Q}(x)\,\dd x = \int_\R q(x)\,x^{j+N-1}\,x^{-2(N-1)}\,w_{\bm Q}(x)\,\dd x,
        \]
        which is the second identity in \eqref{eq:orth-conditions-pq}.

        We have $\deg p \leq 2N$ and $\deg q \leq 2N$ by the definitions of $p$ and $q$ and the fact that $\bm P_N$ is monic of degree $N$. Moreover, the only term of $f$ in \eqref{eq:fg-def} containing $x^{-N}$ comes from $\msf z(x)^N$, so that $p(0) = (-1/\msf a)^N$. Similarly, the only term of $g$ containing $x^{N+1}$ comes from $x\,\msf z(x)^N$, so that $q$ has degree $2N$, with leading coefficient $\msf a^{-N}$.
    \end{proof}

\section{Steepest descent analysis}\label{sec:steepest}

\subsection{Riemann--Hilbert problems for \texorpdfstring{$p$}{p} and \texorpdfstring{$q$}{q}} \hfill

    For every even integer $n$, let $P_{m,n}$ be the $m$-th monic orthogonal polynomial for the positive weight $x^n w_{\bm Q}(x)$ on $\R$, and $h_{m,n} \defeq \int_\R P_{m,n}(x)^2\, x^n w_{\bm Q}(x)\,\dd x > 0$. By Corollary~\ref{cor:orth-pq},
    \[
        q(x) = \msf a^{-N} P_{2N,-2(N-1)}(x).
    \]
    The polynomial $p$ can be characterized in terms of a kernel polynomial.

    \begin{lemma}\label{lem:CD-normalization}
        Let
        \[
            K_n(x,y) \defeq \sum_{k=0}^{n-1} \frac{P_{k,-2N}(x)\, P_{k,-2N}(y)}{h_{k,-2N}}
        \]
        be the Christoffel--Darboux kernel of the weight $x^{-2N} w_{\bm Q}(x)$. Then
        \begin{equation}\label{eq:p-CD}
            p(x) = p(0)\, \frac{K_{2N+1}(x,0)}{K_{2N+1}(0,0)}, \qquad K_{2N+1}(0,0) > 0.
        \end{equation}
    \end{lemma}
    \begin{proof}
        Write $\omega(x) \defeq x^{-2N} w_{\bm Q}(x)$ and $\langle u, v \rangle_\omega \defeq \int_\R u(x) v(x)\, \omega(x) \, \dd x$. The weight $\omega$ is positive on $\R \setminus \{0\}$, and all its moments are finite because $w_{\bm Q}(x)$ decays faster than any power of $x$ both as $x \to 0$ and as $|x| \to \infty$. Hence $\langle \cdot, \cdot \rangle_\omega$ is positive definite on the space $\Pi_{2N}$ of polynomials of degree at most $2N$.

        By Corollary~\ref{cor:orth-pq}, $\deg p \leq 2N$, and \eqref{eq:orth-conditions-pq} says that $\langle p, x^m \rangle_\omega = 0$ for $1 \leq m \leq 2N$. Thus $p$ lies in the orthogonal complement of $x\Pi_{2N-1}$ in $\Pi_{2N}$, which is a one-dimensional linear space.

        On the other hand, $K_{2N+1}(\cdot,0) \in \Pi_{2N}$ has the reproducing property $\langle K_{2N+1}(\cdot,0), \pi \rangle_\omega = \pi(0)$ for $\pi \in \Pi_{2N}$. Taking $\pi(x) = x^m$ with $1 \leq m \leq 2N$ shows that $K_{2N+1}(\cdot,0)$ lies in the same orthogonal complement, and it is not identically zero because
        \[
            K_{2N+1}(0,0) = \sum_{k=0}^{2N} \frac{P_{k,-2N}(0)^2}{h_{k,-2N}} \geq \frac{P_{0,-2N}(0)^2}{h_{0,-2N}} = \frac{1}{h_{0,-2N}} > 0.
        \]
        Therefore $p = c\, K_{2N+1}(\cdot,0)$ for some constant $c$, and evaluating at $x = 0$ gives $c = p(0)/K_{2N+1}(0,0)$.
    \end{proof}

    For even $n$ and $m \geq 1$, the matrix
    \[
        \bm Y^{(n,m)}(z) \defeq
        \begin{pmatrix}
            P_{m,n}(z)    &   \dfrac{1}{2 \pi \ii} \displaystyle\int_{\R} \frac{P_{m,n}(y)}{y-z} y^n w_{\bm Q}(y)\, \dd y \\[2ex]
            -\dfrac{2\pi \ii}{h_{m-1,n}} P_{m-1,n}(z)  &   -\dfrac{1}{h_{m-1,n}} \displaystyle\int_{\R} \frac{P_{m-1,n}(y)}{y-z} y^n w_{\bm Q}(y) \, \dd y
        \end{pmatrix}
    \]
    is the unique solution of the Fokas--Its--Kitaev RH problem \cite{FokasItsKitaev90,FokasItsKitaev91} for the weight $x^n w_{\bm Q}$.

    With $\bm Y = \bm Y^{(-2N,2N+1)}$ the Christoffel--Darboux formula with $y=0$ reads
    \[
        K_{2N+1}(x,0) = \frac{\bm Y_{11}(0)\, \bm Y_{21}(x) - \bm Y_{21}(0)\, \bm Y_{11}(x)}{2\pi\ii\, x}, \qquad
        K_{2N+1}(0,0) = \frac{\bm Y_{11}(0)\, \bm Y_{21}'(0) - \bm Y_{21}(0)\, \bm Y_{11}'(0)}{2\pi\ii},
    \]
    the second identity being obtained by letting $x \to 0$ in the first one. Together with \eqref{eq:p-CD} and \eqref{eq:pq-at-zero}, this expresses $p$ through the first column of the single RHP solution $\bm Y^{(-2N,2N+1)}$.

    We now replace $\bm Q$ by $N\bm Q$, that is, we consider the varying weight $\bm W_{N\bm Q} = \ee^{-N\bm Q}\bm B$, for which $w_{N\bm Q}(x) = \lambda(x)\,\ee^{-N\msf Q(x)}$. With the external field $V$ of \eqref{eq:intro-V}, we have
    \begin{equation}\label{eq:bj-def}
        x^{-2(N-j)} w_{N\bm Q}(x) = \lambda_j(x)\, \ee^{-2NV(x)}, \qquad \lambda_j(x) \defeq x^{2j}\,\lambda(x),
    \end{equation}
    so that $\lambda_0 = \lambda$ and $\lambda_1 = \msf B$. The functions $\lambda_j$ are Laurent polynomials that are positive on $\R\setminus\{0\}$ and do not depend on $N$.

    Our objective is to compute the asymptotics of $p(x)$ and $q(x)$. This specializes the RHP to $m = 2N+k$ and $n = -2(N-j)$ with $j, k \in \{0,1\}$.
    \begin{rhp}\label{rhp:Y}
        Given $N \in \Z^*_+$, find $\bm Y = \bm Y^{[j,k]}: \C \setminus \R \to \C^{2 \times 2}$ satisfying the following conditions:
        \begin{enumerate}[(i)]
            \item $\bm Y$ is an analytic function.
            \item $\bm Y$ admits continuous boundary values at $\R$, that relate to each other through the identity
            \[
                \bm Y_+(x) = \bm Y_-(x) \bm J_{\bm Y}(x), \quad
                \bm J_{\bm Y}(x) = \begin{pmatrix}
                    1   &   \lambda_j(x) \ee^{- 2N V(x)}   \\
                    0   &   1
                \end{pmatrix}
                = \bm I + \lambda_j(x) \ee^{-2N V(x)}\bm E_{12}, \quad x \in \R.
            \]
            \item As $z \to \infty$,
            \[
                \bm Y(z) =
                    \left( \bm I + \frac{\bm Y_1}{z} + \Boh\left( \frac{1}{z^2}\right) \right)
                    z^{(2N+k) \sp_3}, \quad \bm Y_1 = \bm Y_1^{[j,k]}.
            \]
        \end{enumerate}
    \end{rhp}

    By uniqueness, $\bm Y^{[j,k]} = \bm Y^{(-2(N-j),\, 2N+k)}$, and the identity $q = \msf a^{-N}P_{2N,-2(N-1)}$ reads
    \begin{equation}\label{eq:q-via-Y-rescaled}
        q(x) = \msf a^{-N}\, \bm Y^{[1,0]}_{11}(x).
    \end{equation}

    In the same way, \eqref{eq:p-CD}, \eqref{eq:pq-at-zero} and the Christoffel--Darboux formula give
    \begin{equation}\label{eq:p-via-Y-rescaled}
        p(x) = \left( -\frac{1}{\msf a}\right)^N \frac{\bm Y^{[0,1]}_{11}(0)\, \bm Y^{[0,1]}_{21}(x) - \bm Y^{[0,1]}_{21}(0)\, \bm Y^{[0,1]}_{11}(x)}{x \left( \bm Y^{[0,1]}_{11}(0)\, \big(\bm Y^{[0,1]}_{21}\big)'(0) - \bm Y^{[0,1]}_{21}(0)\, \big(\bm Y^{[0,1]}_{11}\big)'(0)\right)}.
    \end{equation}
    Hence only the pairs $(j,k) \in \{(0,1), (1,0)\}$ are needed.

\subsection{Equilibrium measure and related quantities}\label{sec:eq-measure} \hfill

    To prepare ourselves for the asymptotic analysis we now compile the necessary functions and auxiliary results. The quantities used below are given in terms of the equilibrium measure of $\R$ under the external field $V$.

    For a finite positive measure $\sigma$ with compact support in $\C$ we denote by
    \[
        \U^\sigma(z) \defeq \int\log\frac{1}{|z - s|}\,\dd\sigma(s), \qquad \CT^\sigma(z) \defeq \int\frac{\dd\sigma(s)}{s - z}, \quad z \in \C\setminus\supp\sigma,
    \]
    its logarithmic potential and its Cauchy transform. Recall that the equilibrium measure $\mu$ is the minimizer of the energy functional $\I^V$, defined in Section~\ref{sec:intro}, over all probability measures on $\R$. Since $V$ is admissible, the equilibrium measure exists, is unique and has compact support, and there is a constant $\ell$ such that the Euler--Lagrange variational conditions hold
    \[
        \begin{split}
        2 \U^\mu + V &= \ell \quad \text{on } \supp \mu, \\
        2 \U^\mu + V &\geq \ell \quad \text{on } \R \setminus \supp \mu,
        \end{split}
    \]
    in the quasi-everywhere sense \cite[Theorem~I.1.5]{SaffTotik2024}. Consider the \textit{effective potential}
    \begin{equation}\label{eq:F-def}
        F(x) \defeq 2\U^\mu(x) + V(x) - \ell, \qquad x \in \R \setminus \{0\}.
    \end{equation}
    By Lemma~\ref{lem:quad}, $\supp\mu$ is a compact subset of $\R\setminus\{0\}$, $\U^\mu$ is continuous on $\C$, $F \geq 0$ on $\R\setminus\{0\}$, and $F = 0$ on $\supp\mu$.

    \begin{definition}[One-cut and regular measures]\label{def:regular}
        We say that $\mu$ is \emph{one-cut} if $\supp \mu$ is an interval $[a_1,b_1]$. A one-cut $\mu$ is called \emph{regular} if $\mu$ is absolutely continuous with respect to the Lebesgue measure on $\R$ and satisfies
        \begin{enumerate}[\rm (R1)]
            \item $F > 0$ on $\R \setminus \left( [a_1,b_1] \cup \{0\} \right)$,
            \item the density $\mu'$ is positive on $(a_1,b_1)$, and
            \item the density $\mu'$ vanishes as a square root at $a_1$ and at $b_1$.
        \end{enumerate}
    \end{definition}

    By Assumption~\ref{ass:onecut}, $\mu$ is one-cut and regular, with $[a_1,b_1] \subset (-\infty, 0)$. We put
    \begin{equation}\label{eq:m-r-u}
        m \defeq \frac{a_1 + b_1}{2}, \qquad r \defeq \frac{b_1 - a_1}{2}, \qquad \msf h(z) \defeq \left( (z - a_1)(z - b_1) \right)^{\frac12}, \quad u \defeq \msf h(0) = \sqrt{a_1 b_1} > 0,
    \end{equation}
    where $\msf h$ is analytic in $\C \setminus [a_1,b_1]$ and $\msf h(z) = z + \Boh(1)$ as $z \to \infty$.  With $d_\pm$ and $\msf Q_{d_+}$ as in \eqref{eq:Q-Laurent}, we also put
    \begin{equation}\label{eq:k-def}
        \msf k \defeq \frac{4}{d_+\,\msf Q_{d_+}} > 0, \qquad\text{so that}\qquad \frac{V'(z)}{2} = \frac{z^{d_+ - 1}}{\msf k}\left(1 + \Boh\left(z^{-1}\right)\right), \quad z \to \infty.
    \end{equation}

    The proof of the following proposition can be found in Appendix~\ref{app:support}.

    \begin{proposition}\label{prop:onecut}
        If $\supp\mu = [a_1,b_1] \subset (-\infty,0)$, then $0 < u < |m|$, and there is a real monic polynomial $\msf q$ of degree $d_+ + d_- - 1$ such that
        \begin{equation}\label{eq:Cauchy-h-q}
            \CT^\mu(z) + \frac{V'(z)}{2} = \frac{\msf h(z)\, \msf q(z)}{\msf k\, z^{d_- + 1}}, \qquad z \in \C \setminus \left( [a_1,b_1] \cup \{0\} \right),
        \end{equation}
        and
        \begin{equation}\label{eq:density-onecut}
            \mu'(x) = \frac{\sqrt{(x - a_1)(b_1 - x)}\; \msf q(x)}{\pi\, \msf k\, x^{d_- + 1}}, \qquad x \in [a_1,b_1].
        \end{equation}
        In particular $\msf q \leq 0$ on $[a_1,b_1]$, and (R2) and (R3) hold if, and only if, $\msf q < 0$ on $[a_1,b_1]$.
    \end{proposition}

    In particular, under Assumption~\ref{ass:onecut}, \eqref{eq:Cauchy-h-q} holds and $\msf q < 0$ on $[a_1,b_1]$.

    The complexified logarithmic potential (also known as $\msf g$-function)
    \[
        \msf g(z) \defeq \int \log (z-w) \, \dd \mu(w), \quad z \in \C \setminus (-\infty, b_1],
    \]
    with the principal branch of the logarithm, is analytic, satisfies $\msf g' = -\CT^\mu$ and $\msf g(z) = \log z + \Boh(z^{-1})$ as $z \to \infty$. Moreover, $\msf g$ has continuous boundary values on $(-\infty, b_1)$, related by
    \begin{equation}\label{eq:g-jumps}
        \msf g_+(x) + \msf g_-(x) = -2 \U^\mu(x), \qquad
        \msf g_+(x) - \msf g_-(x) = 2 \pi \ii\, \mu([x, b_1]), \qquad x \in \R,
    \end{equation}
    so that $\msf g_+ - \msf g_- = 2\pi\ii$ on $(-\infty,a_1)$ and $\msf g_+ = \msf g_-$ on $(b_1,+\infty)$. In particular, $\ee^{\msf g}$ extends to an analytic and non-vanishing function on $\C\setminus[a_1,b_1]$, which we still denote by $\ee^{\msf g}$. It satisfies $|\ee^{\msf g}| = \ee^{-\U^\mu}$, $\ee^{\msf g(x)} = x\left(1 + \Boh(x^{-1})\right)$ as $x \to \infty$, and $\msf g(0) = \int\log|s|\,\dd\mu(s) \in \R$.

    In terms of $\msf g$, the Euler--Lagrange conditions and (R1) read
    \begin{equation}\label{eq:euler-lagrange}
        \msf g_+(x) + \msf g_-(x) - V(x) + \ell = -F(x)
        \begin{cases}
            = 0, & x \in [a_1,b_1], \\
            < 0, & x \in \R \setminus \left( [a_1,b_1] \cup \{0\} \right).
        \end{cases}
    \end{equation}

    We shall also make use of the auxiliary functions (sometimes called $\phi$-functions or phase functions)
    \begin{equation}\label{eq:phi-functions}
        \phi_{a_1}(z) \defeq \int_{a_1}^z \frac{\msf h(s)\, \msf q(s)}{\msf k\, s^{d_- + 1}}\, \dd s, \qquad
        \phi_{b_1}(z) \defeq \int_{b_1}^z \frac{\msf h(s)\, \msf q(s)}{\msf k\, s^{d_- + 1}}\, \dd s,
    \end{equation}
    defined for $z \in \C \setminus [a_1, +\infty)$ and for $z \in \C \setminus \bigl( (-\infty, b_1] \cup [0, +\infty) \bigr)$, respectively, the path of integration lying in these simply connected domains. By \eqref{eq:Cauchy-h-q} the integrand is $\CT^\mu + V'/2$. It has a pole at $s = 0$ and a non-zero period around $[a_1,b_1]$, which is why the domains are chosen as above.

    The following lemmas compile a series of properties used in the asymptotic analysis. Their proofs are standard for one-cut regular measures, so we omit them.

    \needspace{6\baselineskip}
    \begin{lemma}\label{lem:phi-properties} \hfill

        \begin{enumerate}[\rm (i)]
            \item For $x \in (a_1,b_1)$, $\phi_{a_1,\pm}(x) = \pm \pi \ii\, \mu([a_1,x])$ and $\phi_{b_1,\pm}(x) = \mp \pi \ii\, \mu([x,b_1])$. In particular
            \[
                \msf g_+(x) - \msf g_-(x) = -2\phi_{b_1,+}(x) = 2\phi_{b_1,-}(x).
            \]
            \item $\phi_{b_1}(z) - \phi_{a_1}(z) = \mp\pi\ii$ for $z \in \C_\pm$.
            \item $F(x) = 2\phi_{a_1}(x)$ for $x < a_1$, and $F(x) = 2 \phi_{b_1}(x)$ for $b_1 < x < 0$.
            \item Put
            \[
                c_{a_1} \defeq \frac{\sqrt{b_1 - a_1}\; \msf q(a_1)}{\msf k\, a_1^{d_- + 1}}, \qquad
                c_{b_1} \defeq \frac{\sqrt{b_1 - a_1}\; \msf q(b_1)}{\msf k\, b_1^{d_- + 1}}.
            \]
            Then $c_{a_1} > 0$, $c_{b_1} > 0$, and
            \[
                \phi_{a_1}(z) = \frac{2 c_{a_1}}{3} (a_1 - z)^{\frac32} \psi_{a_1}(z), \qquad
                \phi_{b_1}(z) = \frac{2 c_{b_1}}{3} (z - b_1)^{\frac32} \psi_{b_1}(z),
            \]
            with principal branches, where $\psi_{a_1}$ and $\psi_{b_1}$ are analytic in neighborhoods of $a_1$ and $b_1$, respectively, real on the real line, and equal to $1$ at $a_1$, respectively $b_1$.
        \end{enumerate}
    \end{lemma}

    \begin{lemma}[Local coordinates at the endpoints]\label{lem:zeta-maps}
        For $p \in \{a_1,b_1\}$ and $\delta > 0$ let $U_{p,\delta} \defeq \{z \in \C : |z - p| < \delta\}$, and let $c_p$ and $\psi_p$ be as in Lemma~\ref{lem:phi-properties}(iv). There is $\delta_0 > 0$ such that the closed disks $\overline{U_{a_1,\delta_0}}$ and $\overline{U_{b_1,\delta_0}}$ are disjoint and contained in $\{\Re z < 0\}$, and such that the functions
        \begin{equation}\label{eq:zeta-maps}
            \zeta_{b_1}(z) \defeq (2c_{b_1})^{\frac23} (z - b_1)\, \psi_{b_1}(z)^{\frac23}, \qquad
            \zeta_{a_1}(z) \defeq (2c_{a_1})^{\frac23} (a_1 - z)\, \psi_{a_1}(z)^{\frac23},
        \end{equation}
        with the principal branch of $\psi_p^{2/3}$, are analytic and injective in $U_{b_1,\delta_0}$ and $U_{a_1,\delta_0}$, respectively, and real on the real line. Moreover:
        \begin{enumerate}[\rm (i)]
            \item $\zeta_{b_1}$ maps $U_{b_1,\delta_0} \cap \C_\pm$ into $\C_\pm$ and $(b_1 - \delta_0, b_1)$ into $(-\infty,0)$, whereas $\zeta_{a_1}$ maps $U_{a_1,\delta_0} \cap \C_\pm$ into $\C_\mp$ and $(a_1, a_1 + \delta_0)$ into $(-\infty,0)$;
            \item with principal branches,
            \begin{equation}\label{eq:zeta-phi}
                \zeta_{b_1}(z)^{\frac32} = 3\phi_{b_1}(z), \quad z \in U_{b_1,\delta_0} \setminus (b_1 - \delta_0, b_1], \qquad
                \zeta_{a_1}(z)^{\frac32} = 3\phi_{a_1}(z), \quad z \in U_{a_1,\delta_0} \setminus [a_1, a_1 + \delta_0).
            \end{equation}
        \end{enumerate}
        These properties persist if $\delta_0$ is decreased.
    \end{lemma}

\subsection{Steepest descent analysis}\label{sec:onecut} \hfill

    The steepest descent analysis of $\bm Y$ follows the standard scheme of \cite{Egregium,deift_exp_weights}. For the purposes of this work the most important part of the whole analysis is contained in Lemma~\ref{lem:first-columns}, where we characterize the first columns of the global parametrix, which will describe $f$ and $g$ asymptotically in Section~\ref{sec:asympt1}. We omit the dependence on $j$ and $k$ of the transformations whenever it is not essential.

    We perform the analysis using $z$ as a variable, but note that the whole procedure is done in the $x$-plane.

\subsubsection{First and second transformations: the $\msf g$-function and the opening of lenses} \hfill

    By \eqref{eq:g-jumps} and \eqref{eq:euler-lagrange}, and since $\msf g(z) = \log z + \Boh(z^{-1})$ as $z \to \infty$, the matrix-valued function
    \[
        \bm T(z)
        = \ee^{ N \ell \sp_3} \bm Y(z) \ee^{- 2 N \msf g(z) \sp_3} \ee^{-N \ell \sp_3}
    \]
    satisfies the following RHP.

    \begin{rhp}\label{rhp:T}
        Seek $\bm T: \C \setminus \R \to \C^{ 2 \times 2}$ such that:
        \begin{enumerate}[(I)]
            \item $\bm T$ is analytic.
            \item $\bm T$ admits the boundary values
            \[
                \bm T_+(x) = \bm T_-(x) \bm J_{\bm T}(x), \quad x \in \R
            \]
            where
            \[
                \begin{split}
                \bm J_{\bm T}(x)
                 &\defeq \ee^{-2N(\msf g_+(x) - \msf g_-(x))\sp_3} + \lambda_j(x) \ee^{2N(\msf g_+(x) + \msf g_-(x) - V(x) + \ell)} \bm E_{12} \\
                &=\begin{dcases}
                    \begin{pmatrix}
                        \ee^{- 2 N( \msf g_+(x) - \msf g_-(x))} &   \lambda_j(x) \\
                        0   &   \ee^{ 2 N( \msf g_+(x) - \msf g_-(x))}
                    \end{pmatrix}, \quad &x \in (a_1,b_1),\\
                    \begin{pmatrix}
                       1    &   \lambda_j(x)  \ee^{-2N F(x)}   \\
                       0    &   1
                    \end{pmatrix}, \quad &x \in \R \setminus [a_1,b_1],
                \end{dcases}
                \end{split}
            \]
            with $\bm J_{\bm T}(0) = \bm I$.
            \item As $z \to \infty$,
            \[
                \bm T(z) = \left( \bm I + \frac{\bm T_1}{z} + \Boh( z^{-2}) \right) z^{k \sp_3}.
            \]
        \end{enumerate}
    \end{rhp}

    By Lemma~\ref{lem:phi-properties}(i), for $x \in (a_1,b_1)$ we have
    \[
        \bm J_{\bm T}(x) =
            \ee^{4N \phi_{b_1,+}(x) \sp_3} + \lambda_j(x) \bm E_{12}, \qquad \phi_{b_1,-}(x) = -\phi_{b_1,+}(x).
    \]
    The opening of lenses is motivated by the standard decomposition
    \[
        \begin{pmatrix}
            \ee^{4N \phi_{b_1,+}} &   \lambda_j  \\
            0   &   \ee^{-4N \phi_{b_1,+}}
        \end{pmatrix}
        =
        \begin{pmatrix}
            1   &   0\\
            \lambda_j^{-1}\,\ee^{-4N \phi_{b_1,+}}    & 1
        \end{pmatrix}
        \begin{pmatrix}
            0   &   \lambda_j  \\
            - \lambda_j^{-1}   &   0
        \end{pmatrix}
        \begin{pmatrix}
            1   &   0\\
            \lambda_j^{-1}\,\ee^{4N \phi_{b_1,+}}    & 1
        \end{pmatrix},
    \]
    where all the functions are evaluated at $x \in (a_1,b_1)$.

    The functions $\lambda_0$ and $\lambda_1$ are analytic in $\C\setminus\{0\}$, and their zeros are those of $\msf B$, which lie off the real line. We fix a neighborhood $\mcal V_{\lambda}$ of $[a_1,b_1]$, contained in $\{\Re z < 0\}$, in which $\msf B$ does not vanish. The lenses and the disks around $a_1$ and $b_1$ below are chosen inside $\mcal V_{\lambda}$, which is possible by Lemma~\ref{lem:lens} and by taking $\delta$ small.

    Let $\delta_0$, $\zeta_{a_1}$ and $\zeta_{b_1}$ be as in Lemma~\ref{lem:zeta-maps}. For $0 < \delta < \delta_0$, the lips of the lenses inside the disks $U_{a_1,\delta}$ and $U_{b_1,\delta}$ are defined by
    \begin{equation}\label{eq:lens-near-endpoints}
        \Gamma_{\bm S}^\pm \cap U_{b_1,\delta} \defeq U_{b_1,\delta} \cap \zeta_{b_1}^{-1}\left( \ee^{\pm \frac{2\pi\ii}{3}} [0,+\infty) \right), \qquad
        \Gamma_{\bm S}^\pm \cap U_{a_1,\delta} \defeq U_{a_1,\delta} \cap \zeta_{a_1}^{-1}\left( \ee^{\mp \frac{2\pi\ii}{3}} [0,+\infty) \right),
    \end{equation}
    where $\zeta_p^{-1}(\cdot)$ denotes preimages under $\zeta_p \colon U_{p,\delta_0} \to \C$. By Lemma~\ref{lem:zeta-maps}(i), the sets \eqref{eq:lens-near-endpoints} with the sign $+$ lie in the closed upper half plane, and those with the sign $-$ in the closed lower half plane. Outside the two disks, the lenses are chosen according to the following lemma, whose proof is standard for one-cut regular measures using the Cauchy--Riemann equations.

    \begin{lemma}\label{lem:lens}
        If $\delta \in (0,\delta_0)$ is small enough, there are curves $\Gamma_{\bm S}^\pm$ from $a_1$ to $b_1$ in the upper, respectively lower, half plane, contained in $\{\Re z < 0\}$ and satisfying \eqref{eq:lens-near-endpoints}, such that $\Re \phi_{b_1} < 0$ on $\Gamma_{\bm S}^\pm \setminus \{a_1,b_1\}$ and on the regions $\Omega_{\bm S}^\pm$ enclosed between $\Gamma_{\bm S}^\pm$ and $[a_1,b_1]$. In particular, for every $\delta' > 0$ there is $\eta > 0$ with $\Re\phi_{b_1} \leq -\eta$ on $\Gamma_{\bm S}^\pm \setminus \left( U_{a_1,\delta'} \cup U_{b_1,\delta'} \right)$. Moreover, given a neighborhood $\mcal V$ of $[a_1,b_1]$, after decreasing $\delta_0$ if necessary, the curves $\Gamma_{\bm S}^\pm$ can be chosen inside $\mcal V$.
    \end{lemma}

    Introducing curves $\Gamma_{\bm S}^\pm$ chosen as in Lemma~\ref{lem:lens}, we define
    \[
        \bm S(z) = \bm T(z) \times
        \begin{dcases}
            \bm I, \quad &\text{outside the lenses} \\
            \bm I \mp \lambda_j(z)^{-1} \ee^{4 N \phi_{b_1}(z)} \bm E_{21}, \quad &z \in \Omega_{\bm S}^{\pm}
        \end{dcases}
    \]

    $\bm S$ satisfies the following RHP.

    \begin{rhp}\label{rhp:S} Find a matrix-valued function $\bm S: \C \setminus \Gamma_{\bm S} \to \C^{2 \times 2}$, $\Gamma_{\bm S} \defeq \R \cup \Gamma_{\bm S}^+ \cup \Gamma_{\bm S}^-$, such that
        \begin{enumerate}[(I)]
            \item $\bm S$ is analytic.
            \item $\bm S$ has boundary values given by
            \[
                \bm S_+(z) = \bm S_-(z) \bm J_{\bm S}(z), \quad  \bm J_{\bm S}(z) = \begin{dcases}
                \bm I + \lambda_j(z) \ee^{-2N F(z)} \bm E_{12}, \quad &z \in \R \setminus [a_1,b_1], \\
                \lambda_j(z) \bm E_{12} - \lambda_j(z)^{-1} \bm E_{21}, \quad &z \in (a_1,b_1), \\
                \bm I + \lambda_j(z)^{-1} \ee^{4N \phi_{b_1}(z)} \bm E_{21}, \quad &z \in \Gamma_{\bm S}^+ \cup \Gamma_{\bm S}^-,
                \end{dcases}
            \]
            with the lips oriented from $a_1$ to $b_1$, and $\bm J_{\bm S}(0) = \bm I$.
            \item As $z \to \infty,$
            \[
                \bm S(z) = \left( \bm I + \frac{\bm T_1}{z} + \Boh\left( \frac{1}{z^2} \right) \right) z^{k \sp_3}.
            \]
        \end{enumerate}
    \end{rhp}

\subsubsection{Global parametrix} \hfill

    By dropping all jump matrices that go to $\bm I$ as $N \to \infty$ in the RHP for $\bm S$ we are led to the global parametrix RHP.

    \begin{rhp}\label{rhp:G}
        Find the analytic matrix-valued function $\bm G: \C \setminus [a_1,b_1] \to \C^{2 \times 2}$ satisfying
        \begin{enumerate}[(i)]
            \item $\bm G$ has boundary values on $(a_1,b_1)$ given by
            \[
                \bm G_+(x) = \bm G_-(x) \left(\lambda_j(x) \bm E_{12} - \lambda_j(x)^{-1} \bm E_{21} \right), \quad x \in (a_1,b_1).
            \]
            \item As $z \to \infty$,
            \[
                \bm G(z) = \left(\bm I + \frac{\bm G_1}{z} + \Boh \left( \frac{1}{z^{2}} \right) \right) z^{k \sp_3}.
            \]
            \item $\bm G$ has, at worst, fourth root singularities at $z = a_1, b_1$.
        \end{enumerate}
    \end{rhp}

    For the construction of $\bm G$ we require Szeg\H{o} functions (see for instance \cite[Section~5]{Kuijlaars2004a}), which deal with the jump $\lambda_j$, and a normalization that captures the $z^{k \sp_3}$ behavior at infinity. The normalization is given by the inverse Joukowski map
    \[
        \varphi(z)
            \defeq \frac{z - m + \msf h(z)}{r}, \qquad \varphi(z) = \frac{2z}{r} + \Boh(1), \quad z \to \infty,
    \]
    with $m$, $r$ and $\msf h$ as in \eqref{eq:m-r-u}.  $\varphi$ maps $\C \setminus [a_1,b_1]$ conformally onto $|w| > 1$ and satisfies
    \[
        \varphi'(z) = \frac{\varphi(z)}{\msf h(z)}, \quad \text{in particular} \quad \varphi'(0) = \frac{\varphi(0)}{u},
        \qquad \frac{1}{\varphi(z)} = \frac{z-m-\msf h(z)}{r}.
    \]
    The Szeg\H{o} function of $x^2$ that we consider is given by
    \begin{equation}\label{eq:szego-def}
        D(z) \defeq
            c \varphi(0) \left(1 - \frac{1}{\varphi(z) \varphi(0)} \right)^2 = \frac{2 \varphi(0)z^2}{r (\varphi(z) - \varphi(0))^2}, \quad
        c
            \defeq \frac{r}{2} =  \frac{b_1 - a_1}{4}.
    \end{equation}
    It is analytic  in $\C \setminus [a_1,b_1]$ and  satisfies
    \[
        D_+(x) D_-(x) = x^2, \quad x \in (a_1,b_1), \quad
        D(z) \neq 0, \quad \forall z \in \C \setminus [a_1,b_1], \quad D(\infty) = c \varphi(0).
    \]
    The Szeg\H{o} function of $\lambda$ is
    \begin{equation}\label{eq:szego-b}
        D_{\lambda}(z) \defeq \exp\left(\frac{\msf h(z)}{2\pi}\int_{a_1}^{b_1}\frac{\log\lambda(s)}{\sqrt{(s - a_1)(b_1 - s)}}\,\frac{\dd s}{z - s}\right), \qquad z \in \C\setminus[a_1,b_1].
    \end{equation}
    Since $\lambda$ is real analytic and positive on $[a_1,b_1]$, the function $D_{\lambda}$ is analytic and does not vanish in $\C\setminus[a_1,b_1]$, it extends continuously and without zeros to $[a_1,b_1]$ from both sides, it is positive on $\R\setminus[a_1,b_1]$, and
    \[
        D_{\lambda,+}(x)\,D_{\lambda,-}(x) = \lambda(x), \quad x \in (a_1,b_1), \qquad
        D_{\lambda}(\infty) = \exp\left(\frac{1}{2\pi}\int_{a_1}^{b_1}\frac{\log\lambda(s)\,\dd s}{\sqrt{(s - a_1)(b_1 - s)}}\right) > 0.
    \]
    In particular, $|D_{\lambda}|$ is bounded above and below by positive constants on $\C\setminus[a_1,b_1]$. These properties determine $D_{\lambda}$ uniquely, and $D^jD_{\lambda}$ is the Szeg\H{o} function of $\lambda_j = x^{2j}\lambda$.

    From the jump relations one checks that
    \begin{equation}\label{eq:G-def}
        \bm G(z) = \left(D(\infty)^{j}D_{\lambda}(\infty)\right)^{\sp_3} c^{k \sp_3}\, \bm M(z)\, \varphi(z)^{k \sp_3} \left(D(z)^{j}D_{\lambda}(z)\right)^{-\sp_3},
    \end{equation}
    where
    \[
        \bm M(z) \defeq \bm U_0 \gamma(z)^{\sp_3} \bm U_0^{-1}, \qquad \bm U_0 \defeq \frac{1}{\sqrt{2}} \begin{pmatrix}
            1   &   \ii\\
            \ii &   1
        \end{pmatrix}, \qquad \gamma(z) \defeq \left( \frac{z-b_1}{z-a_1}\right)^{\frac{1}{4}} \sim 1 \quad (z \to \infty),
    \]
    is the unique solution of RHP~\ref{rhp:G}.

    We are interested in the first columns of $\bm G^{[0,1]}$ and $\bm G^{[1,0]}$, which determine the asymptotics of  $p$ (through \eqref{eq:p-via-Y-rescaled}) and $q$ (through \eqref{eq:q-via-Y-rescaled}).  The following lemma compiles all the relations we will require from $\bm G$ after the asymptotic analysis.

    \begin{lemma}\label{lem:first-columns}
        On $\C \setminus [a_1,b_1]$ the following hold.
        \begin{enumerate}[\rm (i)]
            \item $\bm M_{11} = \frac{1}{2}\left( \gamma + \gamma^{-1}\right)$ and $\bm M_{21} = \frac{\ii}{2}\left( \gamma - \gamma^{-1}\right) = -\ii\, \bm M_{11}/\varphi$ do not vanish and $\bm M_{11}(0) > 0$.
            \item  The expressions $\dfrac{D(\infty)}{D(z)} = \left( \dfrac{r}{2} \dfrac{ \varphi(z) - \varphi(0)}{z}\right)^2$, $\dfrac{\varphi(z) - \varphi(0)}{z}$ and $D_{\lambda}(z)$ do not vanish.
            \item The first columns of $\bm G^{[j,k]}$ for $(j,k) = (0,1)$ and $(1,0)$ are
            \[
                \begin{split}
                \bm G^{[0,1]}(x)\bm e_1
                    &= \frac{\bm M_{11}(x)}{D_{\lambda}(x)} \begin{pmatrix} c\,D_{\lambda}(\infty)\, \varphi(x) \\ -\ii/\left(c\,D_{\lambda}(\infty)\right) \end{pmatrix}, \qquad \\
                \bm G^{[1,0]}(x)\bm e_1
                    &=
                    \frac{D_{\lambda}(\infty)}{D_{\lambda}(x)}
                    \,  \bm M_{11}(x) \left( \dfrac{r}{2} \dfrac{ \varphi(x) - \varphi(0)}{x}\right)^2 \begin{pmatrix} 1 \\ -\ii\, \left(D(\infty)D_{\lambda}(\infty)\right)^{-2} \varphi(x)^{-1}\end{pmatrix}.
                \end{split}
            \]
            \item For $\bm G = \bm G^{[0,1]}$,
            \[
                \bm G_{11}(0)\, \bm G_{21}(z) - \bm G_{21}(0)\, \bm G_{11}(z) = \frac{\ii\, \bm M_{11}(0)\, \bm M_{11}(z) \left( \varphi(z) - \varphi(0)\right)}{D_{\lambda}(0)\,D_{\lambda}(z)},
            \]
            and its derivative at $z=0$ equals $\dfrac{\ii\, \bm M_{11}(0)^2 \varphi(0)}{u\,D_{\lambda}(0)^2} = \dfrac{\ii\, r\, \varphi(0)^2}{2u^2D_{\lambda}(0)^2} \neq 0$.
        \end{enumerate}
    \end{lemma}
    \needspace{4\baselineskip}
    \begin{proof} \hfill

        (i) The expressions for $\bm M_{11}$ and $\bm M_{21}$ come from the definition of $\bm M$.  Let $w \defeq \gamma^2 = \frac{\sqrt{z-b_1}}{\sqrt{z-a_1}}$. Then expanding in polar coordinates it is easy to check that $\Re w > 0$, and in particular $w \neq -1$. Writing $\msf h(z) = \sqrt{z-a_1}\sqrt{z-b_1}$  and using $(z - m)^2 - \msf h(z)^2 = r^2$, we obtain
        \[
            \frac{w - 1}{w + 1} = \frac{\left(\sqrt{z-b_1} - \sqrt{z-a_1}\right)^2}{(z - b_1) - (z - a_1)} = \frac{2(z-m) - 2\msf h(z)}{a_1 - b_1} = -\frac{z - m - \msf h(z)}{r} = - \frac{1}{\varphi(z)}.
        \]
        Therefore $\bm M_{21} = \frac{\ii}{2}\left(\gamma + \gamma^{-1}\right) \frac{\gamma^2 - 1}{\gamma^2 + 1} = -\ii\, \bm M_{11}/\varphi$. Finally, $\bm M_{11} = 0$ implies $w = \gamma^2 = -1$, a contradiction. Moreover, $\gamma(0) = (b_1/a_1)^{1/4} > 0$, so that $\bm M_{11}(0) = \frac12\left(\gamma(0) + \gamma(0)^{-1}\right) = \sqrt{\frac{u - m}{2u}}$.

        (ii) The formula for $\frac{D(\infty)}{D}$ follows from the second representation for $D$ in \eqref{eq:szego-def}. As $\varphi$ is injective on $\C\setminus[a_1,b_1]$ and $0 \notin [a_1,b_1]$, $\varphi(z) \neq \varphi(0)$ for $z \neq 0$. At $z = 0$ the quotient tends to $\varphi'(0) = \varphi(0)/u \neq 0$. Finally, $D_{\lambda}$ does not vanish by \eqref{eq:szego-b}.

        (iii) For every $2\times2$ matrix $\bm K$ we have $X^{\sp_3} \bm K\, Y^{\sp_3} \bm e_1 = (X \bm K_{11} Y, X^{-1} \bm K_{21} Y)^T$. The result follows from the expressions of $\bm M_{11}$ and $\bm M_{21}$ obtained in (i) and the formula for $\bm G$ for the chosen $[j,k]$.

        (iv) The identity follows from (iii), noticing that the factors $D_{\lambda}(\infty)$ cancel. Its derivative at $0$ follows by noticing that the factor $\varphi(z) - \varphi(0)$ vanishes at $0$, and using that $\varphi'(0) = \frac{\varphi(0)}{u}$. Finally, we have $w(0) = \sqrt{b_1/a_1}  > 0$, so that
        \[
            \bm M_{11}(0)^2 = \frac{1}{4}\left( w(0) + 2 + w(0)^{-1}\right), \qquad w(0) + w(0)^{-1} = \frac{|a_1| + |b_1|}{u} = -\frac{2m}{u}.
        \]
        Hence $\bm M_{11}(0)^2 = \frac{u - m}{2u}$, and $\varphi(0) = \frac{u - m}{r}$ gives the second expression for the derivative at $z=0$.
    \end{proof}

\subsubsection{Local parametrices}\label{sec:local-param} \hfill

    Let $\delta$ be as in Lemma~\ref{lem:lens}. The local parametrices $\bm P_{a_1}$ and $\bm P_{b_1}$ in the disks $U_{a_1,\delta}$ and $U_{b_1,\delta}$ are constructed in the standard way from the Airy parametrix of \cite[Section~7]{deift_exp_weights}, using the conformal maps $\zeta_{a_1}$ and $\zeta_{b_1}$ of Lemma~\ref{lem:zeta-maps} and the branch of $\lambda_j^{1/2}$ that is analytic in $U_{a_1,\delta}\cup U_{b_1,\delta} \subset \mcal V_{\lambda}$ and positive on the real line. For $p \in \{a_1,b_1\}$, $\bm P_p$ is analytic in $U_{p,\delta}\setminus\Gamma_{\bm S}$, it has the same jumps as $\bm S$ on $U_{p,\delta}\cap\Gamma_{\bm S}$, it remains bounded as $z \to p$, $\det\bm P_p \equiv 1$, and
    \begin{equation}\label{eq:matching}
        \bm P_p(z) = \left( \bm I + \Boh\left(\frac{1}{N} \right)\right) \bm G(z), \qquad N \to \infty,
    \end{equation}
    uniformly for $z \in \partial U_{p,\delta}$.

\subsubsection{Final transformation and error analysis} \hfill

    Let
    \[
        \bm P(z) \defeq
        \begin{dcases}
            \bm P_{a_1}(z), \quad &z \in U_{a_1,\delta}, \\
            \bm P_{b_1}(z), \quad &z \in U_{b_1,\delta},
        \end{dcases}
    \]
    and consider
    \[
        \bm R(z) \defeq
        \begin{dcases}
            \bm S(z) \bm G(z)^{-1}, \quad
                &z \in \C \setminus \left( \Gamma_{\bm S} \cup \overline{U_{a_1,\delta}} \cup \overline{U_{b_1,\delta}}\right), \\
            \bm S(z) \bm P(z)^{-1}, \quad
                &z \in \left( U_{a_1,\delta} \cup U_{b_1, \delta}\right) \setminus \Gamma_{\bm S}.
        \end{dcases}
    \]
    Both $\bm S$ and $\bm G$ have the same jump on $(a_1,b_1)$, and $\bm S$ and $\bm P$ have the same jumps inside the disks (Section~\ref{sec:local-param}), so $\bm R$ is analytic across those contours. Its isolated singularities at $a_1$ and $b_1$ are removable. Let
    \[
        \Gamma_{\bm R} \defeq \partial U_{a_1,\delta} \cup \partial U_{b_1,\delta}
            \cup \left( \left( \Gamma_{\bm S}^+ \cup \Gamma_{\bm S}^- \cup \left( \R \setminus [a_1,b_1] \right) \right) \setminus
            \left( U_{a_1,\delta}\cup U_{b_1,\delta}\right)\right),
    \]
    where the circles are oriented clockwise and the other contours as before. $\bm R$ solves the following RHP.

    \needspace{8\baselineskip}
    \begin{rhp}\label{rhp:R}
        Find $\bm R : \C \setminus \Gamma_{\bm R} \to \C^{2\times 2}$ such that
        \begin{enumerate}[\rm (i)]
            \item $\bm R$ is analytic;
            \item $\bm R$ admits continuous boundary values along $\Gamma_{\bm R}$ related by $\bm R_+(z) = \bm R_-(z) \bm J_{\bm R}(z)$, where
            \[
                \bm J_{\bm R}(z) =
                    \begin{dcases}
                        \bm P(z) \bm G(z)^{-1}, \quad &z \in \partial U_{a_1,\delta} \cup \partial U_{b_1,\delta}, \\
                        \bm I + \lambda_j(z)\, \ee^{- 2 N F(z)}\, \bm G(z) \bm E_{12} \bm G(z)^{-1},
                            \quad &z \in \R \setminus [a_1 - \delta, b_1 + \delta], \\
                        \bm I + \lambda_j(z)^{-1} \ee^{4N \phi_{b_1}(z)}\, \bm G(z) \bm E_{21} \bm G(z)^{-1},
                            \quad &z \in \left( \Gamma_{\bm S}^+ \cup \Gamma_{\bm S}^- \right) \setminus \left( U_{a_1,\delta}\cup U_{b_1,\delta}\right)
                    \end{dcases}
            \]
            with $\bm J_{\bm R}(0) = \bm I$;
            \item as $z \to \infty$,
            \[
                \bm R(z) = \bm I + \frac{\bm R_1}{z} + \Boh \left(\frac{1}{z^2} \right).
            \]
        \end{enumerate}
    \end{rhp}

    \begin{lemma}\label{lem:small-norm}
        There exists $\eta > 0$ such that, as $N \to \infty$,
        \[
            \left\| \bm J_{\bm R} - \bm I \right\|_{L^1 \cap L^\infty(\partial U_{a_1,\delta} \cup \partial U_{b_1,\delta})}
                = \Boh \left( \frac{1}{N}\right), \quad
            \left\| \bm J_{\bm R} - \bm I \right\|_{L^1 \cap L^\infty \left(\Gamma_{\bm R} \setminus
            \left(\partial U_{a_1,\delta} \cup \partial U_{b_1,\delta}\right) \right)}
                = \Boh \left( \ee^{- \eta N} \right).
        \]
    \end{lemma}

    \begin{proof}
        On the circles the estimate is the matching condition \eqref{eq:matching}. On $\R \setminus [a_1 - \delta, b_1 + \delta]$ it comes from the exponential decay of $|\lambda_j(z)|\, \ee^{-2NF(z)}$: by (R1) and \eqref{eq:F-def}, $F$ is positive on this set and grows like $\frac12\msf Q$ as $z \to 0$ and as $|z| \to \infty$, whereas $\lambda_j$ is a Laurent polynomial and $\bm G$ grows at most polynomially at infinity. On the lenses outside the disks it is a consequence of Lemma~\ref{lem:lens}, since $\lambda_j^{-1}$ and $\bm G$ are bounded there.
    \end{proof}

    By the standard small norm theory for RHPs, Lemma~\ref{lem:small-norm} implies that, for $N$ large, $\|\bm R_\pm - \bm I\|_{L^2(\Gamma_{\bm R})} = \Boh(N^{-1})$ and $\bm R = \bm I + \Boh(N^{-1})$ uniformly at a fixed positive distance from $\Gamma_{\bm R}$. This bound extends up to $\Gamma_{\bm R}$ by deforming the contour wherever the jump matrix extends analytically to a neighborhood of $\Gamma_{\bm R}$ and satisfies there the estimates of Lemma~\ref{lem:small-norm}. This is the case on the lips of the lenses and on $\left(\R\setminus[a_1 - \delta, b_1 + \delta]\right)\cap\{|z| \geq \rho_0\}$, for every $\rho_0 > 0$. Near the circles, the bound follows by changing the radius $\delta$ slightly. Near $z = 0$, however, the contour cannot be properly deformed, since $\Re F(z)$ behaves like $|z|^{-d_-}\cos(d_-\arg z)$ as $z \to 0$. We obtain the following.
    \begin{theorem}\label{thm:small-norm}
        For every $\rho_0 > 0$, as $N \to \infty$,
        \[
            \bm R(z) = \bm I + \Boh\left( \frac{1}{N}\right) \quad \text{uniformly for } z \in \C\setminus\Gamma_{\bm R} \text{ with } |z| \geq \rho_0, \qquad
            \left\| \bm R_{\pm} - \bm I \right\|_{L^2(\Gamma_{\bm R})} = \Boh\left( \frac{1}{N}\right).
        \]
    \end{theorem}

    Theorem~\ref{thm:small-norm} does not control $\bm R$ near $z = 0$. This is not needed: in Section~\ref{sec:asympt1} we only use the first column of $\bm T$, which is analytic at $z = 0$, and we estimate it there by the maximum principle.

    This finishes the asymptotic analysis. We now turn to the proofs of the main results.

\section{Proof of main results}\label{sec:asympt1}

    In this section $\bm P_N(\zeta)$ is the $N$-th monic MVOP for the varying weight $\bm W_{N\bm Q} = \ee^{-N\bm Q}\bm B$, and $f=f_N$, $g=g_N$, $p=p_N$, $q=q_N$ are the functions attached to it in Section~\ref{sec:prelim}, with $\bm Q$ replaced by $N\bm Q$.

    We put $\wt a_1 \defeq \iota(a_1)$ and $\wt b_1 \defeq \iota(b_1)$, so that $\iota([a_1,b_1]) = [\wt a_1,\wt b_1] \subset (0,+\infty)$, and
    \begin{equation}\label{eq:f-g-rescaled}
        f(x) = \frac{p(x)}{x^N} = \bm P_{11}(\zeta) + x\,\bm P_{12}(\zeta), \qquad g(x) = \frac{q(x)}{x^{N-1}} = \bm P_{21}(\zeta) + x\,\bm P_{22}(\zeta), \qquad \zeta = \msf z(x).
    \end{equation}
    The branch points of $\mcal R$ are $\zeta = \pm\frac{2\ii}{\msf a}$, that is, $x = \pm\ii$.

    For $\eta > 0$ let
    \[
        \mcal V_\eta \defeq \left\{x \in \C : \operatorname{dist}\left(x,[a_1,b_1]\right) < \eta\right\}.
    \]

\subsection{Proof of Theorem~\ref{thm:intro-asymptotics-fg}}\label{sec:asymptotics-pq} \hfill

    With $\bm M$, $\varphi$, $m$, $r$, $u = \sqrt{a_1b_1}$ and $D_{\lambda}$ as in Section~\ref{sec:onecut}, define, for $x \in \C\setminus[a_1,b_1]$,
    \begin{equation}\label{eq:amplitudes}
        \begin{split}
        A_f(x) &\defeq \frac{D_{\lambda}(0)}{D_{\lambda}(x)}\,\frac{\bm M_{11}(x)}{\bm M_{11}(0)}\,\frac{u}{\varphi(0)}\,\frac{\varphi(x) - \varphi(0)}{x}, \\
        A_g(x) &\defeq \frac{D_{\lambda}(\infty)}{D_{\lambda}(x)}\,\bm M_{11}(x)\left( \frac{r}{2}\,\frac{\varphi(x) - \varphi(0)}{x}\right)^2 = \bm G^{[1,0]}_{11}(x),
        \end{split}
    \end{equation}
    with the values at $x = 0$ given by continuity. A formula for $A_f$ in terms of $\bm G$ is presented in \eqref{eq:B-G} below. By Lemma~\ref{lem:first-columns}, $A_f$ and $A_g$ are analytic and do not vanish in $\C\setminus[a_1,b_1]$. Moreover, they have finite non-zero limits at infinity, and
    \begin{equation}\label{eq:amplitudes-special}
        A_f(0) = A_g(\infty) = 1.
    \end{equation}

    By \eqref{eq:f-g-rescaled}, $f(x) = x^{-N}p(x)$ and $g(x) = x^{-(N-1)}q(x)$. Hence Theorem~\ref{thm:intro-asymptotics-fg} follows from the following asymptotic formulas for the polynomials $p$ and $q$: for every $\eta > 0$, as $N \to \infty$,
    \begin{align}
        p(x) &= \left(-\frac{1}{\msf a}\right)^{N}\ee^{2N(\msf g(x) - \msf g(0))}\,A_f(x)\left(1 + \Boh\left(\frac1N\right)\right), \label{eq:p-asymptotics}\\
        q(x) &= \msf a^{-N}\ee^{2N\msf g(x)}\,A_g(x)\left(1 + \Boh\left(\frac1N\right)\right), \label{eq:q-asymptotics}
    \end{align}
    uniformly for $x \in \C\setminus\mcal V_\eta$. Indeed, a compact set $K \subset \C\setminus\left([a_1,b_1]\cup\{0\}\right)$ is at positive distance from $[a_1,b_1]$, so that $K \subset \C\setminus\mcal V_\eta$ for $\eta > 0$ small enough, and since $0 \notin K$ we may divide \eqref{eq:p-asymptotics} by $x^N$ and \eqref{eq:q-asymptotics} by $x^{N-1}$ on $K$. The same division shows that the formulas of Theorem~\ref{thm:intro-asymptotics-fg} hold uniformly for $x \in \C\setminus\left(\mcal V_\eta\cup\{0\}\right)$, for every $\eta > 0$. This stronger form is used in Section~\ref{sec:det}.

    \begin{proof}[Proof of Theorem~\ref{thm:intro-asymptotics-fg}]
        We prove \eqref{eq:p-asymptotics} and \eqref{eq:q-asymptotics}. In this proof $\|\cdot\|$ denotes the Euclidean norm of $\C^2$ and the induced matrix norm, and $C$ denotes a constant that depends on $\eta$ but not on $N$ or $x$, and may change from one occurrence to the next.

        We may assume $\eta < |b_1|$, so that $0 \notin \overline{\mcal V_\eta}$, and we fix $\rho_0 \in (0,|b_1| - \eta)$, so that the closed disk $\{|x| \leq \rho_0\}$ is contained in $\C\setminus\mcal V_\eta$. In the construction of Section~\ref{sec:onecut}, the disks and the lenses may be chosen inside $\mcal V_\eta\cap\mcal V_{\lambda}$. We fix such a choice and therefore
        \[
            \bm T(x) = \bm R(x)\,\bm G(x), \qquad x \in \C\setminus\left(\mcal V_\eta\cup\R\right),
        \]
        and by Theorem~\ref{thm:small-norm}, $\|\bm R(x) - \bm I\| \leq C/N$ for all such $x$ with $|x| \geq \rho_0$.

        By the definition of $\bm T$, its first column is
        \begin{equation}\label{eq:T-first-column}
            \bm T\bm e_1 = \left(\bm Y_{11}\,\ee^{-2N\msf g},\ \bm Y_{21}\,\ee^{-2N(\msf g + \ell)}\right)^T,
        \end{equation}
        which is an analytic function on $\C\setminus[a_1,b_1]$. Therefore, by the way we chose the contours and by continuity,
        \begin{equation}\label{eq:T-G-first-column}
            \left\|\bm T(x)\bm e_1 - \bm G(x)\bm e_1\right\| \leq \frac CN\left\|\bm G(x)\bm e_1\right\|, \quad x \in \C \setminus \mcal V_\eta,\ |x| \geq \rho_0.
        \end{equation}
        By Lemma~\ref{lem:first-columns}(i)--(iii), $\bm G\bm e_1$ is continuous and does not vanish on $\C\setminus[a_1,b_1]$, so that $\|\bm G\bm e_1\|$ is bounded above and below by positive constants on the disk $\{|x| \leq \rho_0\}$. The maximum principle on this disk, applied to each entry of $\bm T\bm e_1 - \bm G\bm e_1$, shows that \eqref{eq:T-G-first-column} holds, with a larger $C$, for all $x \in \C\setminus\mcal V_\eta$.

        To prove \eqref{eq:q-asymptotics} fix $(j,k) = (1,0)$. Then by \eqref{eq:q-via-Y-rescaled} and \eqref{eq:T-first-column}
        \begin{equation}\label{eq:q-asympt}
            q(x) = \msf a^{-N}\,\bm Y^{[1,0]}_{11}(x) = \msf a^{-N}\ee^{2N\msf g(x)}\,\bm T_{11}(x).
        \end{equation}

        Lemma~\ref{lem:first-columns}(ii)--(iii) and \eqref{eq:amplitudes} give
        \begin{equation}\label{eq:G10-first-column}
            \bm G^{[1,0]}(x)\bm e_1 = A_g(x)\left(1,\,-\ii \left(D(\infty)D_{\lambda}(\infty)\right)^{-2}\varphi(x)^{-1}\right)^T,
        \end{equation}
        and since $|\varphi| > 1$ on $\C\setminus[a_1,b_1]$, \eqref{eq:T-G-first-column} gives
        \begin{equation}\label{eq:G10-estimate}
            \left| \bm T_{11}(x) - \bm G^{[1,0]}_{11}(x) \right| \leq \frac CN\left\|\bm G^{[1,0]}(x)\bm e_1\right\| \leq \frac{ C \left(1 + \left(D(\infty)D_{\lambda}(\infty)\right)^{-2}\right)}{N}|A_g(x)|, \qquad x \in \C\setminus\mcal V_\eta.
        \end{equation}
        Thus \eqref{eq:q-asymptotics} is verified by plugging \eqref{eq:G10-first-column} and \eqref{eq:G10-estimate} into \eqref{eq:q-asympt}.

        For the proof of \eqref{eq:p-asymptotics}, fix $(j,k) = (0,1)$. For a $2 \times 2$ matrix-valued function $\bm X$ whose first column $\bm X\bm e_1$ is analytic near $0$, as is the case for $\bm Y$, $\bm G$ and, by \eqref{eq:T-first-column}, $\bm T$, we introduce the notation
        \begin{equation}\label{eq:B-def}
            \mcal B_{\bm X}(x) \defeq \bm X_{11}(0)\,\bm X_{21}(x) - \bm X_{21}(0)\,\bm X_{11}(x)
            = \det\left(\bm X(0)\bm e_1,\ \bm X(x)\bm e_1\right).
        \end{equation}
        Note that $\mcal B_{\bm X}(0) = 0$ for every such $\bm X$.
        By \eqref{eq:T-first-column},
        \begin{equation}\label{eq:BY-BT}
            \mcal B_{\bm Y}(x) = \ee^{2N(\msf g(x) + \msf g(0) + \ell)}\,\mcal B_{\bm T}(x).
        \end{equation}
        Differentiating \eqref{eq:BY-BT} at $x=0$ and using $\mcal B_{\bm T}(0)=0$, the contribution of the exponential factor vanishes, so that
        \begin{equation}\label{eq:BYp-BTp}
            \mcal B_{\bm Y}'(0) = \ee^{2N(2\msf g(0) + \ell)}\,\mcal B_{\bm T}'(0).
        \end{equation}
        Substituting \eqref{eq:BY-BT} and \eqref{eq:BYp-BTp} into \eqref{eq:p-via-Y-rescaled} yields
        \begin{equation}\label{eq:p-via-T}
            \begin{split}
                p(x)
                &= \left(- \msf a^{-1} \right)^N \frac{\mcal B_{\bm Y}(x)}{x \mcal B_{\bm Y}'(0)} \\
                &= \left(- \msf a^{-1}\right)^N\ee^{2N(\msf g(x) - \msf g(0))}\,
                   \frac{\mcal B_{\bm T}(x)}{x\,\mcal B_{\bm T}'(0)},
                \qquad x \in \C\setminus[a_1,b_1],
            \end{split}
        \end{equation}
        where $\mcal B_{\bm T}'(0) \neq 0$ because $\mcal B_{\bm Y}'(0) = 2\pi\ii\,K_{2N+1}(0,0) \neq 0$.

        On the other hand, by Lemma~\ref{lem:first-columns}(iv) and \eqref{eq:amplitudes},
        \begin{equation}\label{eq:B-G}
            \mcal B_{\bm G}(x) = \frac{\ii\,\bm M_{11}(0)\,\bm M_{11}(x)\left(\varphi(x) - \varphi(0)\right)}{D_{\lambda}(0)\,D_{\lambda}(x)},
            \qquad
            \frac{\mcal B_{\bm G}(x)}{x\,\mcal B_{\bm G}'(0)} = A_f(x).
        \end{equation}
        Thus, by \eqref{eq:p-via-T} and \eqref{eq:B-G}, the proof reduces to showing that
        \begin{equation}\label{eq:goal-f}
            \frac{\mcal B_{\bm T}(x)}{x} = \frac{\mcal B_{\bm G}(x)}{x}\left(1 + \Boh\left(\frac1N\right)\right)
        \end{equation}
        uniformly for $x \in \C\setminus\mcal V_\eta$, where at $x = 0$ both quotients are understood as derivatives. To this end, let
        \[
            e_N(x) \defeq \frac{\mcal B_{\bm T}(x) - \mcal B_{\bm G}(x)}{x},
        \]
        which is analytic in $\C\setminus[a_1,b_1]$, since $\mcal B_{\bm T}(0) = \mcal B_{\bm G}(0) = 0$ by \eqref{eq:B-def}. We need to prove that
        \[
            e_N(x) = A_f(x)\, \Boh \left( \frac{1}{N}\right)
        \]
        uniformly for $x \in \C \setminus \mcal V_\eta$. We begin with two elementary bounds. First, since $|D_\lambda|$ is bounded above on $\C\setminus\mcal V_\eta$, \eqref{eq:B-G} gives
        \[
            \left| \bm M_{11}(x) \right| \frac{|\varphi(x) - \varphi(0)|}{|x|}
            \leq C \left| \frac{\mcal B_{\bm G}(x)}{x} \right|,
            \qquad x \in \C\setminus\mcal V_\eta.
        \]
        Second, introduce
        \[
            \begin{split}
                \bm v(x) &\defeq \bm G(x)\bm e_1
                = \frac{\bm M_{11}(x)}{D_{\lambda}(x)}\left(c\,D_{\lambda}(\infty)\,\varphi(x),\,-\ii/\left(c\,D_{\lambda}(\infty)\right)\right)^T, \\
                \bm w(x) &\defeq \bm T(x)\bm e_1 - \bm v(x),
            \end{split}
        \]
        so that $\|\bm w(x)\| \leq \frac CN\|\bm v(x)\|$ on $\C\setminus\mcal V_\eta$ by \eqref{eq:T-G-first-column}. By the bilinearity of the determinant,
        \begin{equation}\label{eq:B-difference}
            \begin{split}
                \left|\mcal B_{\bm T}(x) - \mcal B_{\bm G}(x) \right|
                &= \left|\det\left(\bm w(0),\bm w(x)\right) + \det\left(\bm v(0),\bm w(x)\right)
                   + \det\left(\bm w(0),\bm v(x)\right) \right| \\
                   &\leq \frac CN\,\|\bm v(0)\|\,\|\bm v(x)\|
                \leq \frac CN\left|\bm M_{11}(x)\right|\left(1 + |\varphi(x)|\right),
                \qquad x \in \C\setminus\mcal V_\eta,
            \end{split}
        \end{equation}
        where we used that $|D_{\lambda}|$ is bounded below by a positive constant.

        We now combine the two bounds. The function $\frac{|\varphi(x) - \varphi(0)|}{1 + |\varphi(x)|}$ is continuous and positive on the closed set $\{x \in \C\setminus\mcal V_\eta : |x| \geq \rho_0\}$, because $\varphi$ is injective, and it tends to $1$ as $x \to \infty$, hence it is bounded below there by a positive constant. By \eqref{eq:B-G} and \eqref{eq:B-difference}, $e_N$ therefore satisfies
        \[
            |e_N(x)|
            \leq \frac{C}{N} |\bm M_{11}(x)| \frac{1 + |\varphi(x)|}{|\varphi(x) - \varphi(0)|}
                 \frac{|\varphi(x) - \varphi(0)|}{|x|}
            \leq \frac{C'}{N} \left| \frac{\mcal B_{\bm G}(x)}{x}\right|,
            \qquad x \in \C \setminus \mcal V_\eta, \quad |x| \geq \rho_0.
        \]
        It remains to handle the disk $\{|x| \leq \rho_0\}$. By the maximum principle, $|e_N| \leq \frac CN\max_{|y| = \rho_0}\left|\mcal B_{\bm G}(y)/y\right|$ there. Since $\mcal B_{\bm G}(x)/x$ is continuous and does not vanish on this disk (Lemma~\ref{lem:first-columns}(i) and~(ii)), the bound $|e_N(x)| \leq \frac CN\left|\mcal B_{\bm G}(x)/x\right|$ holds for $|x| \leq \rho_0$ as well. This proves \eqref{eq:goal-f}.

        Finally, since $\ee^{\msf g}$, $A_f$ and $A_g$ do not vanish, \eqref{eq:p-asymptotics} and \eqref{eq:q-asymptotics} show that, for $N$ large, all the zeros of $p$ and of $q$ lie in $\mcal V_\eta$. Equivalently, $f$ and $g$ do not vanish in $\C\setminus(\mcal V_\eta\cup\{0\})$.
    \end{proof}

\subsection{Properties of \texorpdfstring{$\SW$}{W}}\label{sec:W} \hfill

    In this section we study the function $\SW$ given in \eqref{eq:intro-W}, using the facts from potential theory collected in Appendix~\ref{app:prelim}. We start with some identities for push-forwards by $\iota$ and $\msf z$, which we also use in the proof of Theorem~\ref{thm:det-asymptotics}.

    Let $\wt\mu \defeq \iota_*\mu$ be the push-forward of $\mu$ under $\iota$, that is, $\wt\mu(B) \defeq \mu(\iota(B))$ for Borel sets $B \subset \R\setminus\{0\}$. By Lemma~\ref{lem:quad}, $\wt\mu$ is a probability measure with the bounded density $y \mapsto \mu'(\iota y)\,y^{-2}$, and $\supp\wt\mu = \iota(\supp\mu)$ is a compact subset of $\R\setminus\{0\}$. Under Assumption~\ref{ass:onecut}, $\supp\wt\mu = [\wt a_1,\wt b_1] \subset (0,+\infty)$, and we put $\mcal K \defeq [a_1,b_1]\cup[\wt a_1,\wt b_1]$.

    The following lemma is straightforward to verify, using
    \begin{equation}\label{eq:iota-difference}
        \iota x - \iota y = \frac{x - y}{xy}, \qquad x, y \in \C\setminus\{0\},
    \end{equation}
    for (i) and (ii), the factorization $\msf z(x) - \msf z(s) = \frac{\iota s}{\msf a}\,(x - s)(\iota x - s)$ for (iii), and (iii) with $\sigma = \mu$ for (iv). Here $\msf z_*\sigma$ denotes the push-forward under $\msf z$.
    \begin{lemma}\label{lem:iota-potentials} \hfill

        \begin{enumerate}[\rm (i)]
            \item For   $x \in \C\setminus\{0\}$,
            \begin{equation}\label{eq:potential-iota}
                \U^{\wt\mu}(\iota x) = \U^\mu(x) + \log|x| + \int\log|s|\,\dd\mu(s).
            \end{equation}
            \item For every $y \in \C\setminus\left(\supp\wt\mu\cup\{0\}\right)$,
            \[
                \CT^\mu(\iota y) = y + y^2\,\CT^{\wt\mu}(y).
            \]
            \item For any probability measure $\sigma$ with compact support in $\C \setminus \{0\}$ and $x \in \C\setminus\{0\}$,
            \begin{equation}\label{eq:projected-potential}
                \U^{\msf z_*\sigma}(\msf z(x)) = \U^\sigma(x) + \U^\sigma(\iota x) - \U^\sigma(0) + \log\msf a.
            \end{equation}
            \item Let $\wh{\msf g}$ be the $\msf g$-function of $\msf z_*\mu$ defined in Theorem~\ref{thm:det-asymptotics}. Then $\ee^{\wh{\msf g}}$ extends to an analytic and non-vanishing function in $\C\setminus\msf z([a_1,b_1])$, and
            \begin{equation}\label{eq:projected-g}
                \ee^{\wh{\msf g}(\msf z(x))} = \frac{1}{\msf a}\,\ee^{\msf g(x) + \msf g(\iota x) - \msf g(0)}, \qquad x \in \C\setminus\left([a_1,b_1]\cup[\wt a_1,\wt b_1]\cup\{0\}\right).
            \end{equation}
        \end{enumerate}
    \end{lemma}

    We now prove some simple properties of $\SW$.

    \begin{lemma}\label{lem:W-properties} The function $\SW$ satisfies the following properties.
        \begin{enumerate}[\rm (i)]
            \item \[\SW(x) = \U^\mu(x) - \U^{\wt\mu}(x) + \int\log|s|\,\dd\mu(s), \quad x \in \C.\] In particular, $\SW$ is continuous on $\C$ and harmonic in $\C\setminus\left(\supp\mu\cup\supp\wt\mu\right)$.
            \item $\SW(\iota x) = -\SW(x)$ for $x \in \C\setminus\{0\}$. In particular $\SW(\pm\ii) = 0$, and $\Gamma_{\SW}\setminus\{0\}$ is invariant under $\iota$. Moreover,
            \[
                \SW(0) = \U^\mu(0) = -\int\log|s|\,\dd\mu(s) = - \SW(\infty), \quad \SW(\infty) \defeq \lim_{x \to \infty} \SW(x).
            \]
        \end{enumerate}
    \end{lemma}

    \begin{proof}
        (i) For $s \neq 0$,
        \[
            \log\left|\frac{1 + sx}{x - s}\right| = \log|s| + \log|x - \iota s| - \log|x - s|.
        \]
        Integrating against $\mu$ gives the formula. Since $\mu$ and $\wt\mu$ have bounded densities and compact supports in $\R\setminus\{0\}$, their potentials are continuous on $\C$ and harmonic off the supports.

        (ii) For $x \neq 0$,
        \[
            \frac{1 + s\,\iota x}{\iota x - s} = -\frac{x - s}{1 + sx}.
        \]
        Hence the integrand changes sign under $x \mapsto \iota x$, and therefore $\SW(\iota x) = -\SW(x)$. Since $\iota(\pm \ii) = \pm \ii$, this gives $\SW(\pm \ii) = -\SW(\pm \ii) = 0$. Finally, $\SW(\infty) = -\SW(0)$ follows from (i) and from $\U^\mu(x) - \U^{\wt\mu}(x) \to 0$ as $x \to \infty$.
    \end{proof}

    In general $\SW$ may change sign on $[a_1,b_1]$. However, it cannot be non-positive on the whole support.

    \begin{lemma}\label{lem:W-positive-somewhere}
        There is $x \in [a_1,b_1]$ such that $\SW(x) > 0$.
    \end{lemma}

    \begin{proof}
        Suppose that $\SW \leq 0$ on $[a_1,b_1]$, and let $D \defeq \{x \in \overline{\C} : \SW(x) > 0\}$, with $\SW(\infty)$ as in Lemma~\ref{lem:W-properties}(ii). Then $D$ does not intersect $[a_1,b_1]$.

        By Lemma~\ref{lem:W-properties}(i), $\SW$ is continuous on $\C$ and has a limit at $\infty$, so we may regard it as a continuous function on the Riemann sphere $\overline{\C}$ and $D$ as an open subset of $\overline{\C}$. Also by Lemma~\ref{lem:W-properties}(i), $\Delta\SW = 2\pi(\wt\mu - \mu)$ in the sense of distributions, and $\SW$ is harmonic in a neighborhood of $\infty$. Hence $\SW$ is subharmonic in $D$.

        Since $\SW$ is continuous and $D = \{\SW > 0\}$, we have $\SW = 0$ on the boundary of $D$ in $\overline{\C}$. By the maximum principle, $\SW \leq 0$ in $D$, so that $D$ is empty.

        Therefore $\SW \leq 0$ on $\C$. By Lemma~\ref{lem:W-properties}(ii) we also have $\SW \geq 0$ on $\C\setminus\{0\}$, so $\SW \equiv 0$. Then $\mu = \wt\mu$, because $\Delta\SW = 2\pi(\wt\mu - \mu)$. This is impossible, since $\supp\mu \subset (-\infty,0)$ and $\supp\wt\mu \subset (0,+\infty)$.
    \end{proof}

    Next we show that $\SW$ is regular outside $[a_1,b_1]\cup[\wt a_1,\wt b_1]$.

    \begin{lemma}[$\SW$ has no critical points]\label{lem:W-critical}
        For $x = \xi + \ii\upsilon \in \C\setminus\left([a_1,b_1]\cup[\wt a_1,\wt b_1]\right)$,
        \begin{equation}\label{eq:gradient-W}
            \partial_\xi \SW(x) - \ii\,\partial_\upsilon \SW(x) = \CT^\mu(x) - \CT^{\wt\mu}(x) = \int\frac{\iota s - s}{(s - x)(\iota s - x)}\,\dd\mu(s),
        \end{equation}
        and this function does not vanish. In particular, since its modulus is $|\nabla\SW(x)|$, the gradient of $\SW$ does not vanish.
    \end{lemma}

    \begin{proof}
        The function $\uppsi \defeq \CT^\mu - \CT^{\wt\mu}$ is analytic in $\C\setminus([a_1,b_1]\cup[\wt a_1,\wt b_1])$, with $\uppsi(\bar x) = \overline{\uppsi(x)}$. The first equality in \eqref{eq:gradient-W} follows from Lemma~\ref{lem:W-properties}(i) and \eqref{eq:gradient-potential}, applied to $\mu$ and to $\wt\mu$. The second one is by the definition of $\wt\mu$.

        We first determine the behavior of $\uppsi$ on the real line. Since $\iota s - s > 0$ for $s < 0$, the last integral in \eqref{eq:gradient-W} shows that $x^2\uppsi(x) \to \uppsi_\infty \defeq \int(\iota s - s)\,\dd\mu(s) > 0$ as $x \to \infty$. For real $x$ outside $(a_1,b_1)\cup(\wt a_1,\wt b_1)$ the integrand is real, and $(s - x)(\iota s - x)$ is positive if $x \leq a_1$ or $x \geq \wt b_1$, and negative if $b_1 \leq x \leq \wt a_1$, for $\mu$-almost every $s$. At the four endpoints the integral converges, because the density of $\mu$ vanishes as a square root at $a_1$ and $b_1$ by (R3). This determines the sign of $\uppsi$ on the real line, and in particular $\uppsi \neq 0$ there.

        It remains to show that $\uppsi$ has no zeros in the upper half plane. We do so by the argument principle.

        By \eqref{eq:Cauchy-h-q}, $\CT^\mu = \frac{\msf h\,\msf q}{\msf k\,z^{d_- + 1}} - \frac{V'}{2}$, where the poles of the two terms at $z = 0$ cancel, since $\CT^\mu$ is analytic there. Hence $\CT^\mu$ extends continuously to the closed upper half plane, with boundary values $\CT^\mu_+(x) = -\frac{V'(x)}{2} + \pi\ii\,\mu'(x)$ on $[a_1,b_1]$, where $\mu' > 0$ on $(a_1,b_1)$ by (R2). By Lemma~\ref{lem:iota-potentials}(ii), $\CT^{\wt\mu}(y) = y^{-2}\left(\CT^\mu(\iota y) - y\right)$. Since $\iota$ maps the upper half plane onto itself, and $\CT^{\wt\mu}$ is analytic at $0$, $\CT^{\wt\mu}$ also extends continuously to the closed upper half plane, with $\Im\CT^{\wt\mu}_+(y) = \pi\,y^{-2}\mu'(\iota y) > 0$ for $y \in (\wt a_1,\wt b_1)$. Therefore $\uppsi$ extends continuously to the closed upper half plane, its boundary values $\uppsi_+$ on the real line do not vanish, and
        \[
            \Im\uppsi_+ > 0 \ \text{ on } (a_1,b_1), \qquad
            \Im\uppsi_+ < 0 \ \text{ on } (\wt a_1,\wt b_1), \qquad
            \uppsi_+ = \uppsi \in \R\setminus\{0\} \ \text{ elsewhere on } \R.
        \]

        We follow a continuous branch of $\arg\uppsi_+$ along the real line, from $-\infty$ to $+\infty$. On $(-\infty,a_1]$, $\uppsi > 0$ and the argument is constant, say $0$. Along $[a_1,b_1]$, $\uppsi_+$ goes from the positive real axis to the negative one through the upper half plane, so the argument increases to $\pi$. On $[b_1,\wt a_1]$, $\uppsi < 0$. Along $[\wt a_1,\wt b_1]$, $\uppsi_+$ goes from the negative real axis to the positive one through the lower half plane, so the argument increases to $2\pi$, and then it remains constant on $[\wt b_1,+\infty)$. Thus the argument of $\uppsi_+$ increases by $2\pi$ along the real line.

        On the semicircle $x = R\ee^{\ii\vartheta}$, $0 \leq \vartheta \leq \pi$, with $R$ large, $\uppsi(x) = \uppsi_\infty x^{-2}\left(1 + \Boh(R^{-1})\right)$, so the argument of $\uppsi$ changes by $-2\pi + \Boh(R^{-1})$.

        By the argument principle, applied to the half disk $\{|x| < R,\ \Im x > 0\}$, the number of zeros of $\uppsi$ in the upper half plane is $\frac{1}{2\pi}\left(2\pi - 2\pi\right) = 0$.
    \end{proof}

    By Lemma~\ref{lem:W-critical}, outside $[a_1,b_1]\cup[\wt a_1,\wt b_1]$ the set $\Gamma_{\SW}$ is a union of analytic arcs.

    \begin{lemma}\label{lem:W-zero-set} \hfill

        \begin{enumerate}[\rm (i)]
            \item $\SW$ has finitely many zeros on $[a_1,b_1]\cup[\wt a_1,\wt b_1]$.
            \item $\Gamma_{\SW}$ and $\Gamma = \msf z(\Gamma_{\SW}\setminus\{0\})$ have zero area.
        \end{enumerate}
    \end{lemma}

    \begin{proof}
        (i) Since $F = 0$ on $[a_1,b_1]$, we have $2\U^\mu = \ell - V$ there, and hence
        \[
            \SW(x) = \frac{\ell - V(x)}{2} - \U^\mu(\iota x) + \log|x|, \qquad x \in [a_1,b_1].
        \]
        The right-hand side is real analytic on $(-\infty,0)$. By Lemma~\ref{lem:W-positive-somewhere}, this function does not vanish identically on $[a_1,b_1]$, so it has finitely many zeros there. The same holds on $[\wt a_1,\wt b_1]$ since $\SW\circ\iota = -\SW$.

        (ii) By Lemma~\ref{lem:W-properties}(i), $\SW$ is harmonic in $\C\setminus \mcal K$, and by Lemma~\ref{lem:W-critical} its gradient does not vanish there. Since $\Gamma_{\SW}$ is the zero set of $\SW$, the implicit function theorem shows that every point of $\Gamma_{\SW}\setminus \mcal K$ has a neighborhood in which $\Gamma_{\SW}$ is an analytic arc. Countably many of these neighborhoods cover $\Gamma_{\SW}\setminus \mcal K$, so $\Gamma_{\SW}\setminus \mcal K$ has zero area, and so has $\Gamma_{\SW}$, since $\mcal K \subset \R$. The same holds for $\Gamma$, because $\msf z$ is smooth in $\C\setminus\{0\}$.
    \end{proof}

\subsection{Proof of Theorem~\ref{thm:det-asymptotics}}\label{sec:det} \hfill

    Let
    \begin{equation}\label{eq:A-det}
        A_{\det}(x) \defeq \frac{x\,A_f(\iota x)\,A_g(x) - \iota x\,A_f(x)\,A_g(\iota x)}{x - \iota x}, \qquad x \in \C\setminus\left([a_1,b_1]\cup[\wt a_1,\wt b_1]\cup\{0,\ii,-\ii\}\right).
    \end{equation}
    The numerator vanishes at $x = \pm\ii$, where $\iota x = x$, so $A_{\det}$ extends analytically to $x = \pm \ii$. It also extends analytically to $x = 0$, with $A_{\det}(0) = A_f(0)A_g(\infty) = 1$ by \eqref{eq:amplitudes-special}. Note that $A_{\det}\circ\iota = A_{\det}$, and therefore $A_{\det}(x) \to 1$ as $x \to \infty$.

    Since $A_{\det}\circ\iota = A_{\det}$, and $\msf z(x') = \msf z(x)$ if, and only if, $x' \in \{x,\iota x\}$, there is a function $\wh A_{\det}$ on $\overline{\C}\setminus\msf z([a_1,b_1])$ such that
    \begin{equation}\label{eq:A-det-hat}
        A_{\det} = \wh A_{\det}\circ\msf z, \qquad \wh A_{\det}(\infty) = 1.
    \end{equation}
    It is analytic: near every point other than the branch points $\pm\frac{2\ii}{\msf a}$ and $\infty$, the map $\msf z$ has analytic local inverses, and at these three points $\wh A_{\det}$ is continuous, so that its singularities there are removable.

    The exponential factor in \eqref{eq:det-asymptotics} is the one of Theorem~\ref{thm:intro-asymptotics-fg}, combined over the two points $x$ and $\iota x$. Indeed, by Lemma~\ref{lem:iota-potentials}(iv),
    \begin{equation}\label{eq:dN}
        \msf d_N(x) \defeq \msf a^{-2N}\ee^{2N(\msf g(x) + \msf g(\iota x) - \msf g(0))} = \ee^{2N\wh{\msf g}(\msf z(x))}.
    \end{equation}

    We prove Theorem~\ref{thm:det-asymptotics} in the variable $x$. This is enough because, by \eqref{eq:A-det-hat} and \eqref{eq:dN}, $A_{\det} = \wh A_{\det}\circ\msf z$ and $\msf d_N = \ee^{2N\wh{\msf g}\circ\msf z}$, so that for $\zeta = \msf z(x)$ the formula \eqref{eq:det-asymptotics} reads
    \begin{equation}\label{eq:det-asymptotics-x}
        \det\bm P_N(\msf z(x)) = \msf d_N(x)\,A_{\det}(x)\left(1 + \Boh\left(\frac1N\right)\right).
    \end{equation}
    More precisely, for $\eta > 0$ let
    \[
        X_\eta \defeq \left\{x \in \C\setminus\{0\} : x \notin \mcal V_\eta \text{ and } \iota x \notin \mcal V_\eta\right\}
    \]
    be the set of $x \neq 0$ such that $x$ and $\iota x$ lie at distance at least $\eta$ from $[a_1,b_1]$. Let $O$ be a neighborhood of $\msf z([a_1,b_1])$, and choose $\eta \in (0,|b_1|)$ with $\msf z(\mcal V_\eta) \subset O$. Every $\zeta \in \C\setminus O$ is of the form $\zeta = \msf z(x)$ with $x \neq 0$, and since $\zeta = \msf z(x) = \msf z(\iota x)$ is not in $\msf z(\mcal V_\eta)$, neither $x$ nor $\iota x$ lies in $\mcal V_\eta$, that is, $x \in X_\eta$. Hence, to prove \eqref{eq:det-asymptotics}, it suffices to prove \eqref{eq:det-asymptotics-x} uniformly on $X_\eta$, for every $\eta > 0$. In the same way, the positivity and the non-vanishing of $\wh A_{\det}$ follow from those of $A_{\det}$.

    Recall also that, by \eqref{eq:det-fg} with $\bm Q$ replaced by $N\bm Q$,
    \[
        \det\bm P_N(\msf z(x)) = \frac{f(\iota x)g(x) - f(x)g(\iota x)}{x - \iota x}.
    \]

    \begin{proof}[Proof of Theorem~\ref{thm:det-asymptotics}]
        We first establish \eqref{eq:det-asymptotics-x} uniformly on $X_\eta$ with an absolute error, and then verify the positivity, non-vanishing, relative error and normalized zero counting measure convergence statements.

        If $x$ and $\iota x$ are outside $\mcal V_\eta$, then by \eqref{eq:f-g-rescaled}, \eqref{eq:p-asymptotics} and \eqref{eq:q-asymptotics}, which hold uniformly on $\C\setminus\mcal V_\eta$,
        \begin{align*}
            f(x)\,g(\iota x) &= \iota x\,\msf d_N(x)\,A_f(x)A_g(\iota x)
                \left(1 + \Boh\left(\frac1N\right)\right), \\
            f(\iota x)\,g(x) &= x\,\msf d_N(x)\,A_f(\iota x)A_g(x)
                \left(1 + \Boh\left(\frac1N\right)\right).
        \end{align*}
        Hence
        \[
            \frac{\det\bm P_N\left(\msf z(x)\right)}{\msf d_N(x)}
            = \frac{x\,A_f(\iota x)A_g(x)\left(1 + \Boh(N^{-1})\right)
                 - \iota x\,A_f(x)A_g(\iota x)\left(1 + \Boh(N^{-1})\right)}
                {x - \iota x}.
        \]

        We now compare the right-hand side with $A_{\det}(x)$. We may assume that $\eta < |b_1|$ and that, for some $\rho_1 > 0$, the closed disks $\{|x \mp \ii| \leq \rho_1\}$ are contained in $X_\eta$. The functions $A_f$ and $A_g$ are bounded on $\C\setminus\mcal V_\eta$, and $\frac{|x| + |\iota x|}{|x - \iota x|} = \frac{1 + |x|^2}{|1 + x^2|}$ is bounded on the sets $\{|x \mp \ii| \geq \rho_1\}$. Hence the right-hand side equals $A_{\det}(x) + \Boh(N^{-1})$ uniformly for $x \in X_\eta$ with $|x - \ii| \geq \rho_1$ and $|x + \ii| \geq \rho_1$.

        On each disk $\{|x \mp \ii| \leq \rho_1\}$ the function $\det\bm P_N(\msf z(x))/\msf d_N(x) - A_{\det}(x)$ is analytic and is $\Boh(N^{-1})$ on the boundary, so by the maximum principle it is $\Boh(N^{-1})$ on the disk. Therefore $\det\bm P_N(\msf z(x)) = \msf d_N(x)\left(A_{\det}(x) + \Boh(N^{-1})\right)$ uniformly on $X_\eta$.

        \textit{Positivity:} For real $x \notin [a_1,b_1]$ the number $\gamma(x) = \left(\frac{x - b_1}{x - a_1}\right)^{1/4}$ of Lemma~\ref{lem:first-columns} is positive, so $2\bm M_{11}(x) = \gamma(x) + \gamma(x)^{-1} > 0$. Moreover $\varphi$ is real on $\R\setminus[a_1,b_1]$, with $\varphi > 1$ on $(b_1,+\infty)$ and $\varphi < -1$ on $(-\infty,a_1)$, and $\varphi' = \varphi/\msf h > 0$ on both half-lines. Hence $\frac{\varphi(x) - \varphi(0)}{x} > 0$ for every real $x \notin [a_1,b_1]$ (for $x = 0$ this is $\varphi'(0) > 0$), and, since $D_{\lambda} > 0$ on $\R\setminus[a_1,b_1]$, we have $A_f > 0$ and $A_g > 0$ on $\R\setminus[a_1,b_1]$. Consequently, for real $x \notin [a_1,b_1]\cup[\wt a_1,\wt b_1]\cup\{0\}$ both terms $x\,A_f(\iota x)A_g(x)$ and $-\iota x\,A_f(x)A_g(\iota x)$ in the numerator of \eqref{eq:A-det} have the sign of $x$, and so does $x - \iota x = x + \frac1x$. Hence $A_{\det}(x) > 0$.

        \textit{Non-vanishing:} Suppose that $A_{\det}(x_0) = 0$. Then $x_0 \notin \R$, by the previous part. Since $A_{\det} \not\equiv 0$, by Hurwitz's theorem the functions $\det\bm P_N(\msf z(x))/\msf d_N(x)$ have zeros $x_N \to x_0$. For $N$ large $x_N \notin \R$, hence $\msf z(x_N)$ is a non-real zero of $\det\bm P_N$. This is impossible as the zeros of $\det\bm P_N$ are real, since the weight $\bm W_{N\bm Q}$ is positive definite (see, e.g., \cite{DPS}).

        \textit{Relative error:} The function $A_{\det}$ is continuous and does not vanish on the closed set $\C\setminus(\mcal V_\eta\cup\iota(\mcal V_\eta))$, which contains $X_\eta$, and $A_{\det}(x) \to 1$ as $x \to \infty$. Hence $|A_{\det}|$ is bounded below by a positive constant there, and the absolute error $\Boh(N^{-1})$ above is also a relative error. This proves \eqref{eq:det-asymptotics-x} uniformly on $X_\eta$, and hence \eqref{eq:det-asymptotics}.

        \textit{Convergence of roots:} Let $O$ be a neighborhood of $\msf z([a_1,b_1]) = \msf z([\wt a_1,\wt b_1])$, and choose $\eta \in (0,|b_1|)$ with $\msf z(\mcal V_\eta) \subset O$. As in the previous part, $|A_{\det}|$ is bounded below by a positive constant on $\C\setminus(\mcal V_\eta\cup\iota(\mcal V_\eta))$. Thus, for $N$ large, $\det\bm P_N(\msf z(x)) \neq 0$ whenever $x \notin \mcal V_\eta\cup\iota(\mcal V_\eta)\cup\{0\}$. Since every $\zeta \in \C$ is of the form $\msf z(x)$, all the zeros of $\det\bm P_N$ lie in $\msf z(\mcal V_\eta) \subset O$.

        To identify the limit, let $\sigma_N$ be the normalized zero counting measure of $\det \bm P_N$. Since $\det\bm P_N$ is monic of degree $2N$, \eqref{eq:log-poly} gives
        \[
            \frac{1}{2N}\log\left|\det\bm P_N(\zeta)\right| = -\U^{\sigma_N}(\zeta).
        \]
        For every $\zeta \in \C\setminus\msf z([a_1,b_1])$, Lemma~\ref{lem:iota-potentials}(iii) gives
        \[
            \U^{\sigma_N}(\zeta)
            = \U^{\msf z_*\mu}(\zeta)
              - \frac{1}{2N}\log\left|A_{\det}(x) + \Boh\left(N^{-1}\right)\right|
            \longrightarrow \U^{\msf z_*\mu}(\zeta), \qquad \msf z(x) = \zeta.
        \]
        We already proved that for every compact neighborhood $L$ of $\msf z([a_1,b_1])$ the measures $\sigma_N$ are supported in $L$ for $N$ large. Hence every subsequence of $(\sigma_N)$ has a further subsequence $(\sigma_{N_k})$ that converges weakly to a probability measure $\sigma$ whose support is contained in $L$, see Appendix~\ref{app:prelim}. Since $L$ is arbitrary we have $\supp\sigma \subset \msf z([a_1,b_1])$. Given $\zeta \notin \msf z([a_1,b_1])$, we may choose $L$ with $\zeta \notin L$, and then $\U^{\sigma_{N_k}}(\zeta) \to \U^\sigma(\zeta)$, again by Appendix~\ref{app:prelim}. Hence $\U^\sigma = \U^{\msf z_*\mu}$ on $\C\setminus\msf z([a_1,b_1])$, hence Lebesgue-a.e.\ in $\C$. Therefore $\sigma = -\frac{1}{2\pi}\Delta\U^\sigma = \msf z_*\mu$. Since this limit does not depend on the subsequence, $\sigma_N \to \msf z_*\mu$ weakly.
    \end{proof}

\subsection{Proof of Theorem~\ref{thm:intro-asymptotics}}\label{sec:proof-asymptotics} \hfill

    By Lemma~\ref{lem:divided-differences}, each entry of $\bm P_N(\msf z(x))$ is given by a difference of the values of $f$, or of $g$, at the two points $x$ and $\iota x$ of $\mcal R$ over $\zeta = \msf z(x)$. By \eqref{eq:p-asymptotics}, for $x \neq 0$ such that $x$ and $\iota x$ are away from $[a_1,b_1]$, the ratio $\left|\frac{f(\iota x)}{f(x)} \right|$ is of order $\ee^{2N\SW(x)}$, where
    \[
        \SW(x) = \U^\mu(x) - \U^\mu(\iota x) + \log|x| = \int\log\left|\frac{1 + sx}{x - s}\right|\dd\mu(s)
    \]
    is the function introduced in \eqref{eq:intro-W}. The same holds for $\left|\frac{g(\iota x)}{g(x)}\right|$.

    The properties of $\SW$ that we use below are proven in Section~\ref{sec:W}.

    We start by writing \eqref{eq:divided-differences}, with $\bm Q$ replaced by $N\bm Q$, as
    \begin{equation}\label{eq:PN-two-sheets}
        \bm P_N\left(\msf z(x)\right) = \frac{1}{x - \iota x}
        \begin{pmatrix}
            -\iota x\,f(x)\left(1 - \tau_{11}(x)\right) & f(x)\left(1 - \tau_{12}(x)\right) \\[1ex]
            -\iota x\,g(x)\left(1 - \tau_{21}(x)\right) & g(x)\left(1 - \tau_{22}(x)\right)
        \end{pmatrix}
    \end{equation}
    for $x \in \C\setminus\{0,\ii,-\ii\}$ with $f(x) \neq 0$ and $g(x) \neq 0$, where 
    \begin{equation}\label{eq:tau-def}
        \tau_{11}(x) \defeq \frac{x\,f(\iota x)}{\iota x\,f(x)}, \qquad \tau_{12}(x) \defeq \frac{f(\iota x)}{f(x)}, \qquad
        \tau_{21}(x) \defeq \frac{x\,g(\iota x)}{\iota x\,g(x)}, \qquad \tau_{22}(x) \defeq \frac{g(\iota x)}{g(x)},
    \end{equation}
    which are rational functions of $x$. 

    For each $\zeta \in \C\setminus\Gamma$ the two solutions of $\msf z(x) = \zeta$ do not lie on $\Gamma_{\SW}$, and since $\SW(\iota x) = -\SW(x)$, exactly one of them, $x_-(\zeta)$, satisfies $\SW(x_-(\zeta)) < 0$. The function $x_-$ is analytic in $\C\setminus\Gamma$. Indeed, the branch points $\pm\frac{2\ii}{\msf a} = \msf z(\pm\ii)$ lie on $\Gamma$, because $\SW(\pm\ii) = 0$, so near each point of $\C\setminus\Gamma$ the equation $\msf z(x) = \zeta$ has two analytic solutions. $\SW$ does not vanish on them, and by continuity $x_-$ coincides with one of them.

    The proof of Theorem~\ref{thm:intro-asymptotics} combines Theorem~\ref{thm:intro-asymptotics-fg} at the point $x_-(\zeta)$ with the following upper bound at the point $\iota x_-(\zeta)$, which may lie close to, or be an element of, the support of $\mu$.

    \begin{lemma}\label{lem:upper-bound}
        For every compact set $L \subset \C$ and every $\epsilon > 0$ there is $N_0$ such that, for $N \geq N_0$ and $y \in L$,
        \[
            |p(y)| \leq \msf a^{-N}\ee^{2N(\epsilon - \U^\mu(y) - \msf g(0))}, \qquad |q(y)| \leq \msf a^{-N}\ee^{2N(\epsilon - \U^\mu(y))}.
        \]
    \end{lemma}

    \begin{proof}
        We prove the bound for $p$. The proof for $q$ is the same, with \eqref{eq:q-asymptotics} in place of \eqref{eq:p-asymptotics}. Let
        \[
            u_N \defeq \frac{1}{2N}\log|p|, \qquad
            \psi \defeq -\U^\mu - \msf g(0) - \frac12\log\msf a.
        \]
        We have to show that $u_N \leq \psi + \epsilon$ on $L$.

        The function $u_N$ is subharmonic in $\C$, and $\psi$ is continuous in $\C$. Fix $R > 0$ such that $L\cup[a_1,b_1] \subset \{|y| < R\}$, and let $C_1$ be the maximum of $|\psi|$ on $\{|y| \leq R + 1\}$.

        Let $\eta \in (0,1)$. By \eqref{eq:p-asymptotics} there is $C_\eta$ such that, for $N$ large, $|u_N - \psi| \leq C_\eta/N$ on $\{|y| \leq R + 1\}\setminus\mcal V_\eta$. In particular $u_N \leq C_1 + 1$ on the circle $\{|y| = R + 1\}$, and by the maximum principle also on the disk $\{|y| \leq R + 1\}$, for $N$ large.

        Let $y \in L$ and $0 < r < 1$, so that $D(y,r) \subset \{|w| < R + 1\}$, and put $E \defeq D(y,r)\cap\mcal V_\eta$. Since $\mcal V_\eta \subset \{|\Im w| < \eta\}$, the set $E$ is contained in a rectangle with sides $2r$ and $2\eta$, so that its area satisfies $|E| \leq 4\eta r$. We bound $u_N$ separately on $E$ and on $D(y,r)\setminus E$. On $D(y,r)\setminus E$ we have $u_N \leq \psi + C_\eta/N$. On $E$ we have $u_N \leq C_1 + 1 \leq \psi + 2C_1 + 1$, because $-\psi \leq C_1$ on $\{|w| \leq R + 1\}$. By the sub-mean value property,
        \begin{align*}
            u_N(y) &\leq \frac{1}{\pi r^2}\int_{D(y,r)} u_N\,\dd A
            = \frac{1}{\pi r^2}\left(\int_{D(y,r)\setminus E} u_N\,\dd A + \int_E u_N\,\dd A\right) \\
            &\leq \frac{1}{\pi r^2}\left(\int_{D(y,r)\setminus E}\Big(\psi + \frac{C_\eta}{N}\Big)\,\dd A
            + \int_E \left(\psi + 2C_1 + 1\right)\dd A\right) \\
            &\leq \frac{1}{\pi r^2}\int_{D(y,r)} \psi\,\dd A + \frac{C_\eta}{N}
            + \frac{(2C_1 + 1)\,|E|}{\pi r^2} \\
            &\leq \frac{1}{\pi r^2}\int_{D(y,r)} \psi\,\dd A + \frac{C_\eta}{N}
            + \frac{4\eta}{\pi r}\left(2C_1 + 1\right).
        \end{align*}
        The first term on the right-hand side is at most $\psi(y) + \omega(r)$, where $\omega$ is the modulus of continuity of $\psi$ on $\{|y| \leq R + 1\}$.

        It remains to choose the parameters. First choose $r$ with $\omega(r) < \frac\epsilon3$, then $\eta$ with $4\eta(2C_1 + 1) < \frac{\pi r\epsilon}{3}$, and finally $N_0$ such that $C_\eta/N < \frac\epsilon3$ for $N \geq N_0$. With these choices, $u_N \leq \psi + \epsilon$ on $L$ for $N \geq N_0$.
    \end{proof}

    \begin{proof}[Proof of Theorem~\ref{thm:intro-asymptotics}]
        Let $K$ be a compact subset of $\C\setminus\left(\Gamma\cup\msf z([a_1,b_1]\cap\{\SW < 0\})\right)$, let $K_- \defeq x_-(K)$, let
        \[
            w_1 \defeq \min_{\zeta \in K}|\SW(x_-(\zeta))| > 0,
        \]
        and fix $c \in (0, 2w_1)$. Then $K_-$ is a compact subset of $\C\setminus\{0,\ii,-\ii\}$, and $\SW \leq -w_1 < 0$ on $K_-$. Moreover $K_-$ does not meet $[a_1,b_1]$. Indeed, $x_-(\zeta) \in [a_1,b_1]$ would imply $\zeta \in \msf z([a_1,b_1]\cap\{\SW < 0\})$, contrary to the choice of $K$. Hence $K_-$ is a compact subset of $\C\setminus\left([a_1,b_1]\cup\{0\}\right)$.

        Comparing \eqref{eq:PN-two-sheets} with Theorem~\ref{thm:intro-asymptotics}, it remains to show that $\tau_{ij}(x) = \Boh(\ee^{-cN})$ uniformly for $x \in K_-$.

        Since $A_f$ and $A_g$ are continuous and do not vanish on $K_-$, Theorem~\ref{thm:intro-asymptotics-fg} gives $c_0 > 0$ such that, for $N$ large,
        \begin{equation}\label{eq:fg-lower}
            |f(x)| \geq c_0\,\msf a^{-N}\ee^{-2N(\U^\mu(x) + \msf g(0))}\,|x|^{-N},
            \qquad
            |g(x)| \geq c_0\,\msf a^{-N}\ee^{-2N\U^\mu(x)}\,|x|^{-N+1},
            \qquad x \in K_-.
        \end{equation}
        In particular $f(x)\,g(x) \neq 0$ for $x \in K_-$, so that \eqref{eq:PN-two-sheets} holds at every $x \in K_-$.

        Fix $\epsilon = w_1 - \frac{c}{2} > 0$. The set $\iota(K_-)$ is a compact subset of $\C\setminus\{0\}$. By Lemma~\ref{lem:upper-bound} with $L = \iota(K_-)$, and since $|f(\iota x)| = |x|^{N}|p(\iota x)|$, we have
        \[
            |f(\iota x)| \leq \msf a^{-N}\ee^{2N(\epsilon - \U^\mu(\iota x) - \msf g(0))}|x|^N,
            \qquad x \in K_-,
        \]
        for $N$ large. Together with the first bound in \eqref{eq:fg-lower} and with $\SW(x) = \U^\mu(x) - \U^\mu(\iota x) + \log|x|$, this gives
        \[
            |\tau_{12}(x)| = \left|\frac{f(\iota x)}{f(x)}\right|
            \leq \frac{1}{c_0}\,\ee^{2N(\SW(x) + \epsilon)}
            \leq \frac{1}{c_0}\,\ee^{-Nc},
            \qquad x \in K_-.
        \]
        The same argument, applied to $g(\iota x) = (\iota x)^{-N+1}q(\iota x)$ in place of $f(\iota x)$ and to the second bound in \eqref{eq:fg-lower}, yields $|\tau_{22}(x)| \leq c_0^{-1}|x|^{-2}\ee^{-Nc}$ for $x \in K_-$.

        Finally, $\tau_{11}(x) = \frac{x}{\iota x}\,\tau_{12}(x)$ and $\tau_{21}(x) = \frac{x}{\iota x}\,\tau_{22}(x)$, with $|x/\iota x| = |x|^2$, and $|x|^{\pm2}$ is bounded on $K_-$. Therefore $|\tau_{ij}(x)| \leq C\,\ee^{-Nc}$ for $x \in K_-$ and $N$ large.
    \end{proof}

\subsection{Proof of Theorem~\ref{thm:intro-zeros-general}}\label{sec:proof-zeros-general} \hfill

    Before we start with the proof we motivate how the limit distribution \eqref{eq:intro-zeros-general} arises. Fix $i,j \in \{1,2\}$ and put
    \[
        \Phi_N(\zeta) \defeq \frac1N\log|\bm P_{ij}(\zeta)|.
    \]
    If $\bm P_{ij} \not\equiv 0$, then, by \eqref{eq:log-poly}, $\frac{1}{2\pi}\Delta\Phi_N = \frac1N\sum\delta_\zeta$, the sum running over the zeros of $\bm P_{ij}$ counted with multiplicity. This measure is $\frac{\deg\bm P_{ij}}{N}$ times the normalized zero counting measure of $\bm P_{ij}$, and we look for the limit of $\Phi_N$.

    Let $\zeta \in \C\setminus\left(\Gamma\cup\msf z([a_1,b_1])\right)$ and $x \defeq x_-(\zeta)$, so that $\SW(x) < 0$ and $x \notin [a_1,b_1]$. By Theorem~\ref{thm:intro-asymptotics}, $\bm P_{ij}(\zeta)$ is, up to a factor $1 + o(1)$, the $(i,j)$ entry of the rank-one matrix
    \[
        \frac{1}{x - \iota x}\left(f(x),\,g(x)\right)^T\left(-\iota x,\,1\right),
    \] which only involves the dominant point $x = x_-(\zeta)$. Hence Theorem~\ref{thm:intro-asymptotics-fg} at the point $x$ gives
    \[
        \Phi_N(\zeta) \longrightarrow -2\,\U^\mu(x) - \log|x| - \log\msf a + \beta_i, \qquad \beta_1 \defeq -2\msf g(0), \quad \beta_2 \defeq 0.
    \]
    We use Lemma~\ref{lem:iota-potentials}(iii) to write the right-hand side as a function of $\zeta$. Since $\U^\mu(0) = -\msf g(0)$, and $|\SW(x)| = -\SW(x) = \U^\mu(\iota x) - \U^\mu(x) - \log|x|$, we obtain
    \[
        -2\,\U^\mu(x) - \log|x| - \log\msf a = -\U^{\msf z_*\mu}(\zeta) + |\SW(x)| + \msf g(0),
    \]
    which leads us to consider the function
    \begin{equation}\label{eq:Phi-def}
        \Phi(\zeta) \defeq -\U^{\msf z_*\mu}(\zeta) + |\SW(x)|, \qquad \msf z(x) = \zeta.
    \end{equation}
    It does not depend on the choice of the solution $x$ of $\msf z(x) = \zeta$, since $\SW(\iota x) = -\SW(x)$. We have proven that, with $c_1 \defeq -\msf g(0)$ and $c_2 \defeq \msf g(0)$,
    \begin{equation}\label{eq:PhiN-limit}
        \lim_{N\to\infty}\Phi_N(\zeta) = \Phi(\zeta) + c_i, \qquad \zeta \in \C\setminus\left(\Gamma\cup\msf z([a_1,b_1])\right).
    \end{equation}
    Therefore it is natural to expect that the measures $\frac{1}{2\pi}\Delta\Phi_N$ converge to $\frac{1}{2\pi}\Delta\Phi$. We identify this measure with the limit measure in \eqref{eq:intro-zeros-general}, and then we prove the weak convergence.

    By Lemma~\ref{lem:iota-potentials},
    \begin{equation}\label{eq:Phi-z}
        \Phi\circ\msf z = |\SW| - \U^\mu - \U^{\wt\mu} - \log|x| - \log\msf a \quad \text{in } \C\setminus\{0\},
    \end{equation}
    where $\wt\mu = \iota_*\mu$ is the push-forward of $\mu$ under $\iota$ (see Section~\ref{sec:W}), so the main step in the computation of $\frac{1}{2\pi}\Delta\Phi$ is the Laplacian of $|\SW|$. Outside $\Gamma_{\SW}$ we have $|\SW| = \pm\SW$, whose Laplacian is given by Lemma~\ref{lem:W-properties}(i). Across $\Gamma_{\SW}$ the function $|\SW|$ has a corner: its normal derivative jumps by $2|\nabla\SW|$, and this jump produces the measure $\nu$ of Theorem~\ref{thm:intro-zeros-general} (see Appendix~\ref{app:prelim}). This leads to the following result.

    \begin{lemma}\label{lem:laplacian-absW}
        The measure $\nu$ is locally finite on $\C$, and, in the sense of distributions,
        \begin{equation}\label{eq:laplacian-absW}
            \frac{1}{2\pi}\Delta|\SW| = \sign(\SW)\,(\wt\mu - \mu) + 2\nu \quad \text{on } \C.
        \end{equation}
    \end{lemma}

    \begin{proof}
        Let $\varpi \defeq \frac{1}{2\pi}\Delta|\SW|$, and let $Z \defeq \Gamma_{\SW}\cap \mcal K$, which is a finite set by Lemma~\ref{lem:W-zero-set}(i). We proceed in three steps. First we show that $\varpi$ is a signed measure without atoms. Then, since the Laplacian is a local operator, we identify $\varpi$ in a neighborhood of each point of $\C\setminus Z$. Finally we deal with the finitely many points of $Z$.

        By Lemma~\ref{lem:W-properties}(i), $\SW = \SW_1 - \SW_2$ with $\SW_1 \defeq \msf g(0) - \U^{\wt\mu}$ and $\SW_2 \defeq -\U^\mu$. The functions $\SW_1$, $\SW_2$ and $\max\{\SW_1,\SW_2\}$ are continuous and subharmonic on $\C$, with $\frac{1}{2\pi}\Delta\SW_1 = \wt\mu$ and $\frac{1}{2\pi}\Delta\SW_2 = \mu$. Since $|\SW| = 2\max\{\SW_1,\SW_2\} - \SW_1 - \SW_2$, we have
        \[
            \varpi = 2\sigma - \wt\mu - \mu,
        \]
        where $\sigma$ is the Riesz measure of $\max\{\SW_1,\SW_2\}$. Hence $\varpi$ is a signed measure. It has no atoms, because $\mu$ and $\wt\mu$ have bounded densities and $\sigma$ has no atoms by Corollary~\ref{cor:riesz-no-atoms}.

        In the open sets $\{\SW > 0\}$ and $\{\SW < 0\}$ we have $|\SW| = \pm\SW$, and thus $\varpi = \pm(\wt\mu - \mu)$.

        Now let $x_0 \in \Gamma_{\SW}\setminus \mcal K$. By Lemma~\ref{lem:W-properties}(i), $\SW$ is harmonic near $x_0$, so $\SW = \Re\Psi$ near $x_0$ for an analytic function $\Psi$ with $\Psi(x_0) = 0$. By \eqref{eq:gradient-harmonic} and \eqref{eq:gradient-W},
        \[
            \Psi' = \partial_\xi\SW - \ii\,\partial_\upsilon\SW = \CT^\mu - \CT^{\wt\mu},
            \qquad \text{so that} \qquad |\Psi'| = |\nabla\SW|,
        \]
        and $\Psi'(x_0) \neq 0$ by Lemma~\ref{lem:W-critical}. Hence, for $\epsilon > 0$ small, $\Psi$ maps a neighborhood $U$ of $x_0$, with $U \cap \mcal K = \emptyset$, conformally onto the disk $\{|\omega| < \epsilon\}$. As $\SW = \Re\Psi$, it maps $\Gamma_{\SW}\cap U$ onto the segment $\{\ii t : |t| < \epsilon\}$, and $x(t) \defeq \Psi^{-1}(\ii t)$, $|t| < \epsilon$, parametrizes $\Gamma_{\SW}\cap U$. Since $x'(t) = \ii/\Psi'(x(t))$, we have
        \begin{equation}\label{eq:arc-length-nu}
            \dd t = |\nabla\SW(x)|\,|\dd x| = 2\pi\,\dd\nu(x)
            \quad \text{along } \Gamma_{\SW}\cap U.
        \end{equation}
        In the variable $\omega = s + \ii t = \Psi(x)$ the function $|\SW| = |\Re\Psi|$ becomes $|s|$. Let $\chi$ be a smooth function with compact support in $U$, and put $\hat\chi \defeq \chi\circ\Psi^{-1}$, a smooth function with compact support in $\{|\omega| < \epsilon\}$, which we extend by zero to $\C$. Since $\Psi$ is conformal, $\Delta\chi = |\Psi'|^2\,(\Delta\hat\chi)\circ\Psi$ and $\dd A(\omega) = |\Psi'(x)|^2\,\dd A(x)$, see Appendix~\ref{app:prelim}. Hence, by the change of variables $\omega = \Psi(x)$,
        \[
            \int_U |\SW|\,\Delta\chi\,\dd A
            = \int_U |\Re\Psi(x)|\,(\Delta\hat\chi)(\Psi(x))\,|\Psi'(x)|^2\,\dd A(x)
            = \int_{\R^2} |s|\,\big(\partial_s^2\hat\chi + \partial_t^2\hat\chi\big)(s + \ii t)\,\dd s\,\dd t.
        \]
        We compute the two terms separately, by Fubini's theorem. For fixed $s$, the function $t \mapsto \hat\chi(s + \ii t)$ has compact support, so that $\int_\R \partial_t^2\hat\chi(s + \ii t)\,\dd t = 0$, and the second term vanishes. For fixed $t$, the function $f(s) \defeq \hat\chi(s + \ii t)$ has compact support, and integration by parts on each half-line gives
        \[
            \int_0^\infty s\,f''(s)\,\dd s = -\int_0^\infty f'(s)\,\dd s = f(0),
            \qquad
            \int_{-\infty}^0 (-s)\,f''(s)\,\dd s = \int_{-\infty}^0 f'(s)\,\dd s = f(0),
        \]
        so that $\int_\R |s|\,f''(s)\,\dd s = 2f(0) = 2\hat\chi(\ii t)$. Therefore
        \[
            \int_U |\SW|\,\Delta\chi\,\dd A
            = 2\int_\R \hat\chi(\ii t)\,\dd t
            = 2\int_{-\epsilon}^{\epsilon}\chi(x(t))\,\dd t
            = 4\pi\int\chi\,\dd\nu.
        \]
        In the second equality we used that $\hat\chi(\ii t) = 0$ for $|t| \geq \epsilon$ and that $\hat\chi(\ii t) = \chi(x(t))$, by the definition of $x(t)$. The third one is \eqref{eq:arc-length-nu}, since $x(t)$ parametrizes $\Gamma_{\SW}\cap U$, which contains the part of the support of $\nu$ in $U$. Dividing by $2\pi$, we obtain $\varpi = 2\nu$ in $U$.

        Every point of $\C\setminus Z$ lies in $\{\SW \neq 0\}$ or in a neighborhood $U$ as above. Since $\mu$ and $\wt\mu$ vanish in these neighborhoods, and $\nu$ vanishes in $\{\SW \neq 0\}$, the identity \eqref{eq:laplacian-absW} holds in $\C\setminus Z$.

        Now let $x_1 \in Z$, and let $D$ be an open disk centered at $x_1$ with $D \cap Z = \{x_1\}$. For every continuous function $\chi$ with compact support in $D\setminus\{x_1\}$ and $0 \leq \chi \leq 1$, the identity in $D\setminus\{x_1\}$ gives $2\int\chi\,\dd\nu \leq |\varpi|(D) + \mu(D) + \wt\mu(D)$. Taking the supremum over $\chi$, and since $\nu(\{x_1\}) = 0$, we get $2\nu(D) \leq |\varpi|(D) + \mu(D) + \wt\mu(D) < \infty$. Hence $\nu$ is locally finite on $\C$. Both sides of \eqref{eq:laplacian-absW} are signed measures without atoms, and they agree on $\C\setminus Z$, so they agree on $\C$.
    \end{proof}

    From \eqref{eq:Phi-z} and Lemma~\ref{lem:laplacian-absW} we obtain the Laplacian of $\Phi$.

    \begin{lemma}\label{lem:Phi}
        In the sense of distributions,
        \begin{equation}\label{eq:Phi-laplacian}
            \frac{1}{2\pi}\Delta\Phi = \msf z_*\nu + 2\,\msf z_*\big(\restr{\mu}{\{\SW < 0\}}\big).
        \end{equation}
    \end{lemma}

    \begin{proof}
        We first verify \eqref{eq:Phi-z}. By Lemma~\ref{lem:iota-potentials}(iii) with $\U^\mu(0) = -\msf g(0)$, and by \eqref{eq:potential-iota} at the point $\iota x$, which gives $\U^\mu(\iota x) = \U^{\wt\mu}(x) + \log|x| - \msf g(0)$,
        \[
            \U^{\msf z_*\mu}(\msf z(x)) = \U^\mu(x) + \U^{\wt\mu}(x) + \log|x| + \log\msf a,
            \qquad x \in \C\setminus\{0\},
        \]
        and \eqref{eq:Phi-z} follows from \eqref{eq:Phi-def}. Applying Lemma~\ref{lem:laplacian-absW} in $\C\setminus\{0\}$ therefore gives
        \[
            M \defeq \frac{1}{2\pi}\Delta(\Phi\circ\msf z)
            = \sign(\SW)\,(\wt\mu - \mu) + \mu + \wt\mu + 2\nu
            = 2\,\restr{\mu}{\{\SW < 0\}} + 2\,\restr{\wt\mu}{\{\SW > 0\}} + 2\nu.
        \]

        We now pass from $\msf z$ back to $\C$. The map $\msf z$ is analytic in $\C\setminus\{0\}$, and it maps each of the sets $\{|x| > 1\}$ and $\{0 < |x| < 1\}$ conformally onto $\C\setminus[-\frac{2\ii}{\msf a},\frac{2\ii}{\msf a}]$, since $\msf z\circ\iota = \msf z$ and $\iota$ exchanges these two sets. The unit circle, which $\msf z$ maps onto the segment $[-\frac{2\ii}{\msf a},\frac{2\ii}{\msf a}]$, and this segment have zero area. Hence, applying the change of variables formula on each of the two sets, with Jacobian $|\msf z'|^2$, we get
        \begin{equation}\label{eq:z-change-variables}
            \int_{\C\setminus\{0\}} F(\msf z(x))\,|\msf z'(x)|^2\,\dd A(x) = 2\int_\C F\,\dd A
        \end{equation}
        for every integrable function $F$ on $\C$.

        Let $\chi$ be a smooth function with compact support in $\C$. Then $\chi\circ\msf z$ is smooth with compact support in $\C\setminus\{0\}$, because $\msf z(x) \to \infty$ as $x \to 0$ and as $x \to \infty$. Moreover $\Delta(\chi\circ\msf z) = |\msf z'|^2\,(\Delta\chi)\circ\msf z$, because $\msf z$ is analytic. We apply \eqref{eq:z-change-variables} to $F = \Phi\,\Delta\chi$, which is integrable because $\Phi$ is locally integrable. Then we use the definition of $M$ with the test function $\chi\circ\msf z$. We obtain
        \[
            \int\Phi\,\Delta\chi\,\dd A
            = \frac12\int_{\C\setminus\{0\}}(\Phi\circ\msf z)\,\big((\Delta\chi)\circ\msf z\big)\,|\msf z'|^2\,\dd A
            = \frac12\int_{\C\setminus\{0\}}(\Phi\circ\msf z)\,\Delta(\chi\circ\msf z)\,\dd A
            = \pi\int\chi\circ\msf z\,\dd M.
        \]
        Since $\int\chi\circ\msf z\,\dd M = \int\chi\,\dd(\msf z_*M)$, this says that $\frac{1}{2\pi}\Delta\Phi = \frac12\msf z_*M$.

        It remains to identify $\frac12\msf z_*M$. Since $\wt\mu = \iota_*\mu$, $\msf z\circ\iota = \msf z$ and $\iota(\{\SW < 0\}) = \{\SW > 0\}$, we have
        \[
            \msf z_*\big(\restr{\wt\mu}{\{\SW > 0\}}\big)
            = \msf z_*\big(\restr{\mu}{\{\SW < 0\}}\big),
        \]
        and \eqref{eq:Phi-laplacian} follows.
    \end{proof}

    In particular, $\frac{1}{2\pi}\Delta\Phi$ is a positive measure. Its total mass is determined by the behavior of $\Phi$ at infinity, which we describe now.

    \begin{lemma}\label{lem:Phi-continuity}
        The function $\Phi$ is continuous on $\C$, and
        \begin{equation}\label{eq:Phi-infinity}
            \Phi(\zeta) = \log|\zeta| + |\SW(0)| + \boh(1), \qquad \zeta \to \infty.
        \end{equation}
    \end{lemma}

    \begin{proof}
        We prove the two assertions separately, first the continuity of $\Phi$ and then its behavior at infinity.

        For the continuity, note that $\msf z_*\mu$ has a bounded density on the segment $\msf z([a_1,b_1])$, because $\mu$ has a bounded density and $\msf z' > 0$ on $[a_1,b_1]$. Hence $\U^{\msf z_*\mu}$ is continuous on $\C$. Since $\msf z_*\mu$ is a probability measure with compact support, we also have $\U^{\msf z_*\mu}(\zeta) = -\log|\zeta| + \boh(1)$ as $\zeta \to \infty$.

        It remains to handle the second term in \eqref{eq:Phi-def}. If $\zeta_n \to \zeta$, the solutions of $\msf z(x) = \zeta_n$ stay in a compact subset of $\C\setminus\{0\}$, and all their limit points are solutions of $\msf z(x) = \zeta$. Since $\SW$ is continuous, this makes the second term in \eqref{eq:Phi-def} a continuous function of $\zeta$.

        Finally, as $\zeta \to \infty$, the solutions of $\msf z(x) = \zeta$ tend to $0$ and to $\infty$, and $|\SW(x)| \to |\SW(0)|$ by Lemma~\ref{lem:W-properties}(i) and~(ii). This proves \eqref{eq:Phi-infinity}.
    \end{proof}

    \begin{lemma}\label{lem:Phi-mass}
        The measure $\frac{1}{2\pi}\Delta\Phi$ is a probability measure. Consequently, $\nu$ is a finite measure, with mass $\mu(\{\SW > 0\}) - \mu(\{\SW < 0\}) > 0$.
    \end{lemma}

    \begin{proof}
        Let $\rho \defeq \frac{1}{2\pi}\Delta\Phi$, which is a positive measure by Lemma~\ref{lem:Phi}. We first show that $\rho(\C) = 1$, and then we compute the total mass of $\nu$.

        Let $\chi$ be a smooth radial function on $\C$ with $0 \leq \chi \leq 1$, non-increasing in $|\zeta|$, equal to $1$ on $\{|\zeta| \leq 1\}$ and to $0$ on $\{|\zeta| \geq 2\}$, and for $R > 0$ put $\chi_R(\zeta) \defeq \chi(\zeta/R)$. These cut-off functions have the following properties.
        \begin{enumerate}[\rm (a)]
            \item For each $\zeta$, $\chi_R(\zeta)$ is non-decreasing in $R$ and equal to $1$ for $R \geq |\zeta|$. Since $\rho$ is positive, the monotone convergence theorem gives $\int\chi_R\,\dd\rho \to \rho(\C)$ as $R \to \infty$.
            \item $\Delta\chi_R(\zeta) = R^{-2}(\Delta\chi)(\zeta/R)$ vanishes outside the annulus $A_R \defeq \{R \leq |\zeta| \leq 2R\}$, and the change of variables $\zeta = R\omega$ gives $\int|\Delta\chi_R|\,\dd A = \int|\Delta\chi|\,\dd A$, which does not depend on $R$.
            \item $\int\Delta\chi_R\,\dd A = 0$, because $\chi_R$ has compact support, and $\frac{1}{2\pi}\int\log|\zeta|\,\Delta\chi_R(\zeta)\,\dd A(\zeta) = \chi_R(0) = 1$, because $\frac{1}{2\pi}\Delta\log|\zeta| = \delta_0$.
        \end{enumerate}
        By Lemma~\ref{lem:Phi-continuity} we can write $\Phi(\zeta) = \log|\zeta| + |\SW(0)| + E(\zeta)$, where $E(\zeta) \to 0$ as $\zeta \to \infty$. By the definition of the distributional Laplacian and (c),
        \[
            \int\chi_R\,\dd\rho
            = \frac{1}{2\pi}\int\Phi\,\Delta\chi_R\,\dd A
            = 1 + \frac{1}{2\pi}\int E\,\Delta\chi_R\,\dd A.
        \]
        By (b), the last term is bounded by $\frac{1}{2\pi}\max_{A_R}|E|\int|\Delta\chi|\,\dd A$, which tends to $0$ as $R \to \infty$. Hence $\int\chi_R\,\dd\rho \to 1$, and (a) gives $\rho(\C) = 1$.

        By Lemma~\ref{lem:Phi}, $\rho = \msf z_*\nu + 2\,\msf z_*\big(\restr{\mu}{\{\SW < 0\}}\big)$. The measure $\msf z_*\nu$ has the same total mass as $\nu$, because $\nu(\{0\}) = 0$, and $2\,\msf z_*\big(\restr{\mu}{\{\SW < 0\}}\big)$ has total mass $2\mu(\{\SW < 0\})$. Hence $\nu(\C) = 1 - 2\mu(\{\SW < 0\})$, and in particular $\nu$ is finite. Since $\mu(\{\SW = 0\}) = 0$ by Lemma~\ref{lem:W-zero-set}(i), we have $\mu(\{\SW > 0\}) + \mu(\{\SW < 0\}) = 1$, and therefore
        \[
            \nu(\C) = \mu(\{\SW > 0\}) - \mu(\{\SW < 0\}).
        \]

        Regarding the mass of $\nu$, it suffices to show that $\nu$ does not vanish since $\nu$ is a positive measure. By Lemma~\ref{lem:W-properties}(ii), $\ii \in \Gamma_{\SW}\setminus \mcal K$. Near $\ii$ the function $\SW$ is harmonic with non-vanishing gradient (Lemmas~\ref{lem:W-properties}(i) and \ref{lem:W-critical}). Hence, by the implicit function theorem, $\Gamma_{\SW}$ contains an analytic arc $\gamma$ through $\ii$ on which the density $\frac{1}{2\pi}|\nabla\SW|$ of $\nu$ is positive. Therefore $\nu(\C) \geq \nu(\gamma) > 0$.
    \end{proof}

    We can now prove the weak convergence.

    \begin{proof}[Proof of Theorem~\ref{thm:intro-zeros-general}]
        Fix $i,j \in \{1,2\}$, and let $\Phi_N$ and $c_i$ be as in \eqref{eq:PhiN-limit}. Since $\bm P_N$ is monic of degree $N$, its diagonal entries are monic of degree $N$, and its off-diagonal entries have degree at most $N - 1$. Let $\sigma_N \defeq \frac1N\sum\delta_\zeta = \frac{1}{2\pi}\Delta\Phi_N$, the sum running over the zeros of $\bm P_{ij}$ counted with multiplicity. Since $\deg\bm P_{ij} \leq N$, we have $\sigma_N(\C) = \deg\bm P_{ij}/N \leq 1$. Let $\rho$ be the measure \eqref{eq:intro-zeros-general}. By Lemmas~\ref{lem:Phi} and~\ref{lem:Phi-mass}, $\rho = \frac{1}{2\pi}\Delta\Phi$ is a probability measure.

        As previously mentioned, the main idea is to pass to the limit in $\sigma_N = \frac{1}{2\pi}\Delta\Phi_N$. By \eqref{eq:PhiN-limit} and Lemma~\ref{lem:W-zero-set}(ii), $\Phi_N \to \Phi + c_i$ almost everywhere, which does not allow us to take Laplacians directly. We therefore upgrade this to convergence in $L^1_{\mathrm{loc}}(\C)$ by means of a compactness theorem for subharmonic functions, which requires a uniform upper bound for $\Phi_N$. Once the Laplacians converge, the fact that $\rho$ is a probability measure will show that no mass of $\sigma_N$ escapes to infinity.

        To prove the upper bound, let $R > \frac{2}{\msf a}$. The set $\Sigma_R \defeq \{x \in \C\setminus\{0\} : |\msf z(x)| = R\}$ is a compact subset of $\C\setminus\{0,\ii,-\ii\}$, invariant under $\iota$, so that $|x - \iota x| = |1 + x^2|/|x|$ is bounded below on $\Sigma_R$. By Lemma~\ref{lem:upper-bound}, $|f|$ and $|g|$ are bounded by $\ee^{CN}$ on $\Sigma_R$ for $N$ large. By Lemma~\ref{lem:divided-differences}, $|\bm P_{ij}| \leq \ee^{C'N}$ on $\{|\zeta| = R\}$, where $C'$ depends on $R$ but not on $N$. Hence, by the maximum principle, $\Phi_N \leq C'$ on $\{|\zeta| \leq R\}$. Since $R$ is arbitrary, the functions $\Phi_N$ are bounded above on every compact set, uniformly in $N$.

        Next we apply compactness. The functions $\Phi_N$ are subharmonic on $\C$ and, since $\bm P_{ij} \not\equiv 0$, not identically $-\infty$. By \cite[Theorem~3.2.12]{Hormander1994}, every subsequence of $(\Phi_N)$ has a further subsequence that converges in $L^1_{\mathrm{loc}}(\C)$. The other alternative in that theorem, namely $\Phi_N \to -\infty$ uniformly on compact sets, is excluded by the almost everywhere convergence. Passing to a subsequence that converges almost everywhere, we see that the limit is $\Phi + c_i$. Therefore $\Phi_N \to \Phi + c_i$ in $L^1_{\mathrm{loc}}(\C)$.

        We can finally pass to the measures. In the sense of distributions, $\sigma_N = \frac{1}{2\pi}\Delta\Phi_N \to \frac{1}{2\pi}\Delta\Phi = \rho$, and since the $\sigma_N$ are positive measures, $\int\chi\,\dd\sigma_N \to \int\chi\,\dd\rho$ for every continuous function $\chi$ with compact support, that is, $\sigma_N \to \rho$ vaguely. Since a vague limit cannot gain mass (see Appendix~\ref{app:prelim}),
        \[
            1 = \rho(\C) \leq \liminf_{N\to\infty}\sigma_N(\C)
            \leq \limsup_{N\to\infty}\sigma_N(\C) \leq 1.
        \]
        Hence $\deg\bm P_{ij}/N = \sigma_N(\C) \to \rho(\C)$, and therefore $\sigma_N \to \rho$ weakly. In particular $\deg\bm P_{ij} \geq 1$ for $N$ large, and the normalized zero counting measure $\frac{N}{\deg\bm P_{ij}}\,\sigma_N$ of $\bm P_{ij}$ converges weakly to $\rho$.
    \end{proof}

    \begin{remark}[Balayage]\label{rem:balayage}
        Assume that $\SW(0) \neq 0$, and let $\mu_\pm$ be the restriction of $\mu$ to $\{\pm\SW > 0\}$. Then
        \begin{equation}\label{eq:balayage}
            {\U^{\msf z_*\nu}} = {\U^{\msf z_*\mu_+}} - {\U^{\msf z_*\mu_-}} + |\SW(0)|
            \quad \text{on } \Gamma.
        \end{equation}
        We prove \eqref{eq:balayage} after this remark.

        The meaning of \eqref{eq:balayage} is clearest when $\SW > 0$ on $[a_1,b_1]$. Then $\mu_+ = \mu$ and $\mu_- = 0$, and $\Gamma$ does not meet the segment $\msf z([a_1,b_1])$, because $\SW$ does not vanish on $[a_1,b_1]\cup[\wt a_1,\wt b_1]$. In this case $\msf z_*\nu$ is a probability measure on $\Gamma$, and \eqref{eq:balayage} says that its potential coincides on $\Gamma$ with the potential of $\msf z_*\mu$, up to an additive constant. By \cite[Theorems~II.4.6 and~II.4.7]{SaffTotik2024}, this property characterizes the balayage of $\msf z_*\mu$ onto $\Gamma$. In other words, in this case the limit of the zeros of the entries of $\bm P_N$ is the balayage onto $\Gamma$ of the limit of the zeros of $\det\bm P_N$.

        In general, \eqref{eq:balayage} says that the potential of $\msf z_*\nu$ coincides on $\Gamma$, up to an additive constant, with the potential of the signed measure $\msf z_*\mu_+ - \msf z_*\mu_-$, which does not charge $\Gamma$ and has the same total mass, by Lemma~\ref{lem:Phi-mass}. In this sense, $\msf z_*\nu$ is the difference of the balayages onto $\Gamma$ of $\msf z_*\mu_+$ and $\msf z_*\mu_-$. This is the form of the limit in \cite[Theorem~2.5]{Groot2021a}, with the role of $\Gamma$ played by the imaginary axis.
    \end{remark}

    \begin{proof}[Proof of \eqref{eq:balayage}]
        Recall that $\SW(0) \neq 0$. Since $\SW\circ\iota = -\SW$, the function $\SW$ tends to $-\SW(0) \neq 0$ at infinity. Hence $\Gamma_{\SW}$ is a compact subset of $\C\setminus\{0\}$, and $\Gamma = \msf z(\Gamma_{\SW}\setminus\{0\})$ is compact. Let $\rho$ be the measure \eqref{eq:intro-zeros-general}, which has compact support. By Lemmas~\ref{lem:Phi} and~\ref{lem:Phi-mass}, $\rho = \frac{1}{2\pi}\Delta\Phi$ is a probability measure.

        It follows that $\Phi + \U^\rho$ has zero Laplacian in the sense of distributions. By Weyl's lemma, $\Phi + \U^\rho$ coincides almost everywhere with a harmonic function $h$ on $\C$. In fact the equality holds everywhere, since $\Phi$, $\U^\rho$ and $h$ are, at every point, the limits of their means over small disks centered at that point (for $\U^\rho$ this follows from its superharmonicity).

        We now identify $h$. Since $\U^\rho(\zeta) = -\log|\zeta| + \boh(1)$ as $\zeta \to \infty$, Lemma~\ref{lem:Phi-continuity} gives $h(\zeta) \to |\SW(0)|$. Hence $h \equiv |\SW(0)|$ by Liouville's theorem, and therefore
        \[
            \U^\rho = |\SW(0)| - \Phi \quad \text{on } \C.
        \]

        Finally, we restrict this identity to $\Gamma$. By \eqref{eq:Phi-def}, $\Phi = -\U^{\msf z_*\mu}$ on $\Gamma$. Since $\rho = \msf z_*\nu + 2\,\msf z_*\mu_-$ and $\mu = \mu_+ + \mu_-$ by Lemma~\ref{lem:W-zero-set}(i), substituting into the identity above gives \eqref{eq:balayage}.
    \end{proof}

\appendix

\section{The equilibrium measure}\label{app:eq-measure} \hfill

    In this appendix we prove Proposition~\ref{prop:onecut} (Appendix~\ref{app:support}) and Theorem~\ref{thm:intro-eq} (Appendix~\ref{app:onecut}).

    In Appendix~\ref{app:support}, $\bm Q$ is any matrix polynomial of the form \eqref{eq:intro-Q} for which the external field $V$ of \eqref{eq:intro-V} is admissible, and $\mu$ is the equilibrium measure of $V$. We keep the notation of Section~\ref{sec:eq-measure}, and we write $\I(\sigma) \defeq \iint\log|x - y|^{-1}\,\dd\sigma(x)\,\dd\sigma(y)$ for the unweighted logarithmic energy, so that $\I^V(\sigma) = \I(\sigma) + \int V\,\dd\sigma$.

\needspace{5\baselineskip}
\subsection{Proof of Proposition~\ref{prop:onecut}}\label{app:support} \hfill

    The external field $V$ has the rational derivative
    \begin{equation}\label{eq:Vprime-general}
        \frac{V'(z)}{2} = \frac{\msf Q'(z)}{4} + \frac{1}{2z}, \qquad z \in \C\setminus\{0\},
    \end{equation}
    which has a pole of order $d_- + 1$ at $z = 0$ and satisfies $\frac{V'(z)}{2} = \frac{z^{d_+-1}}{\msf k}\left(1 + \Boh(z^{-1})\right)$ as $z \to \infty$, by \eqref{eq:k-def}. This allows us to describe the equilibrium measure in terms of a few real parameters. The following lemma is standard for polynomial external fields, see for instance \cite{DKM98}. We present its proof for the convenience of the reader.

    \begin{lemma}\label{lem:quad}
        The measure $\mu$ is absolutely continuous, and its support is the union of $M \leq d_+ + d_-$ disjoint compact intervals $[a_j,b_j] \subset \R\setminus\{0\}$, $a_1 < b_1 < a_2 < \dots < a_M < b_M$. There is a monic polynomial $\hat q$ of degree $d_+ + d_- - M$ with real coefficients such that
        \begin{equation}\label{eq:quad-eq}
            \left(\CT^\mu(z) + \frac{V'(z)}{2}\right)^2 = \frac{\prod_{j=1}^M (z - a_j)(z - b_j)\;\hat q(z)^2}{\msf k^2\,z^{2(d_- + 1)}}, \qquad z \in \C\setminus\left(\supp\mu\cup\{0\}\right),
        \end{equation}
        and
        \[
            \dd\mu(x) = \frac{\sqrt{-\prod_{j=1}^M (x - a_j)(x - b_j)}\;|\hat q(x)|}{\pi\,\msf k\,|x|^{d_- + 1}}\,\dd x, \qquad x \in \supp\mu.
        \]
        In particular, the density of $\mu$ is bounded, $\U^\mu$ is continuous on $\C$, $F \geq 0$ on $\R\setminus\{0\}$, and $F = 0$ on $\supp\mu$.
    \end{lemma}

    \begin{proof}
        Since $V$ is admissible, $\supp\mu$ is compact and $2\U^\mu + V \leq \ell$ on $\supp\mu$ \cite[Theorem~I.1.5]{SaffTotik2024}. As $\U^\mu$ is bounded below on compact sets and $V(x) \to +\infty$ as $x \to 0$, a neighborhood of $0$ does not meet $\supp\mu$.

        We follow the approach of Kamvissis and Rakhmanov \cite{KamvissisRakhmanov05}, based on Schiffer variations. The first variation of $\I^V$ at $\mu$ along the push-forwards of $\mu$ under $x \mapsto x + \epsilon h(x)$, with $h$ real, $C^1$ and with $h$, $h'$ bounded, vanishes. By linearity, this extends to complex $h$. Choosing $h(x) = \frac{1}{x - z}$, $z \in \C\setminus\R$, and completing the square gives
        \[
            \left(\CT^\mu(z) + \frac{V'(z)}{2}\right)^2
            = \frac{V'(z)^2}{4} - \int\frac{V'(x) - V'(z)}{x - z}\,\dd\mu(x),
            \qquad z \in \C\setminus\R.
        \]

        We now identify the right-hand side. By \eqref{eq:Vprime-general}, $V'$ is a Laurent polynomial with a pole of order $d_- + 1$ at $z = 0$, and
        \[
            \frac{V'(z)^2}{4} = \frac{z^{2(d_+-1)}}{\msf k^2}\left(1 + \Boh\left(z^{-1}\right)\right),
            \qquad z \to \infty.
        \]
        The integral is a Laurent polynomial in $z$, with a pole of order at most $d_- + 1$ at $z = 0$, and it is $\Boh(z^{d_+ - 2})$ as $z \to \infty$. Hence the right-hand side of the identity above equals $\frac{P(z)}{\msf k^2z^{2(d_- + 1)}}$, where $P$ is a monic polynomial of degree $2(d_+ + d_-)$ with real coefficients. Since $\C\setminus(\supp\mu\cup\{0\})$ is connected, the identity extends to it by analytic continuation.

        By this identity, $|\CT^\mu(z)| \leq \frac{|V'(z)|}{2} + \frac{|P(z)|^{1/2}}{\msf k\,|z|^{d_-+1}}$ for $z \in \C\setminus\R$, so that $\Im\CT^\mu(x + \ii\varepsilon) = \int\frac{\varepsilon\,\dd\mu(s)}{(s - x)^2 + \varepsilon^2}$ is bounded for $x$ in compact subsets of $\R\setminus\{0\}$ and $0 < \varepsilon \leq 1$. Since $0 \notin \supp\mu$, it follows that $\mu$ is absolutely continuous, with a bounded density.

        Since $V'$ is real on $\R\setminus\{0\}$ and $\Im\CT^\mu_+ \geq 0$, the Stieltjes inversion formula gives
        \[
            \mu'(x) = \frac{\sqrt{\max\{-P(x),0\}}}{\pi\,\msf k\,|x|^{d_- + 1}}
            \quad \text{for almost every } x \in \R.
        \]
        The right-hand side is continuous on $\R\setminus\{0\}$, so $\supp\mu$ is the closure of $\{x \in \R : P(x) < 0\}$, a finite union of disjoint compact intervals whose endpoints are the real zeros of $P$ of odd order. All the other zeros of $P$ have even order: for the real ones this follows from the sign of $P$, and for the non-real ones because $P(z)/z^{2(d_-+1)}$ is the square of the function $\msf k\left(\CT^\mu + \frac{V'}{2}\right)$, which is analytic in $\C\setminus\R$. Since $P$ is real and monic, $P = \prod_{j=1}^M(z - a_j)(z - b_j)\,\hat q^2$ with $\hat q$ real and monic of degree $d_+ + d_- - M$. This gives \eqref{eq:quad-eq}, the formula for the density and $M \leq d_+ + d_-$.

        Finally, $\mu$ has a bounded density and compact support, so $\U^\mu$ is continuous on $\C$. The inequality $2\U^\mu + V \geq \ell$ holds quasi-everywhere on $\R$, hence on a dense subset of $\R\setminus\{0\}$, and, by continuity, at every point of $\R\setminus\{0\}$. Together with $2\U^\mu + V \leq \ell$ on $\supp\mu$ \cite[Theorem~I.1.5]{SaffTotik2024}, this gives $F \geq 0$ on $\R\setminus\{0\}$ and $F = 0$ on $\supp\mu$.
    \end{proof}

    Proposition~\ref{prop:onecut} is simply a specialization of Lemma~\ref{lem:quad}. We now present its proof.

    \begin{proof}[Proof of Proposition~\ref{prop:onecut}]
        Let $\supp\mu = [a_1,b_1] \subset (-\infty,0)$. The function $\msf h$ of \eqref{eq:m-r-u} is analytic in $\C\setminus[a_1,b_1]$, with $\msf h(z) \sim z$ as $z \to \infty$, and it is positive on $(b_1,+\infty) \ni 0$. In particular $0 < \msf h(0) = u = \sqrt{|a_1||b_1|} < |m| = \frac{|a_1| + |b_1|}{2}$.

        By Lemma~\ref{lem:quad} with $M = 1$, there is a real monic polynomial $\msf q$ of degree $d_+ + d_- - 1$ such that
        \[
            \left(\CT^\mu(z) + \frac{V'(z)}{2}\right)^2
            = \left(\frac{\msf h(z)\,\msf q(z)}{\msf k\,z^{d_- + 1}}\right)^2,
            \qquad z \in \C\setminus\left([a_1,b_1]\cup\{0\}\right).
        \]
        The two functions appearing in this identity are analytic in the connected set $\C\setminus([a_1,b_1]\cup\{0\})$, so they agree up to a constant sign. As $z \to \infty$ both are $\frac{z^{d_+ - 1}}{\msf k}\left(1 + \Boh(z^{-1})\right)$, by \eqref{eq:k-def} and $d_+ \geq 2$. Hence the sign is $+$, which proves \eqref{eq:Cauchy-h-q}.

        The boundary values of $\msf h$ on $(a_1,b_1)$ are $\msf h_\pm(x) = \pm\ii\sqrt{(x - a_1)(b_1 - x)}$, so \eqref{eq:density-onecut} follows from \eqref{eq:Cauchy-h-q} and the Sokhotski--Plemelj formula $2\pi\ii\,\mu' = \CT^\mu_+ - \CT^\mu_-$. Since $\mu' \geq 0$, and $x^{d_- + 1} < 0$ on $[a_1,b_1]$ because $d_-$ is even, we get $\msf q \leq 0$ there. Consequently, (R2) holds if, and only if, $\msf q < 0$ on $(a_1,b_1)$, and (R3) holds if, and only if, $\msf q(a_1) \neq 0$ and $\msf q(b_1) \neq 0$.
    \end{proof}

\needspace{5\baselineskip}
\subsection{The Toda deformed Hermite-type weight: proof of Theorem~\ref{thm:intro-eq}}\label{app:onecut} \hfill

    In this subsection we write $\msf a$ and $\msf t$ for $\msf a_0$ and $\msf t_0$, to lighten the notation. Then $\bm Q(\zeta) = \zeta^2\bm I + \msf t\bm\Lambda(\zeta)$, as in \eqref{eq:intro-Q-toda}, so that $V = V_{\msf a,\msf t}$ is the external field \eqref{eq:intro-V-toda}, $d_\pm = 2$, $\msf k = 2\msf a^2$, and
    \begin{equation}\label{eq:Vprime-exact}
        \frac{V'(z)}{2} = \frac{z}{2\msf a^2} + \frac{1}{2z} - \frac{\msf t}{4z^2} - \frac{1}{2\msf a^2z^3}, \qquad z \in \C\setminus\{0\}.
    \end{equation}
    We write $\mu = \mu_{\msf a,\msf t}$, and we assume that $\msf a$ is fixed and that $\msf t > 0$.

    \begin{lemma}\label{lem:toda-onecut}
        If $\supp\mu = [a_1,b_1] \subset (-\infty,0)$, then the polynomial $\msf q$ of Proposition~\ref{prop:onecut} is
        \begin{equation}\label{eq:q-def}
            \msf q(z) = z^3 + m z^2 + \left( m^2 + \frac{r^2}{2} - \msf a^2 \right) z - \frac{1}{u}.
        \end{equation}
    \end{lemma}

    \begin{proof}
        Write $\msf q(z) = z^3 + \msf q_2z^2 + \msf q_1z + \msf q_0$. From $\msf h^2 = (z - m)^2 - r^2$ we get $\msf h(z) = z - m - \frac{r^2}{2z} + \Boh(z^{-2})$ as $z \to \infty$, and therefore
        \[
            \begin{split}
                \frac{\msf h(z)\,\msf q(z)}{2\msf a^2z^3} &= \frac{1}{2\msf a^2}\left[z + \left(\msf q_2 - m\right) + \frac{\msf q_1 - m\msf q_2 - \frac{r^2}{2}}{z} + \Boh\left(z^{-2}\right)\right], \\
                \CT^\mu(z) + \frac{V'(z)}{2} &= \frac{z}{2\msf a^2} - \frac{1}{2z} + \Boh(z^{-2}),
            \end{split}
        \]
        because $\CT^\mu(z) = -\frac1z + \Boh(z^{-2})$ and by \eqref{eq:Vprime-exact}. By \eqref{eq:Cauchy-h-q}, $\msf q_2 = m$ and $\msf q_1 = m^2 + \frac{r^2}{2} - \msf a^2$. Moreover, $\CT^\mu$ is analytic at $0 \notin \supp\mu$, while $\frac{V'(z)}{2} = -\frac{1}{2\msf a^2z^3} + \Boh(z^{-2})$ as $z \to 0$, by \eqref{eq:Vprime-exact}. Comparing the coefficients of $z^{-3}$ in \eqref{eq:Cauchy-h-q} gives $\msf h(0)\,\msf q(0) = u\,\msf q_0 = -1$.
    \end{proof}

    The explicit form of $\msf q$ gives regularity at once when the support is close enough to the origin.

    \begin{lemma}\label{lem:q-negative}
        Assume that $\supp\mu = [a_1,b_1] \subset (-\infty,0)$, with $a_1 \geq -\frac{1}{\msf a}$. Then $\msf q < 0$ on $\left[-\frac{1}{\msf a},0\right)$, and in particular on $[a_1,b_1]$. Hence (R2) and (R3) hold.
    \end{lemma}

    \begin{proof}
        By \eqref{eq:q-def}, for $z < 0$,
        \[
            \msf q(z) = \left[z^3 + m z^2 + \left(m^2 + \frac{r^2}{2}\right)z\right] + \left[\msf a^2|z| - \frac1u\right].
        \]
        The first bracket is negative, because $z < 0$ and $m < 0$. Since $|b_1| < |a_1|$, we have $u = \sqrt{|a_1||b_1|} < |a_1| \leq \frac{1}{\msf a}$. Hence, for $|z| \leq \frac{1}{\msf a}$, the second bracket is at most $\msf a - \frac1u < 0$. The last assertion follows from Proposition~\ref{prop:onecut}.
    \end{proof}

    It remains to localize the support for large $\msf t$, and to prove (R1). Both are more transparent in the variable $y = \iota x$.

    \begin{remark}[The transformed field]\label{rem:Vtilde-convex}
        With $\iota$ as in \eqref{eq:intro-z-iota}, let
        \begin{equation}\label{eq:Vtilde}
            \wt V(y) \defeq V(\iota y) + 2\log|y| = \frac{1}{2\msf a^2}\left(y - \frac1y\right)^2 - \frac{\msf t\,y}{2} + \log|y|, \quad y \in \R\setminus\{0\}, \qquad \wt V(0) \defeq +\infty.
        \end{equation}
        In $\wt V$ the parameter $\msf t$ multiplies a term that is linear in $y$. Completing the square, with $y_{\msf t} \defeq \frac{\msf a^2\msf t}{2}$,
        \begin{equation}\label{eq:Vtilde-square}
            \wt V(y) = \frac{1}{2\msf a^2}\left(y - y_{\msf t}\right)^2 - \frac{\msf a^2\msf t^2}{8} + R(y), \qquad R(y) \defeq -\frac{1}{\msf a^2} + \frac{1}{2\msf a^2y^2} + \log|y|,
        \end{equation}
        where $R$ does not depend on $\msf t$, and $R \geq R_{\min} \defeq \frac12 - \frac{1}{\msf a^2} - \log\msf a$, with equality at $|y| = \frac{1}{\msf a}$. Moreover,
        \[
            \wt V''(y) = \frac{1}{\msf a^2}\left(1 + \frac{3}{y^4}\right) - \frac{1}{y^2} = \frac{y^4 - \msf a^2y^2 + 3}{\msf a^2y^4}, \qquad y \neq 0,
        \]
        does not depend on $\msf t$, and $\wt V'' > 0$ on $\{|y| \geq \msf a\}$.
    \end{remark}

    Let $\wt\mu = \iota_*\mu$ be the push-forward of $\mu$ by $\iota$ (see Section~\ref{sec:W}). By \eqref{eq:iota-difference}, $\I(\iota_*\sigma) = \I(\sigma) + 2\int\log|s|\,\dd\sigma(s)$, and hence, by \eqref{eq:Vtilde},
    \begin{equation}\label{eq:energy-iota}
        \I^{\wt V}(\iota_*\sigma) = \I^V(\sigma)
    \end{equation}
    for every probability measure $\sigma$ with compact support in $\R\setminus\{0\}$. Moreover, by \eqref{eq:potential-iota} at the point $\iota y$, and since $V(\iota y) = \wt V(y) - 2\log|y|$,
    \begin{equation}\label{eq:F-iota}
        F(\iota y) = 2\,\U^{\wt\mu}(y) + \wt V(y) - \wt\ell, \qquad \wt\ell \defeq \ell + 2\int\log|s|\,\dd\mu(s), \qquad y \in \R\setminus\{0\}.
    \end{equation}
    Thus $y \mapsto F(\iota y)$ is the effective potential of $\wt\mu$ in the field $\wt V$: it is continuous on $\R\setminus\{0\}$, non-negative, and zero on $\supp\wt\mu$ (Lemma~\ref{lem:quad}). Integrating \eqref{eq:F-iota} against $\wt\mu$ gives $\wt\ell = 2\,\I(\wt\mu) + \int\wt V\,\dd\wt\mu = \I(\wt\mu) + \I^{\wt V}(\wt\mu)$.

    \begin{lemma}[Localization]\label{lem:localization}
        For every fixed $K > 0$,
        \[
            F(x) > 0, \quad x \in \left( -\infty, -\frac{1}{K} \right] \cup \left( 0, +\infty \right)
        \]
        for sufficiently large $\msf t$. In particular, $\supp\mu \subset \left(-\frac1K,0\right)$ for sufficiently large $\msf t$.
    \end{lemma}

    \begin{proof}
        In this proof $C$ denotes a positive constant that depends only on $\msf a$ and may change from one occurrence to the next. We assume that $\msf t \geq 1$ and $y_{\msf t} = \frac{\msf a^2 \msf t}{2} \geq 2$.

        We start with an upper bound for $\wt\ell$. Let $\sigma$ be the normalized Lebesgue measure on $[y_{\msf t} - 1, y_{\msf t} + 1] \subset [1,+\infty)$, so that $\I(\sigma) = \frac32 - \log 2$. On this interval, \eqref{eq:Vtilde-square} gives
        \[
            \wt V(y) \leq -\frac{\msf a^2\msf t^2}{8} + \log(y_{\msf t} + 1)
            \leq -\frac{\msf a^2\msf t^2}{8} + \log\msf t + C.
        \]
        By \eqref{eq:energy-iota}, applied to $\mu$ and to $\iota_*\sigma$, and by the minimality of $\mu$,
        \[
            \I^{\wt V}(\wt\mu) = \I^V(\mu) \leq \I^V(\iota_*\sigma)
            = \I^{\wt V}(\sigma) \leq -\frac{\msf a^2\msf t^2}{8} + \log\msf t + C.
        \]
        Since $\wt V \geq -\frac{\msf a^2\msf t^2}{8} + R_{\min}$ by \eqref{eq:Vtilde-square}, we get $\I(\wt\mu) = \I^{\wt V}(\wt\mu) - \int\wt V\,\dd\wt\mu \leq \log\msf t + C$, using that $\wt\mu$ is a probability measure. Therefore
        \begin{equation}\label{eq:ell-tilde-bound}
            \wt\ell = \I(\wt\mu) + \I^{\wt V}(\wt\mu)
            \leq -\frac{\msf a^2\msf t^2}{8} + 2\log\msf t + C.
        \end{equation}

        Next we bound the support of $\wt\mu$. Let $M_{\msf t} \defeq \max\{|y| : y \in \supp\wt\mu\}$, attained at $y_*$. Since $|y_* - s| \leq 2M_{\msf t}$ for $s \in \supp\wt\mu$, we have $2\,\U^{\wt\mu}(y_*) \geq -2\log(2M_{\msf t})$. As $F(\iota y_*) = 0$, \eqref{eq:F-iota}, \eqref{eq:ell-tilde-bound} and \eqref{eq:Vtilde-square} give
        \[
            \frac{1}{2\msf a^2}\left(y_* - y_{\msf t}\right)^2 - \frac{\msf a^2\msf t^2}{8} + R_{\min}
            \leq \wt V(y_*) = \wt\ell - 2\,\U^{\wt\mu}(y_*)
            \leq -\frac{\msf a^2\msf t^2}{8} + 2\log\msf t + 2\log(2M_{\msf t}) + C.
        \]
        Suppose that $M_{\msf t} > \msf a^2\msf t = 2y_{\msf t}$. Then $|y_* - y_{\msf t}| \geq M_{\msf t} - y_{\msf t} \geq \frac12M_{\msf t}$, and $\log\msf t \leq \log M_{\msf t} + C$. The display above therefore yields $M_{\msf t}^2 \leq C\left(\log M_{\msf t} + 1\right)$. This bounds $M_{\msf t}$ by a constant that depends only on $\msf a$, which is impossible if $\msf t$ is large. Hence, for $\msf t$ large, $\supp\wt\mu \subset [-\msf a^2\msf t,\msf a^2\msf t]$, and
        \begin{equation}\label{eq:U-tilde-lower}
            2\,\U^{\wt\mu}(y) \geq -2\log\left(|y| + \msf a^2\msf t\right),
            \qquad y \in \R.
        \end{equation}

        We now estimate $F(\iota y)$, first for $y < 0$ and then for $0 < y \leq K$. For $y < 0$ we have $(y - y_{\msf t})^2 \geq y^2 + y_{\msf t}^2$, so that \eqref{eq:Vtilde-square} gives $\wt V(y) \geq \frac{y^2}{2\msf a^2} + R_{\min}$. By \eqref{eq:F-iota}, \eqref{eq:U-tilde-lower} and \eqref{eq:ell-tilde-bound},
        \[
            F(\iota y) \geq \frac{\msf a^2\msf t^2}{8} - 2\log\msf t - C
            + \frac{y^2}{2\msf a^2} - 2\log\left(|y| + \msf a^2\msf t\right)
            \geq \frac{\msf a^2\msf t^2}{8} - 4\log\msf t - C.
        \]
        Here the last inequality uses $\log(|y| + \msf a^2\msf t) \leq \log(1 + |y|) + \log(1 + \msf a^2\msf t)$ and the fact that $\frac{y^2}{2\msf a^2} - 2\log(1 + |y|)$ is bounded from below. The right-hand side is positive for $\msf t$ large.

        For $y > 0$, \eqref{eq:Vtilde} and \eqref{eq:Vtilde-square} give $\wt V(y) = \frac{y^2}{2\msf a^2} + R(y) - \frac{\msf t\,y}{2} \geq R_{\min} - \frac{\msf t\,y}{2}$. Hence, for $0 < y \leq K$, by \eqref{eq:F-iota}, \eqref{eq:U-tilde-lower} and \eqref{eq:ell-tilde-bound},
        \[
            F(\iota y) \geq \frac{\msf a^2\msf t^2}{8} - \frac{K\msf t}{2}
            - 2\log\msf t - 2\log\left(K + \msf a^2\msf t\right) - C,
        \]
        which is positive for $\msf t$ large.

        Finally, $\iota$ maps $(-\infty,0)$ onto $(0,+\infty)$ and $(0,K]$ onto $\left(-\infty,-\frac1K\right]$, so the two cases above prove that $F > 0$ on $\left( -\infty, -\frac{1}{K} \right] \cup \left( 0, +\infty \right)$. Since $F = 0$ on $\supp\mu$, this gives $\supp\mu \subset \left(-\frac1K,0\right)$.
    \end{proof}

    \begin{proof}[Proof of Theorem~\ref{thm:intro-eq}]
        Let $K \defeq \msf a$, and let $\msf t$ be sufficiently large as in Lemma~\ref{lem:localization}. Then $F > 0$ on $\left(-\infty,-\frac1K\right]\cup(0,+\infty)$, $\supp\mu \subset \left(-\frac1K,0\right)$, and $\supp\wt\mu \subset (K,+\infty)$.

        Note that the function $y \mapsto F(\iota y)$ cannot vanish at two points $y_1 < y_2$ of $[K,+\infty)$ with $(y_1,y_2)\cap\supp\wt\mu = \emptyset$. Indeed, on $(y_1,y_2)$ the function $\U^{\wt\mu}$ is strictly convex, because $(\U^{\wt\mu})''(y) = \int\frac{\dd\wt\mu(s)}{(y - s)^2} > 0$ there, and $\wt V'' > 0$ by Remark~\ref{rem:Vtilde-convex}, because $K = \msf a$. By \eqref{eq:F-iota}, $y \mapsto F(\iota y)$ is therefore continuous on $[y_1,y_2]$ and strictly convex on $(y_1,y_2)$. If it vanished at both endpoints, it would be negative on $(y_1,y_2)$, contradicting $F \geq 0$.

        Now suppose that $\supp \mu$ is not an interval. Since $\iota$ maps $(K,+\infty)$ bijectively onto $\left(-\frac1K,0\right)$ and $\supp\wt\mu \subset (K,+\infty)$, the set $\supp\wt\mu$ is then not an interval either. Hence there are $y_1 < y_2$ in $\supp\wt\mu$ with $(y_1,y_2)\cap\supp\wt\mu = \emptyset$. Since $y_1,y_2 \in \supp\wt\mu$, we have $F(\iota y_1) = F(\iota y_2) = 0$, contradicting the observation above. Therefore $\supp\wt\mu = [\wt a_1,\wt b_1] \subset (K,+\infty)$, and $\supp\mu = [a_1,b_1] = \iota([\wt a_1,\wt b_1]) \subset \left(-\frac1K,0\right)$.

        To prove (R1), let $x_0 \in \R\setminus([a_1,b_1]\cup\{0\})$. If $x_0 \leq -\frac1K$ or $x_0 > 0$, then $F(x_0) > 0$ by Lemma~\ref{lem:localization}. Otherwise $y_0 \defeq \iota x_0 \in (K,+\infty)\setminus[\wt a_1,\wt b_1]$. Suppose for contradiction that $F(x_0) = 0$. Then $y_0$ and the endpoint of $[\wt a_1,\wt b_1]$ that is nearest to it are two zeros of $y \mapsto F(\iota y)$ in $(K,+\infty)$ with no point of $\supp\wt\mu$ between them, contradicting the observation above. Hence $F(x_0) > 0$.

        The conditions (R2) and (R3) follow from Lemma~\ref{lem:q-negative}, since $a_1 > -\frac1K = -\frac{1}{\msf a}$.
    \end{proof}

\section{Facts from potential theory}\label{app:prelim} \hfill

    In this appendix we collect the facts from potential theory that we use in Section~\ref{sec:asympt1}. They are taken from \cite{Ransford} and \cite{SaffTotik2024}, and we keep the notation of Section~\ref{sec:asympt1}. We write $x = \xi + \ii\upsilon$, with $\xi, \upsilon \in \R$, and $\dd A$ for the area measure on $\C$. For a measurable map $T$ and a measure $\sigma$, the push-forward $T_*\sigma$ is the measure $T_*\sigma(B) \defeq \sigma(T^{-1}(B))$, so that $\int\chi\,\dd(T_*\sigma) = \int\chi\circ T\,\dd\sigma$.

    Let $u$ be a locally integrable function on an open set $\Omega \subset \C$. Its Laplacian in the sense of distributions is the linear functional
    \[
        \chi \longmapsto \int u\,\Delta\chi\,\dd A,
    \]
    defined on smooth functions $\chi$ with compact support in $\Omega$. If $u$ is smooth, Green's formula shows that this is integration against the function $\Delta u$. If there is a signed measure $\sigma$ on $\Omega$ such that $\frac{1}{2\pi}\int u\,\Delta\chi\,\dd A = \int\chi\,\dd\sigma$ for every such $\chi$, we write $\frac{1}{2\pi}\Delta u = \sigma$. The basic example is $\frac{1}{2\pi}\Delta\log|x - s| = \delta_s$, which gives, for every finite positive measure $\sigma$ with compact support,
    \begin{equation}\label{eq:laplacian-potential}
        \frac{1}{2\pi}\Delta\U^\sigma = -\sigma.
    \end{equation}
    This is how the measures in Section~\ref{sec:asympt1} appear: we identify them as Laplacians of suitable functions. For instance, $\frac{1}{2\pi}\Delta\SW = \wt\mu - \mu$ by Lemma~\ref{lem:W-properties}(i). Finally, by Weyl's lemma, if $\Delta u = 0$ in the sense of distributions, then $u$ coincides almost everywhere with a harmonic function.

    Let $u$ be a subharmonic function on $\Omega$, with $u \not\equiv -\infty$ on each connected component of $\Omega$. Then $\Delta u$ is a positive Radon measure on $\Omega$, and we call $\sigma_u \defeq \frac{1}{2\pi}\Delta u$ the Riesz measure of $u$. By \eqref{eq:laplacian-potential}, $\sigma_u = \sigma$ for $u = -\U^\sigma$. Subharmonic functions satisfy the maximum principle and the sub-mean value property, and the maximum of two subharmonic functions is subharmonic. The Riesz decomposition theorem describes subharmonic functions locally: up to a harmonic function, they are negative potentials of their Riesz measures. We recall it in the form that we need, see \cite[Section~3.7]{Ransford}.

    \begin{theorem}[Riesz decomposition theorem]\label{thm:riesz}
        Let $u$ be a subharmonic function on an open set $\Omega \subset \C$, with $u \not\equiv -\infty$ on each connected component of $\Omega$, and let $D$ be an open disk with $\overline D \subset \Omega$. Then there is a harmonic function $h$ on $D$ such that
        \[
            u(x) = \int_D \log|x - s|\,\dd\sigma_u(s) + h(x), \qquad x \in D,
        \]
        that is, $u = -\U^{\restr{\sigma_u}{D}} + h$ on $D$.
    \end{theorem}

    \begin{corollary}\label{cor:riesz-no-atoms}
        Let $u$ be a continuous subharmonic function on an open set $\Omega \subset \C$. Then its Riesz measure $\sigma_u$ has no atoms.
    \end{corollary}

    \begin{proof}
        Let $x_0 \in \Omega$ and $c \defeq \sigma_u(\{x_0\})$. Fix $r \in (0,\frac12]$ with $\overline{D(x_0,r)} \subset \Omega$, put $D \defeq D(x_0,r)$, and let $h$ be as in Theorem~\ref{thm:riesz}. Write $\restr{\sigma_u}{D} = c\,\delta_{x_0} + \sigma'$, where $\sigma'$ is a positive measure.

        For $x \in D(x_0,\frac r2)$ and $s \in D$ we have $|x - s| < 2r \leq 1$, hence $\log|x - s| < 0$ and therefore $\int_D \log|x - s|\,\dd\sigma'(s) \leq 0$. Also $h$ is bounded on $D(x_0,\frac r2)$, because it is continuous on the closure of this disk. Theorem~\ref{thm:riesz} gives a constant $C$ such that
        \[
            u(x) \leq c\log|x - x_0| + C, \qquad x \in D(x_0,\tfrac r2)\setminus\{x_0\}.
        \]
        If $c > 0$, the right-hand side tends to $-\infty$ as $x \to x_0$, forcing $u(x) \to -\infty$ and contradicting the continuity of $u$ at $x_0$. Hence $c = 0$.
    \end{proof}

    Near each point of $\C\setminus\supp\sigma$, the potential $\U^\sigma$ is the real part of the analytic function $-\int\log(x - s)\,\dd\sigma(s)$, and differentiation under the integral sign gives
    \begin{equation}\label{eq:gradient-potential}
        \partial_\xi\U^\sigma - \ii\,\partial_\upsilon\U^\sigma = \CT^\sigma \quad \text{in } \C\setminus\supp\sigma.
    \end{equation}
    More generally, a real harmonic function $h$ is, near each point, the real part of an analytic function $\Psi$, and by the Cauchy--Riemann equations
    \begin{equation}\label{eq:gradient-harmonic}
        \Psi' = \partial_\xi h - \ii\,\partial_\upsilon h, \qquad \text{so that} \qquad |\Psi'| = |\nabla h|.
    \end{equation}
    In this way, the gradient of a potential, or of a difference of potentials such as $\SW$, is given by Cauchy transforms, which are analytic functions. If $\sigma$ is supported on $\R$, then $\CT^\sigma(\bar x) = \overline{\CT^\sigma(x)}$. Moreover, if $\sigma$ has a H\"older continuous density $\frac{\dd\sigma}{\dd x}$ on an open interval $I \subset \R$, then $\CT^\sigma$ extends continuously from the upper half plane to $I$, and its boundary values satisfy $\Im\CT^\sigma_+ = \pi\frac{\dd\sigma}{\dd x}$ on $I$, by the Sokhotski--Plemelj formula.

    Let $\Psi$ be a conformal map from an open set $U$ onto $\Psi(U)$. For smooth $\chi$ we have $\Delta(\chi\circ\Psi) = |\Psi'|^2\,(\Delta\chi)\circ\Psi$, and the Jacobian of $\Psi$ is $|\Psi'|^2$. Consequently, if $u$ is locally integrable on $\Psi(U)$ and $\frac{1}{2\pi}\Delta u$ is a measure, then
    \[
        \frac{1}{2\pi}\Delta(u\circ\Psi) = \big(\Psi^{-1}\big)_*\Big(\frac{1}{2\pi}\Delta u\Big) \quad \text{in } U.
    \]
    In other words, the measure $\frac{1}{2\pi}\Delta u$ does not depend on the conformal coordinate in which it is computed. We use this in the proof of Lemma~\ref{lem:laplacian-absW}, and a version for the two-to-one map $\msf z$ in the proof of Lemma~\ref{lem:Phi}.

    Let $\gamma$ be a smooth arc, and let $u$ be a continuous function near $\gamma$ which is harmonic on both sides of $\gamma$ and smooth up to $\gamma$ from each side. Applying Green's formula on each side of $\gamma$, one finds that, near $\gamma$,
    \[
        \frac{1}{2\pi}\Delta u = \frac{1}{2\pi}\,[\partial_n u]\,\dd s,
    \]
    which is a measure on $\gamma$. Here $\dd s$ is the arc length, and $[\partial_n u] = \partial_n u|_2 - \partial_n u|_1$ is the jump of the normal derivative of $u$ across $\gamma$, where $n$ is the unit normal pointing from side $1$ to side $2$ (the jump does not depend on the choice of $n$). For $u = |\SW|$ and $\gamma \subset \Gamma_{\SW}\setminus \mcal K$ the jump is $2|\nabla\SW|$. This is the origin of the term $2\nu$ in Lemma~\ref{lem:laplacian-absW}.

    The total mass of $\frac{1}{2\pi}\Delta u$ can be read off from the growth of $u$ at infinity. For instance, if $u$ is smooth and $u(\zeta) = \log|\zeta|$ for $|\zeta| \geq 1$, the divergence theorem gives, for $R \geq 1$,
    \[
        \frac{1}{2\pi}\int_{\{|\zeta| \leq R\}}\Delta u\,\dd A = \frac{1}{2\pi}\int_{\{|\zeta| = R\}}\partial_r u\,|\dd\zeta| = \frac{1}{2\pi}\cdot\frac1R\cdot 2\pi R = 1.
    \]
    Lemma~\ref{lem:Phi-mass} is a version of this for $\Phi$, which behaves like $\log|\zeta|$ at infinity, up to an additive constant, by Lemma~\ref{lem:Phi-continuity}. Since $\Phi$ is only continuous, in the proof we move the derivatives to cut-off functions instead of using the divergence theorem.

    Zeros of polynomials enter through the following identity. Let $p$ be a polynomial of degree $n \geq 1$, with leading coefficient $c$ and zeros $\zeta_1,\dots,\zeta_n$, repeated according to multiplicity, and let $\sigma_p \defeq \frac1n\sum_{k=1}^n\delta_{\zeta_k}$ be its normalized zero counting measure. Since $\log|p(\zeta)| = \log|c| + \sum_{k=1}^n\log|\zeta - \zeta_k|$, we have
    \begin{equation}\label{eq:log-poly}
        \frac1n\log|p| = \frac1n\log|c| - \U^{\sigma_p} \quad \text{on } \C.
    \end{equation}
    In particular, $\log|p|$ is subharmonic, and by \eqref{eq:laplacian-potential} its Riesz measure is $\frac{1}{2\pi}\Delta\log|p| = n\,\sigma_p = \sum_{k=1}^n\delta_{\zeta_k}$. Its total mass $n$ is the growth rate of $\log|p(\zeta)| = n\log|\zeta| + \log|c| + o(1)$ as $\zeta \to \infty$, in agreement with the previous paragraph, and the leading coefficient only contributes a constant, which disappears under the Laplacian. Thus, for a sequence of polynomials $p_N$, the limit of the zeros can be obtained from the limit of $\frac1N\log|p_N|$. In Section~\ref{sec:det} we use \eqref{eq:log-poly} directly, for $p_N = \det\bm P_N$. In Section~\ref{sec:proof-zeros-general} we take Laplacians: if $\frac1N\log|p_N|$ converges in $L^1_{\mathrm{loc}}(\C)$, then the measures $\frac{1}{2\pi}\Delta\big(\frac1N\log|p_N|\big) = \frac{\deg p_N}{N}\,\sigma_{p_N}$ converge in the sense of distributions.

    For finite positive measures $\sigma_N$ and $\sigma$ on $\C$, we say that $\sigma_N \to \sigma$ weakly if $\int h\,\dd\sigma_N \to \int h\,\dd\sigma$ for every bounded continuous function $h$ on $\C$, and vaguely if this holds for every continuous function $h$ with compact support. By Helly's selection theorem, every sequence with $\sup_N\sigma_N(\C) < \infty$ has a vaguely convergent subsequence. Taking $h \equiv 1$, weak convergence implies $\sigma_N(\C) \to \sigma(\C)$. A vague limit can lose mass, which escapes to infinity, but it cannot gain mass: taking the supremum over the functions $h$ with compact support and $0 \leq h \leq 1$ gives $\sigma(\C) \leq \liminf_N\sigma_N(\C)$. If no mass is lost, the convergence is weak. Indeed, assume that $\sigma_N \to \sigma$ vaguely and $\sigma_N(\C) \to \sigma(\C)$. Given $\epsilon > 0$, choose a continuous function $\chi$ with compact support, $0 \leq \chi \leq 1$ and $\int\chi\,\dd\sigma \geq \sigma(\C) - \epsilon$. Then $\limsup_N\int(1 - \chi)\,\dd\sigma_N \leq \epsilon$, and writing $h = h\chi + h(1 - \chi)$ for a bounded continuous function $h$ we obtain $\limsup_N\left|\int h\,\dd\sigma_N - \int h\,\dd\sigma\right| \leq 2\epsilon\|h\|_\infty$. Finally, assume that all the $\sigma_N$ are supported in a fixed compact set $L$. Then vague and weak convergence coincide, and if $\sigma_N \to \sigma$, then $\sigma$ is supported in $L$ and $\U^{\sigma_N}(\zeta) \to \U^\sigma(\zeta)$ for every $\zeta \notin L$, since $s \mapsto \log|\zeta - s|$ is continuous on a neighborhood of $L$.

\section*{Acknowledgments}

    The author thanks Arno Kuijlaars for pointing out his paper with Alan Groot \cite{Groot2021a}, whose central techniques are shared by the present work, Guilherme Silva for proposing the initial problem that grew into the general result of the paper, and Lucia Morey for many fruitful discussions.

\section*{Declarations}

    \subsection*{Funding} This work was partially funded by grants \#2022/12756-1 and \#2020/13183-0, S\~ao Paulo Research Foundation (FAPESP).

    \subsection*{Competing interests} The author has no relevant financial or non-financial interests to disclose.

    \subsection*{Data availability} Data sharing is not applicable to this article, as no datasets were generated or analyzed during the current study.

    \subsection*{Use of artificial intelligence} During the preparation of this work the author used AI agents as research and writing assistants. They were used to check computations and proofs, to propose more detailed versions of some arguments, to reorganize the presentation, to search the literature, and to write the code that draws Figure~\ref{fig:spectral-curve} and the code for the numerical computations behind Figure~\ref{fig:zeros-entries}. The author verified all the material obtained in this way, checked the references against the original sources, and takes full responsibility for the content of the paper.

\raggedbottom
\emergencystretch=2em
\printbibliography
\end{document}